\documentclass[11pt]{article}
\input amssymb.sty 
\input amssym.def 
\input amssym 
\input epsf 
\usepackage{amsfonts}
\input{cyracc.def}
\font\tencyr=wncyr10
\def\cyr{\tencyr\cyracc}

\newtheorem{theorem}{Theorem}[section]

\newtheorem{theorem-definition}[theorem]{Theorem-Definition}
\newtheorem{theorem-construction}[theorem]{Theorem-Construction}
\newtheorem{lemma-definition}[theorem]{Lemma--Definition}
\newtheorem{proposition-definition}[theorem]{Proposition--Definition}
\newtheorem{lemma}[theorem]{Lemma}
\newtheorem{proposition}[theorem]{Proposition}
\newtheorem{corollary}[theorem]{Corollary}
\newtheorem{conjecture}[theorem]{Conjecture}
\newtheorem{definition}[theorem]{Definition}

\begin{document}
\newcommand{\Z}{{\mathbb Z}}
\newcommand{\R}{{\mathbb R}}
\newcommand{\Q}{{\mathbb Q}}
\newcommand{\C}{{\mathbb C}}
\newcommand{\lms}{\longmapsto}
\newcommand{\lra}{\longrightarrow}
\newcommand{\hra}{\hookrightarrow}
\newcommand{\ra}{\rightarrow}
\newcommand{\sgn}{\rm sgn}
\begin{titlepage}
\title{The quantum dilogarithm and %unitary modular groupoids
representations of quantum cluster varieties}   
\author{V.V. Fock, A.B. Goncharov}
\end{titlepage}
\date{\it \small To David Kazhdan for his 60th birthday}
\maketitle
~~~~~~~~~~~~~~~~~~~~~~~~~~~~~~~~~~~~~~~~~~~~~~~~~~~~~~~~~~~~~~~~~~~~~~~
{\cyr ``Loshadp1 sostoit iz tre0h neravnyh polovin''. \vskip 1mm
~~~~~~~~~~~~~~~~~~~~~~~~~~~~~~~~~~~~~~~~~~~~~~~~~~~~~~~~~~~~~~~~~
~~~~~~~~~~~~~~~~~~~~~~~~~~~~~~~~
Konskii0 lechebnik, 

~~~~~~~~~~~~~~~~~~~~~~~~~~~~~~~~~~~~~~~~~~~~~~~~~~~~~~~~~~~~~~~~~~~~~
Sochinenie g. A. de Barra. Moskva 1868.}\footnote{``A horse consists of 
three unequal halves''.  cf.  
  A. de Barr, Horse doctor. Moscow 1868.}

\tableofcontents

\section{Introduction}

\subsection{An outline}

{\it Cluster varieties} %\cite{FG2} 
are relatives of 
cluster algebras \cite{FZI}.   
They are schemes over $\Z$ glued from algebraic tori ${\Bbb G}_m^N$. 
{\it Cluster modular groups}   act
 by automorphisms of cluster varieties. 
 They  include 
classical modular groups  of punctured surfaces, 
also known as the mapping class/Teichm\"uller groups. 
We employ certain extensions of these groups,  called  
{\it saturated cluster modular groups}.

\vskip 3mm
Our main result  is a construction of series of 
$\ast$-representations of quantum 
 cluster ${\cal X}$-varieties. 
By this we mean the following. 
A cluster ${\cal X}$-variety ${\cal X}$ 
is equipped with a natural Poisson structure. 
It admits a canonical non-commutative $q$-deformation. %\cite{FG2}. 
The latter gives rise to a non-commutative $\ast$-algebra ${\Bbb L}({\cal X}_q)$, 
which should be thought of as the algebra of 
regular functions on the $q$-deformation. 
Given a complex number $\hbar$, 
we consider its {\it Langlands modular double}, which is a $\ast$-algebra 
\begin{equation} \label{22.22,22}
{\bf L}_{\cal X}:= {\Bbb L}({\cal X}_q) \otimes_\Z 
{\Bbb L}({\cal X}^{\vee}_{q^{\vee}}), \qquad q = e^{i\pi\hbar}, 
\quad q^{\vee}:= e^{i\pi/\hbar}.
\end{equation}
 Here ${\cal X}^{\vee}$ is the 
Langlands dual cluster ${\cal X}$-variety\footnote{
When ${\cal X}$ is related to a split reductive group $G$, 
${\cal X}^{\vee}$ is  related
 to the Langlands dual group $G^{\vee}$ \cite{FG4}.} for  ${\cal X}$. 
The (saturated) cluster mapping class group $\widehat \Gamma$
 of ${\cal X}$ acts by automorphisms 
of the $\ast$-algebra ${\bf L}_{\cal X}$. 

Assume now that $\hbar$ is a positive real number, so that $|q|=1$. 
We define a Freschet linear space ${\cal S}_{\cal X}$, 
the {\it Schwartz space of ${\bf L}_{\cal X}$}, 
and  a 
$\ast$-representation of ${\bf L}_{\cal X}$  in ${\cal S}_{\cal X}$. 
Using the quantum dilogarithm function  
we construct a projective representation of 
the group $\widehat \Gamma$ in the Schwartz space ${\cal S}_{\cal X}$, intertwining 
the action of $\widehat \Gamma$ on ${\bf L}_{\cal X}$. The space ${\cal S}_{\cal X}$ 
is a dense subspace of a Hilbert space  ${V}_{\cal X}$.  
The group $\widehat \Gamma$ acts by unitary operators in ${V}_{\cal X}$. Summarizing: 
\begin{equation} \label{7.29.07.1}
\mbox{{\it A $\ast$-representation of quantum cluster ${\cal X}$-variety is a triple 
$(\widehat \Gamma, {\bf L}_{\cal X}, {\cal S}_{\cal X}\subset V_{\cal X}$}}). 
\end{equation} 
This representation 
is decomposed according  to the unitary characters of the center of 
 algebra ${\bf L}_{\cal X}$.

Equivalentely,  ${\cal S}_{\cal X}$ is a 
$\widehat \Gamma$-equivariant module over 
the $\ast$-algebra ${\bf L}_{\cal X}$\footnote{This just means that 
${\cal S}_{\cal X}$ is a (projective) module over the semidirect product 
of the group algebra $\Z[\widehat \Gamma]$ and  
${\bf L}_{\cal X}$.}, equipped with 
 a $\widehat \Gamma$-invariant unitary scalar product   
-- the Hilbert space 
$V_{\cal X}$ is the completion of ${\cal S}_{\cal X}$ for the scalar product. 

There is a  space of generalized functions ${\cal S}^*_{\cal X}$, defined 
as the topological dual  of ${\cal S}_{\cal X}$. 
\vskip 3mm
Cluster modular groups 
are automorphism groups of objects of 
 {\it cluster modular groupoids}.  
We construct  representation (\ref{7.29.07.1})  for 
groupoids rather then groups -- this generalisation in fact greatly 
simplifies the construction.

\vskip 3mm
The unitary projective representation of the group 
$\widehat \Gamma$ in the Hilbert space $V_{\cal X}$ 
can be viewed as a rather sophisticated 
analog of the Weil representation of the metaplectic cover of 
the group $Sp(2n,\Z)$. In both cases  representations are given 
by integral operators, intertwining 
representations of certain Heisenberg-type algebras. These 
{\it intertwiners}  in our case are of two types: the simpler ones 
are Weil's intertwiners; the kernels of the new ones are given 
by the quantum dilogarithm function.

\vskip 3mm
The  program of quantization of cluster ${\cal X}$-varieties, 
including a construction of 
intertwiners, was  
initiated in \cite{FG2} (the final version 
appeared in \cite{FG2II}).  However it lacked 
some ingredients, including a proof of crucial 
relations for the intertwiners, so representations of 
$\widehat \Gamma$ were not available. 

The main new features of the present approach are the following. 
We give another construction of the intertwiners, clarifying their structure.
We introduce the Schwartz space ${\cal S}_{\cal X}$. 
Since the algebra ${\bf L}_{\cal X}$ actually acts in  ${\cal S}_{\cal X}$,  
\footnote{${\bf L}_{\cal X}$ acts by unbounded operators, and thus 
does not act on $V_{\cal X}$ -- the space ${\cal S}_{\cal X}$ 
is its maximal domain of definition.}
the claim that 
the intertwiners indeed intertwine this action makes sense, and we prove it. 
We show that this implies the 
relations for the intertwiners. 
In the quasiclassical limit they give functional 
equations for the classical dilogarithm. The simplest instance of this program, 
quantization of the moduli space ${\cal M}_{0,5}$, was implemented in 
 \cite{Go2}. 
Not overshadowed by issues of algebraic nature, 
it may serve as an introduction to the proof of our main result. 
 
\vskip 3mm
 One  of applications of our   construction  
is {\it quantum higher Teichm\"uller theory}. 
Let $\widehat S$ be a surface $S$ with holes and a finite collection of 
marked points at the boundary, considered modulo isotopy. 
Let $G$ be a split reductive group. The pair $(G, \widehat S)$  
gives rise to 
a moduli space ${\cal X}_{G, \widehat S}$ \cite{FG1}, 
related to the moduli space of $G$-local systems on $S$. 
The modular group $\Gamma_S$ of $S$  acts on ${\cal X}_{G, \widehat S}$.

The moduli space ${\cal X}_{G, \widehat S}$, in the case when $G$ has trivial center,  
has a natural  
a cluster ${\cal X}$-variety structure  (see {\it loc. cit.} 
for the case of $G= PGL_n$)\footnote{A similar result for an arbitrary $G$ 
will appear in elsewhere.}. 
This means, in particular, that the moduli space ${\cal X}_{G,\widehat S}$ 
carries a $\Gamma_S$-equivariant collection of rational coordinate systems 
(atlas) 
with the following properties: 

(i) The natural Poisson structure $\{*, *\}$ on 
${\cal X}_{G,\widehat S}$ 
in each of the coordinate system $\{X_i\}$
has a standard quadratic form: $\{X_i, X_j\} = \varepsilon_{ij}X_iX_j$  
for certain 
$\varepsilon_{ij}\in \Z$ depending on the coordinate system. 

(ii) The transition transformations between different 
coordinate systems are  
compositions of certain standard transformations, {\it cluster mutations}. 

The set parametrising coordinate systems of the cluster atlas 
on ${\cal X}_{G,\widehat S}$ includes a subset provided by ideal 
triangulations of $\widehat S$ 
equipped with an additional data at each triangle.

Therefore our general construction 
provides 
a family of infinite dimensional unitary projective representations 
of the saturated cluster modular group $\widehat \Gamma_{G,\widehat S}$ 
related to the pair $(G, \widehat S)$. The group  $\widehat \Gamma_{G,\widehat S}$ 
includes, as a subquotient,  the classical modular group 
$\Gamma_{S}$ of $S$, 
but can be bigger if $G\not =PGL_2$. 

Here is the story for $G = GL_1$. 
The moduli space ${\cal L}_{GL_1, S}$ of $G$-local systems on $G$ 
is an algebraic torus ${\rm Hom}(H_1(S), GL_1)$. 
Assume that $S$ has no holes. Then 
$H_1(S)$ is a lattice with 
a symplectic structure, so ${\cal L}_{GL_1, S}$ 
is a symplectic algebraic torus. The group $\Gamma_S$ acts by its 
automorphisms through the quotient $Sp(H_1(S)) 
\stackrel{\sim}{=} Sp(2g, \Z)$. 
Our projective representation of $\Gamma_S$ in this case is nothing else but   
the projection $\Gamma_S\to Sp(2g, \Z)$ composed with 
the Weil representation of $Sp(2g, \Z)$. 
%If $S$ has punctures, there is a canonical projection 
%$\pi: {\cal L}_{GL_1, S} \to 
%{\Bbb G}_m^{\{\mbox{punctures of $S$}\}}$, and the functions 
%$\pi*f$ give the center of the Poisson algebra of functions on  
%${\cal L}_{GL_1, S}$. 
If $S$ has holes, we get a family of representations 
parametrised by unitary characters of copies of 
$\R^*_+$ assigned to the holes.

\vskip 3mm
We conjecture that  these representations for 
a given $G$ and different $S$ 
form an infinite dimensional modular functor -- understanding 
representations of quantum cluster varieties  as data  (\ref{7.29.07.1}) 
 allows to state this precisely -- see Section 6.2. 
The algebra (\ref{22.22,22}) in this case contains 
the Langlands modular double 
of the $q$-deformed algebra of functions on the moduli 
space of $G$-local systems on $S$. 

\vskip 3mm
Let $K_G$ be the maximal compact group of the complex group $G(\C)$. 
Quantization of the moduli space of $K_G$-local systems on $S$ 
given by Witten \cite{W1}, \cite{W2}, Hitchin \cite{H} 
and others leads to a construction of 
certain local systems  
on the moduli spaces ${\cal M}_{g,n}$, the  
local systems  of conformal blocks 
for the Wess-Zumino-Witten theory. One can view them as 
finite-dimensional projective representations of the modular group of $S$. 

Let us take now instead of the maximal compact subgroup $K_G$ 
a non-compact real form of the complex Lie group $G(\C)$, namely  the split real 
Lie group $G(\R)$. 
The quantum higher Teichm\"uller theory can be viewed 
as a quantization of the moduli space of $G(\R)$-local systems on 
$S$, and provides infinite-dimensional unitary 
projective  representations of the modular  group of $S$. 

The relationship between the two stories is similar 
to the relationship between finite dimensional 
representations of $G$ and principal series representations of the 
real Lie group $G(\R)$. %It survives the quantization. 

\vskip 3mm
To prove relations for the intertwiners we introduce and study 
a geometric object encapsulating 
their properties: the {\it symplectic double} of a cluster  
${\cal X}$-variety and its non-commutative $q$-deformation, 
the {\it quantum double}. The quantum double is determined by the 
holonomic system of difference relations for the intertwiners. 
The symplectic double has a structure of 
the symplectic groupoid related to the Poisson variety ${\cal X}$. 
The relation between the intertwiner and the symplectic double 
is similar to the one 
between the Fourier transform on a 
vector space $V$ and the cotangent bundle 
$T^*V$. 

\vskip 3mm
We define a  {\it canonical representation of the quantum double}. 
It is given by a data similar to (\ref{7.29.07.1}), and 
includes a canonical unitary projective representation of  $\widehat \Gamma$. 
This representations is realized in the Hilbert space $L^2({\cal A}^+)$, where 
${\cal A}^+$ is 
the manifold  of real positive points of the {\it cluster 
${\cal A}$-variety} assigned to ${\cal X}$. In the  
situation related to a pair $(G, \widehat S)$, 
there is a dual moduli space ${\cal A}_{G.\widehat S}$ introduced in \cite{FG1}, 
closely related to the moduli space of unipotent $G$-local systems on $S$.  
It has a cluster ${\cal A}$-variety structure, and 
the canonical representation 
is realized in the Hilbert space  $L^2({\cal A}^+_{G.\widehat S})$. 
\vskip 3mm

The symplectic/quantum double is an object of  independent interest. 
To support this, we show in \cite{FG5} how it 
appears in the higher Teichm\"uller theory. This 
is new even in the classical set-up: we  construct 
a collection of canonical coordinates on a 
certain modification and completion of 
the Teichm\"uller 
space for the double of surface $S$. 
A partial tropicalisation of the symplectic double 
delivers the cluster algebra with principal coefficients 
of Fomin and Zelevinsky \cite{FZIV}. 

\vskip 3mm
This paper 
provides a more transparent treatment of 
a part of \cite{FG2}:  Using the quantum dilogarithm we simplify and clarify 
construction of cluster ${\cal X}$-varieties and their 
$q$-deformations. 

\vskip 3mm
\paragraph{The quantum double and the  intertwiner}  
Here is an outline of this relationship. We use 
the standard cluster terminology, recalled in Section 2.1.  
\footnote{The reader might benefit from taking a brief look at Section 2.1.}
The saturated cluster modular 
groupoid $\widehat {\cal G}$ has a   combinatorial description. Its objects are 
{\it feeds}
 ${\bf i} = (I, \varepsilon_{ij}, d_i)$, where $I$ is a finite set, 
$d_i \in \Q_{>0}$, and 
$\varepsilon_{ij}$, $i,j\in I$,
 is an  integral valued matrix such that   $d_i\varepsilon_{ij}$ 
is skew-symmetric. 
The morphisms are {\it feed cluster transformations} 
modulo certain equivalence relations. 

Denote by ${\Bbb G}_m$ the multiplicative algebraic group  -- 
 one has ${\Bbb G}_m(K) = K^*$ 
for a field $K$. 
Each feed ${\bf i}$ gives rise to a split algebraic
 torus ${\cal A}_{\bf i} := {\Bbb G}_m^I$ equipped with 
canonical coordinates $A_i$, $i \in I$.  A feed cluster transformation 
${\bf c}: {\bf i} \to {\bf i'}$ gives rise to a positive birational map  
${\cal A}_{\bf i} \to {\cal A}_{\bf i'}$, which induces an isomorphism 
of their sets of real positive points. 

To define the canonical representation of the modular groupoid $\widehat {\cal G}$, 
we assign to a feed ${\bf i}$ 
a Hilbert space $L^2({\cal A}^+_{\bf i})$, where ${\cal A}^+_{\bf i} = \R^I$ 
is a vector space with coordinates $a_i := \log A_i$, $A_i \geq 0$. 
The {\it quantum torus algebra ${\bf D}^q_{\bf i} = {\bf D}_{\bf i}$} is represented  
as the algebra of all 
$ \hbar $-difference operators in this space ($q = e^{\pi i\hbar }$). 
It is generated by operators of 
multiplications 
by ${\rm exp}(a_i)$,  and shift by $2 \pi i \hbar$ 
of a given coordinate $a_j$.  
The intertwiner corresponding to a feed cluster transformation 
${\bf c}: {\bf i}\lra {\bf i'}$ 
is a unitary map of Hilbert spaces
\begin{equation} \label{8:40as}
{\bf K}_{{\bf c^o}}: 
L^2({\cal A}^+_{\bf i'}) \lra L^2({\cal A}^+_{\bf i}).
\end{equation}
Its Schwarz kernel 
satisfies a system of difference relations, consisting   
of $\hbar $-difference relations, and  
their ``Langlands dual'' counterparts. If $\hbar$ is irrational, 
it is  characterised uniquely up to a constant  
by  these relations. The 
${\hbar}$-difference relations alone do not characterize it. 

The $\hbar$-difference relations form an ideal 
in ${\bf D}^q_{\bf i'} \otimes {\bf D}^q_{\bf i}$. It 
gives rise to an isomorphism  of the fraction fields 
of algebras ${\bf D}^q_{\bf i'}$ and ${\bf D}^{q}_{\bf i}$. 
The quantum double is obtained by gluing the corresponding quantum tori 
along these maps. Setting $q=1$, we get the symplectic double. 

To show that the intertwiners provide a representation of a 
modular groupoid we need to prove that relations between cluster transformations  
give rise to the identity maps up to constants. 
One of the issues is that  in many important cases, 
e.g. in quantization of higher Teichm\"uller spaces,  we do not know 
an explicit form of all {\it trivial} (i.e. equivalent to the identity)  
cluster transformations.

\vskip 3mm
In the rest of the Introduction we give a more detailed account on 
the results of the paper and  relevant connections. 
Reversing the logic, we discuss the double first, 
and 
the quantum dilogarithm and the intertwiners after that.  

\vskip 3mm
\subsection{The symplectic double and its properties}
Cluster ${\cal X}$-varieties %were defined in \cite{FG2}. They 
are equipped with a 
positive atlas and a Poisson structure. They admit a canonical 
non-commutative $q$-deformation. 
Given a Poisson variety ${\cal X}$, 
we denote by   ${\cal X}^{\rm op}$ the same  variety equipped with 
the opposite Poisson structure $\{*,*\}^{\rm op}:= -\{*,*\}$. 
 
The symplectic double ${\cal D}$ of 
a cluster ${\cal X}$-variety  is a symplectic variety equipped with a 
positive atlas and  a Poisson map
\begin{equation} \label{tyu}
\pi: {\cal D} \lra {\cal X}\times {\cal X}^{\rm op}.
\end{equation}
There is an embedding $j:  {\cal X} \hra {\cal D}$. Its image is a Lagrangian subvariety.  
 The following diagram, where $\Delta_{\cal X}$ is the diagonal 
in ${\cal X}\times {\cal X}^{\rm op}$, is commutative: 
\begin{equation} \label{8:40}
\begin{array}{ccc}
{\cal X}&\stackrel{j}{\hra} &{\cal D} \\
\downarrow && \downarrow \pi\\
\Delta_{\cal X}&\stackrel{}{\hra} & {\cal X}\times {\cal X}^{\rm op}
\end{array}
\end{equation}
 There is a canonical involution $i$ of ${\cal D}$ interchanging the two projections in (\ref{tyu}). 
The pair   $({\cal X}, {\cal D})$ together with the maps $i, j, \pi$ 
has a structure of  the symplectic groupoid \cite{We} 
related to the Poisson space ${\cal X}$, where ${\cal D}$ 
serves as the space of morphisms, and ${\cal X}$ as the space of objects.

\vskip 3mm
A cluster ${\cal X}$-variety ${\cal X}$ comes 
together 
with a dual object ${\cal A}$, called a 
{\it cluster ${\cal A}$-variety}. 
The algebra of regular functions on 
${\cal A}$ coincides with the upper cluster algebra 
 from  \cite{BFZ}.  
There is a map $p: {\cal A} \to {\cal X}$; the triple  
$({\cal A}, {\cal X}, p)$ is called a {\it cluster ensemble}. 
There is a $2$-form $\Omega_{\cal A}$  on 
 ${\cal A}$; 
the map $p$ is the quotient along the null space of this form. 
The form  $\Omega_{\cal A}$ comes from a class in $K_2({\cal A})$. 

\vskip 3mm
The double ${\cal D}$ sits in a commutative diagram: 
\begin{equation} \label{11.5.06.1}
\begin{array}{ccc}
{\cal A}\times {\cal A} && \\
&\searrow \varphi &\\
p \times p\downarrow  && {\cal D}\\
&\swarrow \pi &\\
{\cal X}\times {\cal X}^{\rm op} & & 
\end{array}
\end{equation}
The  modular group $\widehat \Gamma$ acts  
by automorphisms of the triple 
$({\cal A}, {\cal X}, p)$, and of
 diagrams (\ref{8:40}) and (\ref{11.5.06.1}). %(\ref{11.5.06.q1}). 
\vskip 3mm
The very existence of the positive $\widehat \Gamma$-equivariant 
space ${\cal D}$ implies that there  is  a 
  symplectic space 
${\cal D}(\R_{>0})$ of the positive points of 
${\cal D}$, equipped with an action of  
$\widehat \Gamma$. It is isomorphic to $\R^{N}$.
We found another interesting $\widehat \Gamma$-equivariant real subspace 
in ${\cal D}(\C)$, the {\it unitary part} ${\cal D}^U$ of ${\cal D}(\C)$. 

\vskip 3mm
 {\it The quantum double}. 
We define in Section 3 a  non-commutative $\widehat \Gamma$-equivariant $q$-deformation 
${\cal D}_q$ of the double. There is a  map of quantum spaces 
$
\pi: {\cal D}_q \to {\cal X}_q\times {\cal X}^{\rm op}_q, 
$ and an involution $i: {\cal D}_q \to {\cal D}^{\rm op}_q$ interchanging the two components of $\pi$. 
The two real subspaces ${\cal D}(\R_{>0})$ and ${\cal D}^U$ of ${\cal D}(\C)$ 
match two $\ast$-algebra structures on ${\cal D}_q$.

\vskip3mm
{\it Duality conjectures for the double}. 
We conjecture that there exists a canonical basis in the space 
${\Bbb L}_+({\cal D})$ 
of regular positive functions on the positive space ${\cal D}$, as well 
as in its $q$-analog. 
Its vectors should be parametrised  by the 
set of integral tropical points 
of the Langlands dual double ${\cal D}^{\vee}$. 
This is closely related  to the  duality 
conjectures from Section 4 of \cite{FG2}. 
 
\vskip 3mm

{\it Relation with the work of Fomin and Zelevinsky}. 
After we found the symplectic double, we learned about 
a very interesting paper \cite{FZIV} were, among the other things, 
 ``cluster 
algebras with the principal coefficients'' were defined. The corresponding 
dynamics  formulas coincide with the partial tropicalization of the formulas 
which we use in the definition of the double ${\cal D}$. 

\vskip 3mm 
\subsection{Quantum dilogarithm,  quantum double,   and quantization}

The dilogarithm plays a crucial role in our  construction and  
understanding of 
cluster ${\cal X}$-varieties and their
doubles, both on the classical and quantum levels. 
The story goes as follows.

\vskip 3mm

\paragraph{The dilogarithm function and its two quantum versions}

(i) Recall the  {\it classical dilogarithm function} 
$
{\rm Li}_2(x):= -\int_0^x\log(1-t)d\log t.
$ 
We employ its version
$$
{\rm L}_2(x):= \int_0^x\log(1+t)\frac{dt}{t} = -{\rm Li}_2(-x).
$$

(ii) {\it The $q$-exponential as a ``quantum dilogarithm''}. 
Consider the following formal power series 
 $$
{\bf \Psi}^q(x):= \prod_{a=1}^{\infty}(1+q^{2a-1}x)^{-1} = 
\frac{1}{(1+qx)(1+q^3x)(1+q^5x)(1+q^7x)\ldots}.
$$
It is characterized by the 
difference relation
$
{\bf \Psi}^q(q^2x) = (1+qx){\bf \Psi}^q(x). 
$ 
If $|q|<1$ the power series ${\bf \Psi}^q(x)^{-1}$ converge, providing 
an analytic function in $x \in \C$. 
The dilogarithm appears in the quasiclassical limit: 
${\bf \Psi}^q(x)$ admits an asymptotic expansion when $q\to 1^-$:
\begin{equation} \label{sine}
{\bf \Psi}^q(x) \sim {\rm exp}\Bigl(\frac{{\rm L}_2(x)}{\log q^2}\Bigr). 
\end{equation}
Indeed, setting $x=e^t$, $q=e^{-s}$ we have 
$$
\log {\bf \Psi}^q(x) = -\sum_{k=1}^{\infty}\log (1+e^{t-(2k-1)s}) 
\sim_{s\to 0} -\int_0^{e^t}
\log(1+u)\frac{du}{u}/2s = {\rm L}_2(e^t)/-2s. 
$$

The function ${\bf \Psi}^q(x)$ 
appeared in XIX-th century under the names 
$q$-exponential, infinite Pochhammer symbol etc. 
The name $q$-exponential is 
justified by the following power series expansion: 
$$
{\bf \Psi}^q(x)= \sum_{n=0}^{\infty}\frac{x^n}
{(q-q^{-1}) (q^2-q^{-2}) \ldots (q^{n}-q^{-n})}.
$$
Sch\"utzenberger \cite{Sch} found a remarkable relation for 
the function ${\bf \Psi}^q(X)$ involving
 the $2$-dimensional quantum torus algebra. 
Its version  was rediscovered and interpreted as a quantum analog 
of the Abel's pentagon relation for the dilogarithm
 by Faddeev and Kashaev \cite{FK}. 

\vskip 3mm
(iii) {\it The quantum dilogarithm for $|q|=1$.} 
\footnote{One should rather call by quantum dilogarithms the logarithms of the functions $\Psi^q$ and $\Phi^\hbar$. The function $\Phi^\hbar$ is also known as 
the ``non-compact quantum dilogarithm''.}  Let 
${\rm sh}(u):= (e^u-e^{-u})/2$. Then 
$$
\Phi^\hbar(z) := {\rm exp}\Bigl(-\frac{1}{4}\int_{\Omega}\frac{e^{-ipz}}{ {\rm sh} (\pi p)
{\rm sh} (\pi \hbar p) } \frac{dp}{p} \Bigr), \qquad \hbar \in \R_{>0}, 
$$
where the contour $\Omega$ goes along the real axes 
from $- \infty$ to $\infty$ bypassing the origin 
from above.
Let 
$$
q=e^{\pi i \hbar}, \qquad q^{\vee}=e^{\pi i/\hbar}, \quad \hbar \in \R_{>0}.
$$ 
The function $\Phi^\hbar(z)$  is characterized by the {\it two}
 difference relations (Property {\bf B5} in Section 4.2):
$$
\Phi^h(z+2\pi i \hbar)= (1+qe^z)\Phi^\hbar(z), \qquad \Phi^\hbar(z+2\pi i )= (1+q^{\vee}e^{z/\hbar})\Phi^\hbar(z).
$$
It is related in several ways to the dilogarithm (Section 4.2), e.g. 
has an asymptotic expansion 
$$
\Phi^\hbar(z) \sim {\rm exp}\Bigl(\frac{{\rm L}_2(e^{z})}{2\pi i \hbar}\Bigr)\quad \mbox{as $\hbar \to 0$}. 
$$

The functions $\Phi^\hbar(z)$ and ${\bf \Psi}^q(e^{z})$ 
are close relatives: both deliver the classical dilogarithm in 
the quasiclassical limit, both are characterized by difference relations. 
When the Planck constant $\hbar$ is 
a complex number with ${\rm Im}~\hbar > 0$,  they are related by 
an  infinite product presentation 
$$
\Phi^\hbar(z) = \frac{{\Psi}^q(e^{z})}
{{\Psi}^{1/q^{\vee}}(e^{z/\hbar})}. 
$$
In this case $|q|<1$ as well as $|1/q^{\vee}|<1$. (The map $q \lms 1/q^{\vee}$ is 
a modular transformation). 
\vskip 3mm

The function $\Phi^\hbar(z)$ goes back to 
Barnes \cite{Ba}, and reappeared  in second half of XX-th century 
 in the works \cite{Sh}, 
\cite{Bax}, 
\cite{Fad}, \cite{K1}, \cite{CF}  and many others. 
The quantum pentagon relation for the function $\Phi^\hbar(z)$ was suggested in \cite{Fad} 
and proved in different ways/forms in \cite{Wo},\cite{FKV} and \cite{Go2}.

\vskip 3mm
\paragraph{The dilogarithm and cluster varieties} 
Every feed ${\bf i}$ provides a torus ${\cal X}_{\bf i}$, 
called the {\it feed ${\cal X}$-torus}. It
is equipped with canonical coordinates $X_i$, $i \in I$, 
providing an isomorphism 
${\cal X}_{\bf i} \to {\Bbb G}_m^{I}$. There is a Poisson structure defined by 
\begin{equation} \label{UIOT}
\{X_i, X_j\}= \widehat \varepsilon_{ij}X_iX_j, \quad  
\widehat \varepsilon_{ij}:= 
\varepsilon_{ij}d_j^{-1}. 
\end{equation}
A cluster ${\cal X}$-variety 
is obtained from the  feed ${\cal X}$-tori by the following 
gluing procedure. 
Every element $k \in I$ determines a feed ${\bf i'}$, called the 
{\it mutation of the feed ${\bf i}$ in the direction $k$}, and   
a Poisson birational isomorphism 
$\mu_k: {\cal X}_{\bf i} \to {\cal X}_{\bf i'}$. 
A cluster ${\cal X}$-variety is obtained by gluing the cluster 
feed tori ${\cal X}_{\bf i}$ parametrised by the 
feeds related to an initial feed 
by mutations, and taking the affine closure.

The cluster ${\cal A}$- and ${\cal D}$-varieties are constructed similarly 
from the feed ${\cal A}$- and ${\cal D}$-tori. 
In each case the combinatorial set-up is the same: 
the groupoid $\widehat {\cal G}$. 
However the gluing patterns are different. 
The obtained objects carry different structures: 
${\cal X}$-varieties are Poisson, ${\cal D}$-varieties are 
 symplectic with the symplectic structure coming from a class in $K_2$, and ${\cal A}$-varieties 
carry a class in $K_2$.

\vskip 3mm
One of the ideas advocated in this paper is that one should
 decompose the 
mutation birational isomorphism  $\mu_k$ 
into a composition of two maps, called 
``the automorphism part'' and ``the monomial part'' of the mutation. 
We introduce such
 decompositions for the ${\cal A}$-,  ${\cal D}$- and ${\cal X}$-varieties.  
They are compatible with the 
relating them maps 
(Theorem \ref{8.15.05.10}). 
The ``the automorphism part'' of the ${\cal X}$-mutation is a Poisson birational 
automorphism of the feed torus ${\cal X}_{\bf i}$ 
and the ``monomial part'' is a Poisson  isomorphism 
${\cal X}_{\bf i} \to {\cal X}_{\bf i'}$. 
The situation for the ${\cal D}$- and 
${\cal A}$-spaces is similar --  the $K_2$-classes are preserved. 
The ``monomial part'' is a very simple isomorphism. 
The  ``automorphism part'' is governed by the quantum dilogarithm.

\vskip 3mm

Here is how we use the classical and quantum dilogarithm functions 
in the gluing process. 
 The cluster coordinate $X_i$ on a feed 
 torus ${\cal X}_{\bf i}$, via the projection $\pi$, see (\ref{tyu}),  
lifts to 
the functions $X_i$ and $\widetilde X_i$ on the feed 
${\cal D}_{\bf i}$-torus. 
\vskip 3mm

(i) 
Recall that the Hamiltonian flow for a Hamiltonian function $H$ 
on a Poisson space %$(Y, \{*,*\})$ 
  with coordinates $y_i$  is given by
$
{dy_i(t)}/{dt} = \{H, y_i\}
$. 
The 
``automorphism part'' of the 
mutation map $\mu_k$ is the Hamiltonian flow 
for the time $1$ 
for the Hamiltonian functions given by:

\begin{itemize}

\item The dilogarithm function ${\rm L}_2(X_k)$ for the 
cluster ${\cal X}$-variety.

\item The difference of the dilogarithm functions ${\rm L}_2(X_k) - 
{\rm L}_2(\widetilde X_k)$ for the double  ${\cal D}$. 

\end {itemize}

(ii) Similarly, the ``quantum dilogarithm''  ${\bf \Psi}^q(X)$ is used to 
 describe the quantum mutation maps: 
The ``automorphism part'' of the quantum 
mutation in the direction $k$ 
is given by conjugation by  

\begin{itemize}

\item The 
``quantum dilogarithm''  
function ${\bf \Psi}^q(X_k)$ for the 
quantum cluster ${\cal X}$-variety.

\item The ratio of the ``quantum dilogarithms'' 
${\bf \Psi}^q(X_k)/{\bf \Psi}^q(\widetilde X_k)$ 
for the quantum double  ${\cal D}$. 

\end {itemize}

So although both  ${\rm L}_2(x)$ 
and ${\bf \Psi}^q(X)$ are transcendental functions, 
they provide birational transformations, 
``the automorphism parts'' of the mutation maps, and make all properties of the quantum 
${\cal X}$- and ${\cal D}$-spaces very transparent.

\vskip 3mm
\paragraph{Representations of quantum cluster ${\cal X}$-varieties} 
A lattice (= free abelian group) $\Lambda$ with 
a skew-symmetric form $\widehat \varepsilon: \Lambda \times \Lambda \to 
\frac{1}{N}\Z$ 
provides a {\it quantum torus algebra} ${\rm T}$. 
It is an associative $\ast$-algebra over $\Z[q^{1/N}, q^{-1/N}]$ with  
generators $e_{\lambda}, \lambda \in \Lambda$, and relations  
\begin{equation} \label{QTAL}
e_{\lambda} e_{\mu} = q^{-(\lambda, \mu)}e_{\lambda + \mu}, \qquad 
\ast e_{\lambda} = e_{\lambda}, \quad  \ast q = q^{-1}.
\end{equation}
A basis $\{e_i\}$, $i\in I$ in $\Lambda$ provides a set of generators 
$\{X^{\pm 1}_i\}$ 
of the algebra ${\rm T}$, satisfying the relations 
\begin{equation} \label{QTO}
q^{-\widehat \varepsilon_{ij}}X_iX_j = q^{-\widehat \varepsilon_{ji}}X_jX_i,\qquad 
\widehat \varepsilon_{ij}:= \widehat \varepsilon(e_i, e_j).
%\ast X_i = X_i, \quad \ast q = q^{-1}.
\end{equation}
Its quasiclassical limit delivers the Poisson structure (\ref{UIOT}). 
The feed algebras for the quantum ${\cal X}$- and ${\cal D}$-varieties  
are quantum tori algebras related to the corresponding Poisson structures. 

Given a quantum torus ${\rm T}$, we associate to it a Heisenberg algebra 
${\cal H}_{\rm T}$. It is a topological $\ast$-algebra over $\C$ generated by 
elements $x_i$ such that
$$
[x_p, x_q] = 2\pi i \hbar \widehat \varepsilon_{pq}; \quad \ast x_p = x_p, \quad 
q = e^{\pi i \hbar}. 
$$
Setting $X_p:= {\rm exp}(x_p)$ we get an embedding 
${\rm T}\hra {\cal H}_{\rm T}$.

\vskip 3mm
{\it Representations}. The Heisenberg $\ast$-algebra ${\cal H}_{\rm T}$ for $|q|=1$ 
has a family of 
$\ast$-representations
 by unbounded operators in a Hilbert space, 
parametrized by the central characters $\lambda$ 
of  ${\rm T}$. 
 \vskip 3mm
{\bf Example}. The simplest quantum torus algebra  is generated by $X,Y$ with the relation 
$qXY = q^{-1}YX$. The corresponding Heisenberg algebra is generated by 
$x,y$ with the relation $[x,y] = -2\pi i \hbar$. It has a $\ast$-representation in $L^2(\R)$:
$
 x \lms x, \quad y \lms 2\pi i \hbar \frac{\partial}{\partial x}.
$
\vskip 3mm

Denote by 
${\rm V}_{\lambda, {\bf i}}$ a $\ast$-representation 
of the quantum feed ${\cal X}$-torus algebra assigned to a feed ${\bf i}$ and a central character 
$\lambda$.
Given a mutation of feeds $ 
{\bf i} \to {\bf i'}$, we define a unitary operator 
\begin{equation} \label{12.15.06.2}
{\bf K}_{{\bf i'} \to {\bf i}}: 
{\rm V}_{\lambda, {\bf i'}} \lra {\rm V}_{\lambda, {\bf i}}, \qquad 
{\bf K}_{{\bf i'} \to {\bf i}}= {\bf K}^{\sharp} \circ {\bf K}'.
\end{equation}
Here ${\bf K}': {\rm V}_{\lambda, {\bf i'}} \lra {\rm V}_{\lambda, {\bf i}}$ 
is a very simple unitary operator corresponding to the ``monomial part''
 of the 
mutation map. The main hero is the automorphism ${\bf K}^{\sharp}: 
{\rm V}_{\lambda, {\bf i}} \lra {\rm V}_{\lambda, {\bf i}}$, given by 
\begin{equation} \label{12.q15.06.2}
{\bf K}^{\sharp} := {\bf \Phi}^\hbar(\widehat x_k),
\end{equation}
where  $\widehat x_k$ is a self-adjoint operator provided by the image $x_k$ 
in the representation ${\rm V}_{\lambda,  {\bf i}}$. 
The Schwartz kernels of the operators ${\bf K}'$ and ${\bf K}^{\sharp}$ 
are characterised by a Langlands dual pairs of systems of difference equations.
 They 
 describe the relevant  (``monomial''/''automorphism'') 
parts of the mutation map of the feed ${\cal X}$-torus.  
Compositions of elementary intertwiners (\ref{12.15.06.2}) give us the intertwiners 
${\bf K}_{\bf c^o}$ assigned to cluster transformations ${\bf c}$. 

We define the Schwartz space ${\cal S}_{\lambda,  {\bf i}} \subset {\rm V}_{\lambda,  {\bf i}}$ 
as the maximal domain of the 
algebra ${\bf L}$ in 
${\rm V}_{\lambda,  {\bf i}}$. Further, let $\widehat {\cal G}$ be 
the saturated cluster modular groupoid 
of ${\cal X}$. So $\widehat \Gamma$ is the automorphism group of its objects. 

\begin{theorem} \label{MT} The datum $
\Bigl(\widehat {\cal G}, {\bf L}, \{{\cal S}_{\lambda,  {\bf i}} 
\subset {\rm V}_{\lambda,  {\bf i}}\}, 
\{{\bf K}_{{\bf c}^o}\}\Bigr)
$ 
provides a $\ast$-representation of 
quantum cluster ${\cal X}$-variety. 
In particular operators (\ref{12.15.06.2}) provide a   unitary projective 
representation of the groupoid $\widehat {\cal G}$, and hence of the group $\widehat \Gamma$.
\end{theorem}

%More generally, we get a quantum cluster ${\cal X}$-variety given by the data 
%(\ref{7.29.07.1}) in the groupoid set-up. 

\vskip 3mm
{\bf Remark}. The mutation maps are described in (ii) 
via the formal power series ${\bf \Psi}^q(X)$. 
There is a similar description using the function ${\bf \Phi}^\hbar(x)$ and 
(\ref{12.q15.06.2}). 
The advantage of the description via ${\bf \Phi}^\hbar(x)$ 
is that it works in representations, while 
${\bf \Psi}^q(X)$ at $|q|=1$ makes sense  only as a power series, 
and thus does not act in a representation.

\vskip 3mm

\paragraph{Quantization of higher 
Teichm\"uller spaces}
Let $S$ be a surface with punctures. 
According to \cite{FG1},  Section 10, 
the pair of moduli spaces $({\cal A}_{G,S}, {\cal X}_{G,S})$ for $G=SL_m$  
has a cluster ensemble structure. 
Therefore thanks to Theorem \ref{MT} 
we arrive at quantization of  ${\cal X}_{G,S}$. It  includes  
a series of its 
$\ast$-representations, and thus a family of unitary projective 
representations of the saturated cluster modular group 
  $\widehat {\Gamma}_{G,S}$.  
They are parametrised by unitary characters of the group 
$
H(\R_{>0})^{\{\mbox{punctures of $S$}\}}
$, where $ H$ is the Cartan group of $G$. 
By the results of {\it loc. cit.} 
the classical modular group ${\Gamma}_{S}$ is a subgroup 
of the group ${\Gamma}_{G,S}$. The latter is a quotient of 
$\widehat {\Gamma}_{G,S}$. So ${\Gamma}_{S}$ is a subquotient of 
 $\widehat {\Gamma}_{G,S}$. 
We show that the obtained representations of the  
group $\widehat {\Gamma}_{G,S}$ descend to projective 
representations  of ${\Gamma}_{S}$. Here is 
our second main result:
\begin{theorem} \label{MT1} 
The cluster structure of 
$({\cal A}_{G,S}, {\cal X}_{G,S})$ gives rise to a series of  
$\ast$-representations of quantum moduli space ${\cal X}_{G,S}$. They provide 
a series of  infinite dimensional unitary projective 
representations of the classical mapping class group ${\Gamma}_{S}$. 
\end{theorem}

The $G=SL_2$ case was the 
subject of works \cite{K1} (which deals with a single representation with
 trivial central character) and 
 \cite{CF}. 
Unfortunately the argument presented in \cite{CF}  as a 
proof of the pentagon relation 
has a serious problem. As a result, 
the approach to quantization of Teichm\"uller spaces 
advocated in {\it loc. cit.} was put on hold. 
The proof of the pentagon 
relation for the quantum dilogarithm given  in \cite{Go2} serves as  a model 
for the proofs  in Section 5. 

We conjecture that the family 
of these representations for different surfaces $S$ form  an 
infinite dimensional modular functor. 
For $G=SL_2$ the proof is claimed by J. Teschner in \cite{T}.

\vskip 3mm

\paragraph{A canonical representation of the 
quantum double} 
It is defined similarly, 
by using the quantum double ${\cal D}$ instead of ${\cal X}$.  
Recall the quantum torus  algebra  ${\bf D}_{\bf i}$ for the quantum 
double assigned  
to a feed ${\bf i}$. Its 
Heisenberg algebra ${\cal H}_{{\bf D}_{\bf i}}$  
has a unique $\ast$-representation. It has a natural realization in 
  $L^2({\cal A}^+_{\bf i})$ where  ${\bf D}_{\bf i}$ 
acts as the algebra of $\hbar $-difference operators 
described in Section 1.1. 

Take the inverse images of $X_i \otimes 1$ and $1\otimes X_i$ under 
the quantum map $\pi^*$. 
Let $x_i, \widetilde x_i$ be their ``logarithms'', i.e. 
the  elements 
of the Heisenberg algebra ${\cal H}_{{\bf D}_{\bf i}}$. 
Given a mutation $\mu_k: {\bf i} \to {\bf i'}$, we define 
  a unitary operator 
\begin{equation} \label{12.15.06.3}
{\bf K}_{{\bf i'} \to {\bf i}}: 
L^2({\cal A}^+_{\bf i'}) \lra L^2({\cal A}^+_{\bf i}) 
\end{equation}
 as a composition 
${\bf K}_{{\bf i'} \to {\bf i}} = {\bf K}^{\sharp}\circ {\bf K}'$, where ${\bf K}'$ 
is the unitary operator induced on functions by a linear map of spaces 
$ {\cal A}^+_{\bf i} \to {\cal A}^+_{\bf i'}$,  the ``monomial part''
 of the ${\cal D}$-mutation map, and 
\begin{equation} \label{12.15.06.3q}
\quad 
{\bf K}^{\sharp}:= {\bf \Phi}^\hbar(\widehat x_k){{\bf \Phi}^\hbar(\widehat 
{{\widetilde x}_k})}^{-1}.
\end{equation}
The operator (\ref{12.15.06.3}) 
coincides with the one defined in \cite{FG2II} in a different way. 
 The Schwartz kernels of the operators ${\bf K}'$ and ${\bf K}^{\sharp}$ 
are characterised by Langlands dual pairs of systems of difference equations.
 They  
 describes the relevant  (``monomial''/''automorphism'') parts of 
 the coordinate transformations 
for the mutation  ${\bf i} \to {\bf i'}$,   
employed in the definition of the quantum double. 
In particular the 
 symplectic double describes 
the quasiclassical limit of intertwiners (\ref{12.15.06.3}). 
Compositions of elementary intertwiners (\ref{12.15.06.3}) give us the intertwiners 
${\bf K}_{\bf c^o}$, see (\ref{8:40as}). 

\vskip 3mm
The quantum double ${\cal D}_q$ is a functor: 
we assign to a feed  ${\bf i}$  
the quantum torus algebras ${\bf D}_{\bf i}$, 
and to a feed cluster transformation ${\bf c}: 
{\bf i} \lra {\bf i'}$ 
an isomorphism $\gamma_{\bf c^o}$ 
of the fraction fields of algebras ${\bf D}_{\bf i'}$ 
and  ${\bf D}_{\bf i}$. The {\it algebra of regular 
functions ${\Bbb L}({\cal D}_q)$ on  ${\cal D}_{q}$} consists of collections of elements 
  $F_{\bf i}\in {\bf D}_{\bf i}$  identified by these maps: 
 $\gamma_{\bf c^o}(F_{\bf i'}) = F_{\bf i}$. 
The crucial role 
plays the 
algebra ${\bf L}:= {\Bbb L}({\cal D}_q)\otimes 
{\Bbb L}({\cal D}^{\vee}_{q^{\vee}})$ of regular functions on 
the {\it quantum 
Langlands modular double ${\cal D}_{q} \times {\cal D}^{\vee}_{q^{\vee}}$}.  
For each feed ${\bf i}$ it is identified with  a subalgebra ${\bf L}_{\bf i} \subset 
{\bf D}_{{\bf i}, q}\times {\bf D}_{{\bf i^{\vee}}, q^{\vee}}$. 
The map $\gamma_{\bf c^o}$ induces an isomorphism ${\bf L}_{\bf i'}
\stackrel{\sim}{\lra} {\bf L}_{\bf i}$. The algebra ${\bf L}_{\bf i}$ acts by unbounded operators in 
the 
 {Schwartz space} $S_{\bf i}\subset L^2({\cal A}_{\bf i}^+)$,  
defined as the maximal common domain of operators from ${\bf L}_{\bf i}$. 
We prove that 
intertwiners restrict 
to isomorphisms  between the Schwarz spaces. This allows to introduce  
distribution spaces respected
by  the intertwiners. Summarizing, we get

\begin{theorem} \label{KMa} The datum
 $
\Bigl(\widehat {\cal G}, {\bf L}_{\bf i}, {\cal S}_{\bf i} \subset
 L^2({\cal A}_{\bf i}^+), {\bf K}_{{\bf c^o}}\Bigr) 
$ 
provides a {\rm canonical $\ast$-representation of the quantum cluster double}. 
This means that: 

\begin{itemize}

\item 
We have a   unitary projective representation of the groupoid $\widehat {\cal G}$ 
given by the Hilbert spaces $L^2({\cal A}_{\bf i}^+)$ and the unitary maps 
${\bf K}_{{\bf c^o}}$  between them;

\item The maps ${\bf K}_{{\bf c^o}}$ 
preserve the Schwartz spaces ${\cal S}_{\bf i}$ and intertwine on them the action of the groupoid 
$\widehat {\cal G}$ on the algebras ${\bf L}_{\bf i}$: for any 
$s \in {\cal S}_{\bf i'}$ and $A \in {\bf L}_{\bf i'}$ one has 
$$
{\bf K}_{{\bf c^o}} A (s) = \gamma_{\bf c^o}(A) {\bf K}_{{\bf c^o}} (s).
$$ 
\end{itemize}
\end{theorem} 

So for every cluster transformation ${\bf c}: 
{\bf i} \lra {\bf i'}$ there are commutative diagrams
$$
\begin{array}{ccccc}
{\bf L}_{\bf i}&\mbox{acts on} & {\cal S}_{\bf i}&\hookrightarrow &L^2({\cal A}_{\bf i}^+)\\
&&&&\\
\gamma_{\bf c^o}\downarrow &&\downarrow {\bf K}_{{\bf c^o}}&&\downarrow {\bf K}_{{\bf c^o}}\\
&&&&\\
{\bf L}_{\bf i'}&\mbox{acts on} &{\cal S}_{\bf i'}&\hookrightarrow &L^2({\cal A}_{\bf i'}^+)
\end{array} ~~~~~~~~~~~
\begin{array}{ccc}
{\bf L}_{\bf i}&\mbox{acts on} & {\cal S}^*_{\bf i}\\
&&\\
\gamma_{\bf c^o}\downarrow &&\downarrow {\bf K}_{{\bf c^o}}\\
&&\\
{\bf L}_{\bf i'}&\mbox{acts on} &{\cal S}^*_{\bf i'}
\end{array}
$$

The quantum map $\pi$ provides a subalgebra 
$\pi_q^*({\bf L}_{\cal X}\otimes {\bf L}^{\rm op}_{{\cal X}}) \subset {\bf L}$.  
It is isomorphic to the quotient of ${\bf L}_{\cal X}\otimes {\bf L}^{\rm op}_{{\cal X}}$ 
by a central subalgebra isomorphic to  ${\rm Center}({\bf L}_{\cal X})$. 
Restricting the canonical representation to this subalgebra 
and decomposing it  according to its central characters, identified with the characters of 
${\rm Center}({\bf L}_{\cal X})$,  we get, for every feed  ${\bf i}$, 
 a decomposition into an integral of  Hilbert spaces: 
$$
L^2({\cal A}^+_{\bf i}) = \int_{\lambda} {\rm End}(V_{\lambda,{\bf i} }) d\lambda.
$$
These decompositions 
are respected by the intertwiners. 
Thus the canonical representation is 
an integral 
of the endomorphisms of the principal series 
representations of the quantum cluster variety ${\cal X}_q$. 
So {\it the canonical representation is a regular representation 
for the $\ast$-algebra of regular functions on ${\cal X}_q$}.  
Using this one can easily deduce Theorem \ref{MT} from  Theorem \ref{KMa}. 

\vskip 3mm
\subsection{The symplectic double and 
higher Teichm\"uller theory \cite{FG5}} 
Let $S$ be an  oriented 
 topological surface with boundary, 
and $G$ a split semi-simple group over $\Q$ with trivial center. 
Cluster ensembles  provide a framework and 
tool for study of the dual pair 
$({\cal A}_{G,S}, {\cal X}_{G,S})$ of moduli spaces related to $S$ and $G$  
\cite{FG1}.  In 
\cite{FG5}
we show that the same happens with the  double. 

\vskip 3mm
{\it A moduli space ${\cal D}_{G,S}$}. Let $S_{\cal D}$ be the double of $S$. It is a 
topological surface obtained by gluing the surface $S$ 
with its ``mirror'', given by the same surface 
with the opposite orientation, 
along the corresponding boundary components. We introduce a  moduli space ${\cal D}_{G,S}$. 
It is a rather non-trivial
 relative of the moduli space of $G$-local systems on the double $S_{\cal D}$. 
In particular it contains a divisor whose points do not correspond to any kind of 
local systems 
on $S_{\cal D}$. Its irreducible components match 
components of the boundary of $S$.

The modular group $\Gamma_S$ of $S$ 
acts on ${\cal D}_{G,S}$. 
The moduli space 
${\cal D}_{G,S}$ has a natural $\Gamma_S$-invariant symplectic structure and a class in 
$K_2$. There are commutative diagrams similar to (\ref{8:40}) and (\ref{11.5.06.1}). 
We prove that, unlike the 
classical moduli space of $G$-local systems on $S_{\cal D}$,  
the moduli space ${\cal D}_{G,S}$ is rational, and 
equip  it  with a positive 
$\Gamma_S$-equivariant atlas. 
It is described as the symplectic double 
of the cluster ${\cal X}$-variety structure on the moduli 
space ${\cal X}_{G,S}$. 
 This was not known even for $G=PGL_2$.

%We introduce a notion of an integral ${\cal D}$-lamination on $S_{\cal D}$, prove that the set of 
% integral ${\cal D}$-laminations is identified with the set of integral tropical points 
%of the moduli space ${\cal D}_{PGL_2,S}$, and using this define  in Section 7.4 a canonical 
%basis in the space of regular function on ${\cal D}_{PGL_2,S}$. 

\vskip 3mm
{\it The symplectic double and quasifuchsian representations}. 
The space  ${\cal D}_{G,S}(\R_{>0})$ of positive real points 
of ${\cal D}_{G,S}$ is  a 
  symplectic space isomorphic to $\R^{-2\chi(S){\rm dim}G}$. 
We believe that it is closely related to the space of 
{\it framed quasifuchsian} representations $\pi_1(S) \to G(\C)$ modulo 
conjugation. For  $G = PGL_2$ and closed $S$ this 
boils down to the Bers double uniformization theorem.

\vskip 3mm
{\it Relation with the work of Bondal}. We will show elsewhere that 
the symplectic groupoid introduced in \cite{B} is the symplectic double 
in our sense of a certain twisted moduli space ${\cal X}_{G, \widehat S}$.

\vskip 3mm
{\bf Acknowledgments}. We are grateful to Alexei Bondal, Sergey Fomin, 
and  Andrey Zelevinski, and especially to Joseph Bernstein 
for useful conversations. 

V.F. was supported by the grants CRDF 2622; 2660. 
A.G. was 
supported by the  NSF grant  DMS-0400449 and  DMS-0653721. 
The first draft of the paper was written in 
August of 2005, and the last one in 2007 
when A.G. enjoyed the 
hospitality of IHES.  He is grateful to IHES for the hospitality and support. 

We are very 
grateful to the referee for the  
useful comments and suggestions incorporated in the paper.

\section{The symplectic double of a cluster ${\cal X}$-variety}

\subsection{Cluster ${\cal X}$- and ${\cal A}$-varieties} 
For the convinience of the reader we recall below 
the language of cluster ${\cal X}$- and ${\cal A}$-varieties 
in the form it was introduced in \cite{FG2}, 
borrowing from Section 2 
of {\it loc. cit.}, presenting a simplified version, without
 frozen variables. 

We reprove all the results mentioned below 
in Theorems \ref{8.15.05.10} and 
\ref{8.15.05.10gg} and  Lemma \ref{11.18.06.12trtr}. 
Doing the quantum 
case first and using the decomposition of mutations we simplify a 
lot the original proofs given in  {\it loc. cit.}.

\vskip 3mm
A {\em feed} \footnote{A feed is a combinatorial data  obtained by 
excluding cluster coordinates from a seed of \cite{FZI}.} 
${\mathbf i}$ is a triple $(I, \varepsilon, d)$, 
where $I$ is a finite set, 
$\varepsilon$ is a matrix $\varepsilon_{ij}$, where $i,j \in I$, 
such that 
$\varepsilon_{ij} \in {\mathbb Z}$, and $d = \{d_i\}$, where $i \in I$, 
are positive integers, such that the matrix 
$\widehat{\varepsilon}_{ij}=\varepsilon_{ij}d^{-1}_j$ is skew-symmetric.

For a feed ${\mathbf i}$ we associate a torus 
${\mathcal X}_{\mathbf i} = ({\mathbb G}_m)^I$ 
with a Poisson structure given by 
\begin{equation} \label{f1}
\{X_i,X_j\}=\widehat{\varepsilon}_{ij}X_iX_j,
\end{equation} 
where $\{X_i| i\in I\}$ are the 
standard coordinates on the factors. It is called the {\em cluster ${\cal X}$-torus}. 

Let ${\mathbf i}=(I, \varepsilon, d)$ and ${\mathbf i}'=(I', \varepsilon', d')$  be two feeds, and $k\in
I$. A {\em mutation in the direction $k\in I$} is an isomorphism $\mu_k: I\rightarrow I'$ satisfying the following conditions: $d'_{\mu_k(i)}=d_i$, and 
\begin{equation} \label{f2}
\varepsilon'_{\mu_k(i)\mu_k(j)}=
\left\{ \begin{array}{lll} -\varepsilon_{ij} & \mbox{ if }  i=k \mbox{ or } j=k,\\
\varepsilon_{ij}& \mbox{ if } \varepsilon_{ik}\varepsilon_{kj} \leq  0,\\
\varepsilon_{ij} + |\varepsilon_{ik}|\varepsilon_{kj} & \mbox{ if } \varepsilon_{ik}\varepsilon_{kj}> 0.
\end{array}\right.
\end{equation}

An {\em automorphism $\sigma$} 
of a feed ${\mathbf i}=(I, \varepsilon, d)$ is 
an  automorphism  of the set $I$ preserving 
the matrix $\varepsilon$ and the numbers $d_i$. 
Automorphisms and mutations induce rational maps between the corresponding  
cluster ${\cal X}$-tori, denoted by the same symbols $\mu_k$ and $\sigma$ and 
acting on the coordinate functions by the formulae 
$\sigma^*X_{\sigma(i)}=X_i$ 
and
\begin{equation} \label{f3}
\mu_{k}^*X_{\mu_k(i)} = \left\{\begin{array}{lll} X_k^{-1}& \mbox{ if } & i=k, \\
    X_i(1+X_k^{-\sgn(\varepsilon_{ik})})^{-\varepsilon_{ik}} & \mbox{ if } &  i\neq k. \\
\end{array} \right.
\end{equation} 

A {\em cluster transformation} between two feeds (and between two cluster
${\cal X}$-tori) is a composition of automorphisms and mutations. 
Two feeds are called {\em equivalent} if they are related by a cluster transformation. The equivalence class of a feed ${\mathbf i}$ is denoted by $|{\mathbf i}|$.

\vskip 3mm
 Recall that  the affine closure of a scheme $Y$ is the Spectrum of the ring of regular functions on $Y$. 
For instance the affine closure of $\C^2 -\{0\}$ is $\C^2$. 
A {\em cluster ${\cal X}$-variety} is a scheme over $\Z$ obtained by gluing the
 feed ${\cal X}$-tori for the feeds equivalent to a given 
feed ${\mathbf i}$  
via the above birational isomorphisms, and taking 
the affine closure of the obtained scheme. 
It is denoted by ${\cal X}_{|{\mathbf i}|}$, or simply by ${\cal X}$ if the 
equivalence class of feeds is apparent. 
   Every feed 
  provides our cluster ${\cal X}$-variety with a rational  coordinate system. The corresponding rational functions are called {\em cluster coordinates}.

\vskip 3mm
Cluster transformations preserve the Poisson structure. 
In  particular a cluster ${\cal X}$-variety
 has a canonical 
Poisson structure. 
 The mutation formulas (\ref{f2}) 
are recovered from the ones (\ref{f3}) 
and the condition that the latter preserve the Poisson structure.

\vskip 3mm
 Recall now the definition of the 
 {\em cluster ${\cal A}$-variety}. Given a feed ${\bf i}$, we define 
a torus ${\mathcal A}_{\mathbf i} = ({\mathbb G}_m)^I$ 
with the standard  
standard coordinates $\{A_i| i\in I\}$ on the factors. We call it the 
{\em feed ${\cal A}$-torus}.

Automorphisms and mutations give rise to  birational 
maps between the corresponding  
feed ${\cal A}$-tori, which are  given by 
$\sigma^*A_{\sigma(i)}=A_i$ 
and
\begin{equation} \label{f3cc}
\mu_{k(i)}^*A_{\mu_k(i)} = \left\{\begin{array}{lll} A_i
 & \mbox{ if } & i\not =k, \\
    A_k^{-1}\left(\prod_{i|\varepsilon_{ki}>0}A_i^{\varepsilon_{ki}} + 
\prod_{i|\varepsilon_{ki}<0}A_i^{-\varepsilon_{ki}}\right) & \mbox{ if } &  i= k. \\
\end{array} \right.
\end{equation} 
The cluster ${\cal A}$-variety corresponding to 
a feed ${\bf i}$ is a scheme over $\Z$ obtained by gluing 
 all feed ${\cal A}$-tori for the feeds equivalent to a given 
feed ${\mathbf i}$
  using the above birational isomorphisms, and taking 
the affine closure. 
It is denoted by ${\cal A}_{|{\mathbf i}|}$, or simply by ${\cal A}$. 

\vskip 3mm
There is a $2$-form $\Omega$ on the cluster ${\cal A}$-variety 
(\cite{GSV2}, \cite{FG2}), given 
in every cluster coordinate system by 
\begin{equation} \label{f1c}
\Omega_{\mathbf i} =\sum_{i,j\in I} 
\widetilde {\varepsilon}_{ij}d\log A_i\wedge d\log A_j,  
\quad \widetilde {\varepsilon}_{ij} = d_i {\varepsilon}_{ij}.
\end{equation}

There is a class in ${\bf W}_{\cal A} \in K_2({\cal A})$ 
given in every cluster coordinate system by 
$$
 \sum_{i,j\in I} 
\widetilde {\varepsilon}_{ij} \{A_i, A_j\} \in K_2(\Q({\cal A}_{\bf i})).
$$
It is a $K_2$-avatar of the $2$-form $\Omega_{\cal A}$, in the sense that 
$\Omega_{\cal A} = d\log({\bf W}_{\cal A})$, see Section 2.3 below.

\vskip 3mm
There is a map 
$
p: {\cal A} \lra {\cal X}, 
$
given in every cluster coordinate system by $p^*X_k = \prod_{I\in I}A_i^{\varepsilon_{ki}}$. 
It is the quotient of the space ${\cal A}$ 
along the null-foliation of the $2$-form $\Omega_{\cal A}$.

\vskip 3mm
{\it The cluster modular groupoid}. 
The inverse of a mutation is a mutation: $\mu_k\mu_k=id$. 
Feed cluster transformations inducing  
isomorphisms of the feed ${\cal A}$-tori as well as  the feed ${\cal X}$-tori   
are called {\it trivial feed cluster transformations}. 
The {\it cluster modular groupoid} ${\cal G} = {\cal G}_{|{\bf i}|}$
\footnote{Below we skip the subscripts $|{\bf i}|$ whenever possible.}
is a category whose objects are feeds equivalent to a given feed ${\bf i}$, and 
$Hom({\bf i}, {\bf i'})$ is the set of all 
feed cluster transformations from ${\bf i}$ to  ${\bf i'}$ modulo the trivial ones. 
In particular, given  a feed ${\bf i}$, 
the {\it cluster modular group $\Gamma_{\bf i}$} is 
the automorphism group of the object ${\bf i}$ of ${\cal G}$. 
By the very definition, it 
 acts by automorphisms 
of the cluster ${\cal A}$-variety. 
It preserves the class in $K_2$. 
\vskip 3mm

{\it The $(h+2)$-gon relations}. They are crucial
 examples of trivial cluster transformations. 
Denote by $\sigma_{ij}$ the map of feeds induced by  
interchanging $i$ and $j$. Let $h=2,3,4,6$  when $p = 0, 1, 2, 3$ respectively. Then if $\varepsilon_{ij}=-p\varepsilon_{ji}=-p$, 
then $(\sigma_{ij}\circ \mu_i)^{h+2} ={\rm Id}$ on feeds, and 
\begin{equation} \label{K10}
(\sigma_{ij}\circ \mu_i)^{h+2} = \mbox{a trivial cluster transformation}. 
\end{equation}
Relations (\ref{K10}) are affiliated with the rank two 
Dynkin diagrams, i.e. $A_1 \times A_1, A_2, B_2, G_2$. 
The number $h=2,3,4,6$ is the Coxeter 
number of the diagram. 
We do not know any other general procedure 
to generate trivial cluster transformations.

{\it The saturated 
coordinate groupoid} $\widehat {\cal G}$. 
Its objects are the isomorphism classes 
of feeds equivalent to a given one. 
The morphisms admit an explicit description as compositions of  
the cluster transformations (\ref{K10}). 
Its fundamental 
group $\widehat \Gamma$ is 
the {\it saturated  cluster modular group}. 
There is a canonical surjective map $\widehat \Gamma
\to \Gamma$. 
\vskip 3mm
{\it The quantum space ${\cal X}_q$}. It is a canonical non-commutative 
$q$-deformation of the cluster ${\cal X}$-variety. 
We start from the {\it feed quantum torus algebra} ${\rm T}^q_{\bf i}$, defined as 
an associative $\ast$-algebra with generators 
$X^{\pm 1}_i$, $i\in I$ and $q^{\pm 1}$ and relations
$$
q^{-\widehat \varepsilon_{ij}}X_iX_j = q^{-\widehat \varepsilon_{ji}}X_jX_i,\qquad 
\ast X_i = X_i, \quad \ast q = q^{-1}.
$$
Let ${\rm QTor}^*$ be a category whose objects are 
quantum torus algebras and morphisms are subtraction free 
$\ast$-homomorphisms of their fraction fields over $\Q$. 
The quantum space ${\cal X}_q$ is understood 
as a contravariant functor 
$$
\eta^q: 
\mbox{The  modular groupoid $\widehat {\cal G}$} \lra {\rm QTor}^*.
$$
It assigns to 
a feed ${\bf i}$ the quantum torus $\ast$-algebra ${\rm T}_{{\bf i}}$, and 
to a mutation   ${\bf i} \lra {\bf i'}$ a map of the fraction fields 
${\rm Frac}({\rm T}_{{\bf i'}}) \lra {\rm Frac}({\rm T}_{{\bf i}})$, given by 
a $q$-deformation of formulas (\ref{f3}) (
(\cite{FG2}, Section 3)).\footnote{A more transparent and 
less computational
 definition is given in Section 3.}  
The group $\widehat \Gamma_{\bf i}$ 
acts by automorphisms of the quantum 
space ${\cal X}_q$. 

One can 
understood cluster ${\cal A}$- and ${\cal X}$-varieties as 
similar functors: this is equivalent to their  definition  
as schemes. For the quantum spaces this is the only way to 
do it.  

\vskip 3mm
There are tori $H_{\cal X}$ and  $H_{\cal A}$ related to a cluster ensemble. 
They play the following role. 
There is 
a surjective projection $\theta: {\cal X} \lra H_{\cal X}$, as well as its quantum version 
$\theta_q: {\cal X}_q \lra H_{\cal X}$.
The torus $H_{\cal A}$ acts on ${\cal A}$, and 
the map $p:{\cal A} \to {\cal X}$ provides an embedding 
${\cal A}/H_{\cal A} \hra {\cal X}$. There is a canonical isomorphism $H_{\cal A}\otimes \Q = 
H_{\cal X}\otimes \Q$. 
In particular $H_{\cal A}(\R_{>0}) = H_{\cal X}(\R_{>0})$.  

\vskip 3mm
The {\it chiral dual} to a feed ${\bf i} = (I, \varepsilon, d)$ is 
 the 
feed 
${\bf i}^o:= (I, -\varepsilon, d)$. Mutations commute with the chiral duality on feeds. Therefore given a cluster 
${\cal X}$-variety ${\cal X}$ (respectively ${\cal A}$-variety ${\cal A}$), 
there is the 
chiral dual cluster 
${\cal X}$-variety (respectively ${\cal A}$-variety) denoted by 
 ${\cal X}^o$ (respectively ${\cal A}^o$).

The cluster ${\cal X}$-varieties 
${\cal X}$ and ${\cal X}^o$ are canonically isomorphic as schemes; 
for every feed, 
the isomorphism is given by inversion of the coordinates: 
$X_i \lms X_i^{-1}$. However this isomorphism changes the Poisson 
bracket to the opposite one: the Poisson structure 
on ${\cal X}^o$ differs from the one on ${\cal X}$ by the sign. 
The chiral dual cluster ${\cal A}$-variety ${\cal A}^o$ is 
canonically isomorphic to ${\cal A}$. The isomorphism is given by the identity 
maps on each feed torus. However ${\bf W}_{{\cal A}^o} =  -{\bf W}_{{\cal A}}$.

The {\it Langlands dual} to a feed 
${\bf i} = (I, \varepsilon_{ij}, d)$ is the feed 
 ${\bf i}^{\vee} = (I, \varepsilon^{\vee}_{ij}, d^{\vee})$, where 
$\varepsilon^{\vee}_{ij} := -\varepsilon^{\vee}_{ji}$ and $d^{\vee} = d^{-1}$. 
The Langlands duality on feeds commutes with mutations. Therefore 
it gives rise to the Langlands dual cluster 
${\cal A}$-, and ${\cal X}$-varieties, denoted 
${\cal A}^{\vee}$ and ${\cal X}^{\vee}$. 

\subsection{The symplectic double}

We follow the same pattern as in the definition of cluster ${\cal X}$- and 
${\cal A}$-varieties: define feed 
tori, introduce the relevant structures on them,  glue the feed 
tori in a certain specific way, and show that 
the gluing respects the structures.   
We assign to each feed ${\bf i}$ a split torus 
${\cal D}_{\bf i}$, equipped with canonical coordinates  
$(B_i, X_i)$, $i \in I$. 
There is a Poisson structure on the torus 
${\cal D}_{\bf i}$, given in coordinates by 
\begin{equation}\label{zx1}
\{B_i, B_j\}= 0, \quad \{X_i, B_j\}= d_i^{-1}\delta_{ij}X_iB_j, \quad \{X_i, X_j\} = 
\widehat \varepsilon_{ij}X_iX_j, \quad \widehat \varepsilon_{ij}:= \varepsilon_{ij}d_j^{-1}.
\end{equation}

\vskip 3mm
There is a symplectic $2$-form $\Omega_{\bf i}$ on the torus ${\cal D}_{\bf i}$:
\begin{equation}\label{zx101}
\Omega_{\bf i}:=  -\frac{1}{2}\sum_{i,j\in I} d_i \cdot 
\varepsilon_{ij} d\log B_i \wedge d\log B_j -
 \sum_{i \in I} d_i \cdot d\log B_i \wedge d\log X_i.
\end{equation}
So $\Omega_{\bf i}
(\partial_{B_i}\wedge \partial_{B_j}) = \widetilde \varepsilon_{ij}$ 
thanks to the coefficient $1/2$ in front of the first term. 
The following lemma is equivalent to Lemma \ref{zxa} which we prove below. 
 
\begin{lemma}\label{zx}
The Poisson and symplectic structures on ${\cal D}_{\bf i}$ are compatible: 
The Poisson structure (\ref{zx1}) coincides with the one defined by the 
$2$-form $\Omega_{\bf i}$.  
\end{lemma}

{\bf Mutations}. 
Set 
$$
{\Bbb B}_k^+:= \prod_{i|\varepsilon_{ki}>0}B_i^{\varepsilon_{ki}}, \qquad
{\Bbb B}_k^-:= \prod_{i|\varepsilon_{ki}<0}B_i^{-\varepsilon_{ki}}.
$$
Given an element $k \in I$, let us define a birational transformation 
$\mu_k: {\cal D}_{\bf i} \to {\cal D}_{\bf i'}$. Abusing notation, we identify $I'$ and $I$ 
via  the map  
$\mu_k: I \to I'$. Let $\{X'_i, B'_i\}$ 
be the coordinates 
on the torus ${\cal D}_{\bf i'}$. Set 
\begin{equation}\label{zx1qt}
\mu_k^*: X' _{i} \lms \left\{\begin{array}{lll} X_k^{-1}& \mbox{ if $i=k$}  \\
 X_i(1+X_k^{-\sgn(\varepsilon_{ik})})^{-\varepsilon_{ik}}, & \mbox{ if $i \not = k$}.
\end{array} \right.
\end{equation}
\begin{equation}\label{zx1qtrt}
\mu_k^*: B'_{i} \lms 
\left\{\begin{array}{ll} 
B_i & \mbox{ if $i\not =k$}, \\
 \frac{{\Bbb B}^-_k + X_k {\Bbb B}_k^+}{B_k(1+X_k)} &
\mbox{ if $i = k$}.
\end{array} \right.
\end{equation}

\begin{definition}\label{zx4}
The symplectic double ${\cal D}$ is a scheme over $\Z$ obtained by gluing the feed tori ${\cal D}_{\bf i}$ using 
formulas (\ref{zx1qt})-(\ref{zx1qtrt}), and then taking the affine closure.  
\end{definition}

%It is easy to see that mutations for the double ${\cal D}$ 
%commute with the chiral duality. 
%We denote by ${\cal D}^o$ the chiral dual to 
%${\cal D}$. 
Given a Poisson variety ${\cal Y}$, we denote by ${\cal Y}^{\rm op}$ 
the same variety with the opposite Poisson structure. 

Let us introduce a notation
\begin{equation}\label{zx1q1}
\widetilde X_i:= X_i ~ \frac{{\Bbb B}_k^+}{{\Bbb B}_k^-} = X_i
\prod_{j\in I}B_j^{\varepsilon_{ij}}. 
\end{equation}

We denote by $(A_i, A_i^o)$ the coordinates on 
the cluster torus ${\cal A}_{\bf i}\times {\cal A}_{\bf i}$.  
Let 
$p_-, p_+$ be the projections of ${\cal A}\times {\cal A}$ onto the two 
factors. 
The key
 properties of the symplectic double are the following: 
\begin{theorem} \label{8.15.05.10} a) There are  
$\widehat \Gamma$-equivariant positive symplectic 
spaces  ${\cal A}$, ${\cal D}$ and  ${\cal X}$. 

b)  There is a map $\varphi: {\cal A}\times {\cal A} \to {\cal D}$, 
given in any cluster coordinate system by the formulas
$$
\varphi^*(X_i) = \prod_j A_j^{\varepsilon_{ij}}, \quad 
%(A^+_j)^{-\varepsilon_{ij}}
\varphi^*(B_i) = \frac{A_i^o}{A_i}.
$$
It respects the canonical $2$-forms: 
$
\varphi^*\Omega_{\cal D} = 
p_-^*\Omega_{\cal A} -  p_+^*\Omega_{{\cal A}}$. 

c) There is a  Poisson map $\pi: {\cal D} \to {\cal X}\times {\cal X}^{\rm op}$, given in any cluster 
coordinate system by 
$$
\pi^*(X_i \otimes 1) = X_i, \quad \pi^*(1 \otimes X_j) = 
\widetilde X_j.
$$

d) There are  commutative diagrams
$$
\begin{array}{ccc}
{\cal A}\times {\cal A} && \\
&\searrow \varphi &\\
\downarrow p\times p && {\cal D}\\
&\swarrow \pi &\\
{\cal X}\times {\cal X}^{\rm op} & & 
\end{array}
\qquad 
\begin{array}{ccc}
{\cal X}&\stackrel{j}{\hra} &{\cal D} \\
\downarrow && \downarrow \pi\\
\Delta_{\cal X}&\stackrel{}{\hra} & {\cal X}\times {\cal X}^{\rm op} 
\end{array}
$$

e)  The map $j$ is a Lagrangian  embedding. 
The  intersection of its image with each cluster torus is 
given by equations $B_i =1$, $i \in I$. 

f) There is an involutive isomorphism $i: {\cal D}\to {\cal D}^{\rm op}$ 
interchanging the two 
components of the map $\pi$. 
It is given in any cluster 
coordinate system by 
$
i^*(B_i) = B^{-1}_i, \quad i^*(X_i) = \widetilde X_i.
$ 

g) The map $\varphi$ is the quotient by the diagonal action of the 
torus $H_{\cal A}$. The map $\pi$ is the quotient by a free Hamiltonian action of the 
torus $H_{\cal A}$ on ${\cal D}$. 
Its commuting Hamiltonians 
are given by the Poisson composition map
$
{\cal D} \stackrel{\pi}{\lra} {\cal X} \times {\cal X}^{\rm op} {\lra} {\cal X} 
\stackrel{\theta }{\lra} H_{\cal X}  
%\qquad  m(h_1, h_2) = h_1h_2.
$  
followed by characters of $H_{\cal X}$.
% provide a family of commuting Hamiltonians on ${\cal D}$. 
%Their Hamiltonian flows provide a free action of the group $H_{\cal A}$ on ${\cal D}$. 

\end{theorem}

\vskip 3mm
The proof is postponed till Section 3.  
It uses a decomposition of mutations from Section 2.4.

%{\bf Proof}. Since $\pi^*(1\otimes X_i)/\pi^*(X_i\otimes 1)  = 
%\prod_{j}B_j^{\varepsilon_{ij}}$,  the claim 
%boils down to the definitions. It is equivalent to 
%the fact that $\theta\circ p:{\cal A} \to {\cal X} \to H_{\cal X}$ is the 
%projection to the point $e$.

\vskip 3mm
{\bf The symplectic groupoid structure.}  
The data  $({\cal X}, {\cal D}, \pi, i, j)$ describes a 
symplectic groupoid related to the Poisson space ${\cal X}$. A definition of a  
symplectic groupoid  see in \cite{We}. 
Namely, ${\cal X}$ 
is the space of objects, ${\cal D}$ is the space of morphisms.  
The two components $\pi_-$ and $\pi_+$ of $\pi$ 
provide the source and target maps.  There is an involution $i$ on ${\cal D}$,  
reversing  the sign of the symplectic structure and 
interchanging the source and the target maps. 
It is the inversion map. 
There is a partial composition 
${\cal D} \times {\cal D} \to {\cal D}$, defined if 
$\pi_+y_1 = \pi_-y_2$, so that 
of $\pi_-(y_1 \circ y_2) = \pi_-(y_1)$ and $\pi_+(y_1 \circ y_2) = \pi_+(y_2)$. The  subvariety of ${\cal D} 
\times {\cal D}^{\rm op} \times {\cal D}^{\rm op}$  consisting of triples 
$(y_1, y_2, y_3)$ such that $y_3$ is the composition of $y_1$ and $y_2$ is isotropic. 
The Lagrangian subspace 
$j({\cal X})$ is the subspace of the identity morphisms.

\vskip 3mm 
\subsection{Cluster linear algebra } 
Let $\Lambda_{\cal X}$ be a finite rank lattice (a free abelian group) 
 with a skew-symmetric bilinear form 
$(\ast, \ast)_{\cal X}$. Let $\{e_i\}$, $i \in I$ be  a basis in $\Lambda_{\cal X}$. 
Then there is a dictionary translating the notion of a feed into the linear algebra language:
\begin{equation} \label{DICT}
\mbox{Feed = (lattice, skew-symmetric form, basis, multiplier)= 
$\left(\Lambda_{\cal X}, (\ast, \ast)_{\cal X}, \{e_i\}, \{d_i\}\right)$.}  
\end{equation}
Indeed,  $\varepsilon_{ij}:= d_j(e_i, e_j)_{\cal X}$. 
Set $\widehat \varepsilon_{ij}:= (e_i, e_j)_{\cal X}$. 

Consider the following two lattices associated with $\Lambda_{\cal X}$: 
$$
\Lambda^*_{\cal X}:= {\rm Hom}_\Z(\Lambda_{\cal X}, \Z), \qquad 
\Lambda_{\cal D}:= \Lambda_{\cal X}\oplus \Lambda^*_{\cal X}.
$$  
The basis $\{e_i\}$ in $\Lambda_{\cal X}$
 provides the dual basis $\{e^{\vee}_j\}$ in $\Lambda^*_{\cal X}$.
We define a quasidual basis $\{f_j\}$ in $\Lambda^*_{\cal X}$ by setting $f_i:= d_i^{-1}e^{\vee}_i$.
There is a basis 
$\{e_i, f_j\}$ in $\Lambda_{\cal D}$.

\vskip 3mm
{\it Canonical maps}. 
(i) The bilinear form on the lattice $\Lambda_{\cal X}$ provides a canonical, i.e. defined without using a basis, 
  homomorphism
\begin{equation} \label{CMA}
p^*: \Lambda_{\cal X} \lra \Lambda^*_{\cal X}, \qquad e_i \lms \sum_{j \in I} (e_i, e_j)_{\cal X}e^{\vee}_j  = 
\sum_{j \in I}\varepsilon_{ij}f_j. 
\end{equation}

(ii) There is a canonical homomorphism
$$
\pi^*:= ({\rm Id}, {\rm Id}+ p^*): 
\Lambda_{\cal X} \oplus \Lambda_{{\cal X}} \lra \Lambda_{\cal D}, \qquad e_l \oplus  e_r
 \lms e_l + 
e_r + p^*(e_r).
$$

(iii) Denote by  $\{f_i, f_i^o\}$ the basis in $\Lambda^*_{\cal X} \oplus \Lambda^*_{{\cal X}}$. 
There is a canonical homomorphism
$$
\varphi^*:= \left (\matrix{p^*& - {\rm Id}_-\cr 0& {\rm Id}_+\cr}\right ):
\Lambda_{\cal D} \lra \Lambda^*_{\cal X} \oplus \Lambda^*_{{\cal X}}, \qquad e_i 
 \lms p^*(e_i), f_j \lms f_j^o-f_j.
$$

(iv) There is a canonical involution
$$
i^*: \Lambda_{\cal D} \lra 
\Lambda_{\cal D}, \qquad f_j \lms - f_j, \quad e_i \lms e_i + p^*(e_i). 
$$

\vskip 3mm {\it Canonical bilinear forms/ bivectors}. 
The bilinear form $(\ast, \ast)_{\cal X}$ 
can be viewed as an element $\omega^*_{\cal X} \in \Lambda^2 \Lambda^*_{\cal X}$.  
It is written in the basis as 
\begin{equation}\label{AK2}
\omega^*_{\cal X} = \frac{1}{2}\sum_{i,j\in I} \widetilde \varepsilon_{ij} \cdot f_i \wedge f_j, 
\qquad \widetilde \varepsilon_{ij}= d_i \varepsilon_{ij}. 
\end{equation}

We define a canonical  element $\omega_{\cal D} \in \Lambda^2\Lambda_{\cal D}$ by setting   
\begin{equation}\label{zx101a}
\omega_{\cal D} = -\sum_{i\in I}e_i^{\vee} \wedge e_i - \frac{1}{2}\cdot e_i^{\vee} \wedge p^*(e_i) 
=  -\frac{1}{2}\sum_{i,j\in I} 
\widetilde \varepsilon_{ij} \cdot f_i \wedge f_j -
 \sum_{i \in I} d_i \cdot f_i \wedge e_i. 
\end{equation}
It is evidently non-degenerate. Thus it determines a dual element 
$\omega^*_{\cal D} \in \Lambda^2 \Lambda^*_{\cal D}$, which we view as a skew-symmetric bilinear form 
$(\ast, \ast)_{\cal D}$ on the lattice $\Lambda_{\cal D}$. Clearly 
$i^{*}\omega^*_{\cal D} = - \omega^*_{\cal D}$.

\begin{lemma}\label{zxa}
The form $(\ast, \ast)_{\cal D}$ is given explicitly by   
\begin{equation} \label{12.12.04.2}%\widehat \varepsilon_{ij}
(e_i, e_j)_{\cal D} : = (e_i, e_j)_{\cal X}, \quad (e_i, f_j)_{\cal D} = d_i^{-1}
\delta_{ij}, \quad (f_i, f_j)_{\cal D} = 0. 
\end{equation}\end{lemma}

{\bf Proof}. The symplectic structure $\omega_{\cal D}$ on $\Lambda^*_{\cal D}$
provides an isomorphism ${\rm I}: \Lambda^*_{\cal D} \to \Lambda_{\cal D}$, given by 
$\langle a, {\rm I}(b)\rangle = \omega_{\cal D}(a,b)$. 
One easily checks that 
$$
{\rm I}(f^*_i)= d_ie_i + \sum_{j\in I} \widetilde \varepsilon_{ij}e_j, 
\quad {\rm I}(e^*_i)= d_if_i, \qquad {\rm I}^{-1}(f_i) = - d_i^{-1}e_i^*, \quad 
{\rm I}^{-1}(e_i) = d_i^{-1}f_i^* + \sum_{j\in I} \widehat \varepsilon_{ij}e^*_j.
$$
Since by definition $\omega^*_{\cal D}(a, b) = \omega_{\cal D}({\rm I}^{-1}a, {\rm I}^{-1}b)$, we have  
$$
\omega^*_{\cal D}(e_i, e_j) = \omega_{\cal D}(d_i^{-1}f_i^* + \sum_{k} \widehat \varepsilon_{ik}e^*_k, 
d_j^{-1}f_j^* + \sum_{s} \widehat \varepsilon_{js}e^*_s) = - d_i^{-1}d_j^{-1} \widetilde \varepsilon_{ij} + 
\widehat \varepsilon_{ij} + \widehat \varepsilon_{ij} = \widehat \varepsilon_{ij}. 
$$
and $\omega^*_{\cal D}(e_i, f_j) = d^{-1}_j\delta_{ij}$. The lemma is proved. 

\vskip 3mm
Denote the bivector related to the second summand in 
$\Lambda^*_{\cal X}\oplus \Lambda^*_{\cal X}$ by $\omega^{o}_{\cal X}$. 
Denote by $\Lambda^{\rm op}_{{\cal X}}$ the lattice $\Lambda_{{\cal X}}$ equipped 
with the opposite 
bilinear form $(\ast, \ast)^{\rm op}_{{\cal X}} := -(\ast, \ast)_{{\cal X}}$. 

Let $j^*: \Lambda_{\cal D}  = \Lambda_{\cal X}^*\oplus 
\Lambda_{\cal X}\longrightarrow \Lambda_{\cal X}$ be the canonical projection. 

\begin{proposition} \label{LR}
(i) There are commutative diagrams
\begin{equation} \label{8.9.07.1}
\begin{array}{ccc}
\Lambda_{\cal X}&\stackrel{j^*}{\longleftarrow } &\Lambda_{\cal D} \\
\uparrow = && \uparrow \pi^*\\
\Lambda_{\cal X}&\stackrel{({\rm id}, {\rm id})}{\longleftarrow } & \Lambda_{\cal X}\times \Lambda_{\cal X}^{\rm op}
\end{array}
\qquad \qquad \qquad 
\begin{array}{ccc}
\Lambda^*_{\cal X}\oplus \Lambda^*_{\cal X} && \\
&\nwarrow \varphi^* &\\
\uparrow p^*\times p^* && \Lambda_{\cal D}\\
&\nearrow \pi^* &\\
\Lambda_{\cal X}\oplus \Lambda_{\cal X}^{\rm op} & & 
\end{array}
\end{equation}
The map $\pi^*$ respects the bilinear forms in 
$\Lambda_{\cal X} \oplus \Lambda^{\rm op}_{{\cal X}}$ and $\Lambda_{\cal D}$. One has  
\begin{equation} \label{FMW}
\varphi^*(\omega_{\cal D}) = \omega_{\cal X} - \omega^{o}_{\cal X}.
\end{equation}

(ii) The involution $i$ 
 reverses the sign of the symplectic $(\ast, \ast)_{\cal D}$ form in $\Lambda_{\cal D}$. 
\end{proposition}

{\bf Proof}. (i) The commutativity is clear. 
The claim about $\pi^*$ is clear for the left component of $\pi^*$. 
For the right one it follows from (\ref{CH}). 
To check the claim about $\varphi^*$, notice that 
$$
-\frac{1}{2}\sum_{i,j\in I}
\widetilde \varepsilon_{ij} \cdot (f_i^o - f_i)\wedge 
(f_j^o - f_j) - \sum_{i,j\in I}d_i \cdot (f_i^o - f_i)
\wedge \varepsilon_{ij}\cdot f_j = 
\frac{1}{2}\sum_{i,j\in I}\widetilde \varepsilon_{ij}\cdot  (
f_i \wedge f_j -f_i^o \wedge f^o_j).
$$

(ii) Indeed, 
$
(i^*e_i, i^*f_j)_{\cal D} =  (e_i + p^*(e_i), -f_j)_{\cal D} = -(e_i, f_j)_{\cal D}. 
$ and 
\begin{equation} \label{CH}
 (i^*e_i, i^*e_j)_{\cal D} = (e_i + \sum_{s \in I}\varepsilon_{is}f_s, 
e_j+ \sum_{t \in I}\varepsilon_{jt}f_t)_{\cal D} = \widehat \varepsilon_{ij} - 
\widehat \varepsilon_{ij} + \widehat \varepsilon_{ji} = -\widehat \varepsilon_{ij}.  
\end{equation}

\vskip 3mm
{\it Mutated bases}. Set $[\alpha]_+= \alpha$ if $\alpha\geq 0$ and $[\alpha]_+=0$ otherwise. 
Given $k \in I$, 
choose another basis $\{e_i'\}$ in $\Lambda_{\cal X}$ given by 
\begin{equation} \label{12.12.04.2a}
e'_i := 
\left\{ \begin{array}{lll} e_i + [\varepsilon_{ik}]_+e_k
& \mbox{ if } &  i\not = k\\
-e_k& \mbox{ if } &  i = k.\end{array}\right.
\end{equation}
\begin{lemma} \label{12.12.04.1}
The quasidual basis $\{f_i'\}$ in $\Lambda^*_{\cal X}$ is given by 
\begin{equation} \label{12.12.04.2ca}
 \qquad 
f'_i:= 
\left\{ \begin{array}{lll} - f_k + \sum_{j\in J}[-\varepsilon_{kj}]_+f_j
& \mbox{ if } &  i= k\\
f_i& \mbox{ if } &  i \not = k.\end{array}\right. 
\end{equation}
\end{lemma} 

{\bf Proof}. Let $t_k$ be the automorphism of $\Lambda_{\cal X}$ 
given by $e_i \lms e_i'$. 
The matrix of the dual map $t_k^*$ in the dual basis $e_i^{\vee}$ is obtained by 
the transposition of the matrix of the map $t_k$. Thus the map $t_k^*$ acts as 
$e_k^{\vee} \lms - e_k^{\vee} + \sum_j[\varepsilon_{jk}]_+e_j^{\vee}$ and $e_i^{\vee} \lms e_i^{\vee} $ if $i \not = k$. 
Changing to the basis $f_j$ we get the lemma. 

%Assume that $k \not \in \{i,j\}$. 
%Then $
%(e'_i, f'_j)  = ( e_i + [\varepsilon_{ik}]_+e_k, f_j) = (e_i, f_j)$, 
%$$
%(e'_k, f'_i)  = (-e_k, f_i) = 0, 
%\qquad (e'_k, f'_k)  = (-e_k, - f_k + \sum_j[-\varepsilon_{kj}]_+f_j)  = (e_k, f_k), 
%$$ 
%$$
%(e'_i, f'_k)  = (e_i + [\varepsilon_{ik}]_+e_k, - f_k + \sum_j[-\varepsilon_{kj}]_+f_j)
%= d_k^{-1}[\varepsilon_{ik}]_+ - d_i^{-1}[-\varepsilon_{ki}]_+ = 
%[\widehat \varepsilon_{ik}]_+ - [\widehat \varepsilon_{ik}]_+ = 0.
%$$

\begin{lemma} \label{12.12.04.1a}
Set $\widehat \varepsilon'_{ij}:= (e'_i, e'_j)_{\cal X}$. Then one has 
 \begin{equation} \label{9.28.04.5a}
\widehat \varepsilon'_{ij}=
\left\{ \begin{array}{ll} -\widehat \varepsilon_{ij}, & \mbox{ if }  i=k \mbox{ or } j=k\\
\widehat \varepsilon_{ij}  & \mbox{ if } \varepsilon_{ik}\varepsilon_{kj}
 \leq 0\\
\widehat \varepsilon_{ij} + |\varepsilon_{ik}|\widehat  
\varepsilon_{kj} & \mbox{ if } \varepsilon_{ik}\varepsilon_{kj}>0.
\end{array}\right.
\end{equation}
%The bilinear form $(\ast, \ast)$ on $\Lambda_{\cal D}$ 
%in the basis $\{e'_i, f'_j\}$ is given by formula (\ref{12.12.04.2}) with 
%$\widehat \varepsilon_{ij}$ replaced by $\widehat \varepsilon'_{ij}$.
\end{lemma}

{\bf Proof}. Clearly 
$
(e'_i, e'_k)  = (e_i + [\varepsilon_{ki}]_+e_k, - e_k) 
= - \widehat \varepsilon_{ik} = \widehat \varepsilon'_{ik}.  
$ 
Assume that $k \not \in \{i,j\}$. Then 
$$
(e'_i, e'_j) 
= (e_i + [\varepsilon_{ik}]_+e_k, e_j + [\varepsilon_{jk}]_+e_k) 
= 
\widehat \varepsilon_{ij} + [\varepsilon_{ik}]_+\widehat \varepsilon_{kj}
+ \widehat \varepsilon_{ik}[\varepsilon_{jk}]_+  
$$
$$
=\widehat \varepsilon_{ij} + [\varepsilon_{ik}]_+\widehat \varepsilon_{kj}
+ \varepsilon_{ik}[-\widehat \varepsilon_{kj}]_+ = \widehat
\varepsilon'_{ij}. 
$$
The lemma is proved. 

\vskip 3mm
{\bf Remarks}. 1. Since the multipliers $d_i$ do not change under mutations, 
the mutation formula for $\varepsilon_{ij}$ is
 equivalent to formula  (\ref{9.28.04.5a}). 

2. The map $t_k$ is an automorphism of the lattice $\Lambda_{\cal X}$.
However it does not preserve the bilinear form $(*,*)_{\cal X}$ 
on $\Lambda_{\cal X}$.
Its square $t_k^2$ is not 
the identity map in general. 
However it preserves the  bilinear form $(*,*)_{\cal X}$ 
on $\Lambda_{\cal X}$. 

\vskip 3mm

We define a category of algebraic tori (viewed as group schemes) 
as the the category 
opposite to the category  of abelian groups of finite rank: 
A finite rank abelian group $A$ gives rise to the torus 
${\rm Hom}(A, {\Bbb G}_m)$. So if $A$ is finite, it is a torsion group, a product of finite 
multiplicative groups, 
and if $A$ is torsion free, it is a product of copies of ${\Bbb G}_m$. It is an abelian category. 

Having in mind dictionary (\ref{DICT}), 
denote by $\Lambda_{{\cal D}, {\bf i}}$ the lattice $\Lambda_{{\cal D}}$ equipped with the basis 
corresponding to the feed ${\bf i}$, and similarly for the lattices $\Lambda_{{\cal X}}$ and 
$\Lambda^*_{\cal X}$. Then 
for the  feed tori we have: 
\begin{equation} \label{MISq}
{\cal D}_{\bf i} = {\rm Hom}(\Lambda_{{\cal D}, {\bf i}}, ~{\Bbb G}_m), \quad 
{\cal X}_{\bf i} = {\rm Hom}(\Lambda_{{\cal X}, {\bf i}}, ~{\Bbb G}_m), \quad 
{\cal A}_{\bf i} = {\rm Hom}(\Lambda^*_{{\cal X}, {\bf i}}, ~{\Bbb G}_m).
\end{equation} 
The cluster coordinates for these tori are  
provided by the corresponding lattice bases.

We define the cluster tori $H_{\cal X}$ and $H_{\cal A}$ by setting
$$
H_{\cal X}:= {\rm Hom}({\rm Ker}~p^*, ~{\Bbb G}_m), \qquad 
H_{\cal A}:= {\rm Hom}({\rm Coker}~p^*, ~{\Bbb G}_m). 
$$
So $H_{\cal X}$ is torsion free, while $H_{\cal A}$ may have torsion. 
The exact sequence of abelian groups $0 \to {\rm Ker} ~p^* \to 
{\Lambda}_{\cal X} \stackrel{p^*}{\to} 
{\Lambda}_{\cal X}^* \to {\rm Coker} ~p^* \to 0$ 
gives rise to an exact exact sequence of tori
\begin{equation} \label{LEST}
0 \lra H_{\cal A}\lra {\cal A}_{\bf i} \lra {\cal X}_{\bf i} \lra H_{\cal X}\lra 0.
\end{equation} 
%The group $H_{\cal A}$ acts freely on ${\cal A}_{\bf i}$.

\subsection{Decomposition of mutations}

In this subsection we decompose mutations of the 
${\cal A}$-, ${\cal D}$- and ${\cal X}$-spaces into a composition of 
an automorphism and a monomial transformation,  so that
$$
\mbox{Monomial part of the mutation = basis change from Section 2.3.}
$$
Given $k \in I$, the bases changes in the lattices can be interpreted as 
 isomorphisms 
\begin{equation} \label{MIS}
\tau_k: {\Lambda}_{{\cal D}, \bf i'} \lra {\Lambda}_{{\cal D}, \bf i}, \quad 
\tau_k: {\Lambda}_{{\cal X}, \bf i'} \lra {\Lambda}_{{\cal X}, \bf i}, \quad 
\tau_k: {\Lambda}^*_{{\cal X}, \bf i'} \lra {\Lambda}_{{\cal X}, \bf i}.
\end{equation} 

\begin{definition} \label{MUS}
The mutation isomorphisms 
$$
\mu_k': {\cal D}_{\bf i} \lra {\cal D}_{\bf i'}, \quad 
\mu_k': {\cal X}_{\bf i} \lra {\cal X}_{\bf i'}, \quad 
\mu_k': {\cal A}_{\bf i} \lra {\cal A}_{\bf i'}
$$
are the isomorphisms determined by the maps of the lattices (\ref{MIS}) and (\ref{MISq}).
\end{definition}

Lemmas \ref{12.12.04.1} and \ref{12.12.04.1a} imply

\begin{corollary} \label{Cor}
The maps $\mu_k': {\cal D}_{\bf i} \lra {\cal D}_{\bf i'}$ and 
$\mu_k': {\cal X}_{\bf i} \lra {\cal X}_{\bf i'}$ are Poisson maps. 
\end{corollary}

\vskip 3mm
Let us proceed now to decomposition of mutations. 

\vskip3mm
{\it The ${\cal D}$-space}. Let us define a birational 
 automorphism $\mu_k^{\sharp}$ of the feed torus 
${\cal D}_{\bf i}$ which acts on the coordinates as follows: 
\begin{equation} \label{f3**}
B_i \lms B^{\sharp}_{i} := \left\{\begin{array}{lll} B_i& \mbox{ if } & i\not =k, \\
    B_k(1+X_k)(1+\widetilde X_k)^{-1} & \mbox{ if } &  i= k. \\
\end{array} \right.
\end{equation}
\begin{equation} \label{f3*yyy}
X_i \lms X^{\sharp}_{i} := 
    X_i(1+X_k)^{-\varepsilon_{ik}}
\end{equation}
Denote by $(B_i', X_i')$ the coordinates of the feed torus 
${\cal D}_{\bf i'}$. 
The  isomorphism 
$
\mu'_k: {\cal D}_{\bf i} \lra {\cal D}_{\bf i'}
$ 
acts on the coordinates as follows 
(observe that ${\Bbb B}_k^+ = {{\Bbb B}^{\sharp}_k}^+$): 

\begin{equation} \label{11.18.06.10hr}
B'_{i} \lms \left\{\begin{array}{lll} B_i& \mbox{ if } & i\not =k, \\
    {\Bbb B}_k^-/B_k
 & \mbox{ if } &  i= k. \\ 
\end{array} \right. \qquad
X'_{i} \lms \left\{\begin{array}{lll} X_k^{-1}& \mbox{ if } & i=k, \\
    X_i(X_k)^{[\varepsilon_{ik}]_+} & \mbox{ if } &  i\neq k. \\
\end{array} \right.
\end{equation}

\begin{lemma} \label{11.18.06.12} The mutation $\mu_k: {\cal D}_{\bf i} \lra {\cal D}_{\bf i'}$ 
is decomposed into a composition 
$
\mu_k = \mu'_k\circ \mu^{\sharp}_k. 
$
\end{lemma}

{\bf Proof}. We should check that on the level of functions we have 
$\mu^*_k = (\mu'_k)^*\circ (\mu^{\sharp}_k)^*$. Indeed, 
$$
B_k' \lms {\Bbb B}_k^-/B_k \lms 
\frac{1+X_k{\Bbb B}_k^+/{\Bbb B}_k^-}{B_k(1+X_k)}{\Bbb B}_k^- =
\frac{{\Bbb B}_k^- + X_k{\Bbb B}_k^+}{B_k(1+X_k)}.
$$
The rest is an obvious computation. 
The lemma is proved. 

\vskip 3mm
{\it The ${\cal X}$-space}. The decomposition is given by 
formulas (\ref{f3*yyy}) and the right formula in (\ref{11.18.06.10hr}).

\begin{lemma} \label{11.18.06.12trtrd} 
The mutation $\mu_k: {\cal X}_{\bf i} \lra {\cal X}_{\bf i'}$ 
is decomposed into a composition 
$
\mu_k = \mu'_k\circ \mu^{\sharp}_k. 
$
\end{lemma}

{\bf Proof}. We have 
$
X_k' \lms X_k^{-1} \lms X_k^{-1}$, and  
$$
X_i' \lms X_i(X_k)^{[\varepsilon_{ik}]_+}  \lms X_i(1+X_k)^{-\varepsilon_{ik}} 
(X_k)^{[\varepsilon_{ik}]_+} = X_i(1+X_k^{-{\rm sgn}(\varepsilon_{ik})})^{-\varepsilon_{ik}}.
$$

\vskip 3mm
{\it The ${\cal A}$-space}. Set 
$$
{\Bbb A}_k^+:= 
\prod_{j\in I}A_j^{[\varepsilon_{kj}]_+}, \qquad {\Bbb A}_k^-:= 
\prod_{j\in I}A_j^{[-\varepsilon_{kj}]_+}.
$$ We define a birational 
 automorphism $\mu_k^{\sharp}$ of the feed torus 
${\cal A}_{\bf i}$ which acts on the coordinates as follows: 
\begin{equation} \label{f3**}
A_i \lms A^{\sharp}_{i} := \left\{\begin{array}{lll} A_i& \mbox{ if } & i\not =k, \\
    A_k(1+p^*X_k)^{-1} = A_k(1+ {\Bbb A}_k^+/{\Bbb A}_k^-)^{-1} & \mbox{ if } &  i= k. \\
\end{array} \right.
\end{equation}

Denote by $A_i'$ the coordinates of the feed torus 
${\cal A}_{\bf i'}$. 
The isomorphism  
$
\mu'_k: {\cal A}_{\bf i} \lra {\cal A}_{\bf i'}
$ 
acts on them as follows: 
\begin{equation} \label{11.18.06.10sdf}
A'_{i} \lms \left\{\begin{array}{lll} 
A_i& \mbox{ if } & i\not =k, \\
    {\Bbb A}_k^-/A_k
 & \mbox{ if } &  i= k. \\
\end{array} \right.
\end{equation}

\begin{lemma} \label{11.18.06.12trtr} 
The mutation $\mu_k: {\cal A}_{\bf i} \lra {\cal A}_{\bf i'}$ 
is decomposed into a composition 
$
\mu_k = \mu'_k\circ \mu^{\sharp}_k. 
$
\end{lemma}

{\bf Proof}. We need to check only 
how the composition acts on $A_k'$. 
It is given by 
$$
A'_{k} \lms {\Bbb A}_k^-/A_k= 
{\Bbb A}_k^-(1+ \frac{{\Bbb A}_k^+}{{\Bbb A}_k^-})A_k^{-1}
 = \frac{{\Bbb A}_k^+ + {\Bbb A}_k^-}{A_k}.
$$

\vskip 3mm
Decomposition of mutations is compatible
 with the map $p:{\cal A}_{\bf i} \to {\cal X}_{\bf i}$: 
\begin{lemma} \label{11.23.06.1} 
The map $p: {\cal A}_{\bf i} \to {\cal X}_{\bf i}$ 
intertwines each of the components $\mu'_k$,
 $ \mu^{\sharp}_k$ of the mutation map. 
\end{lemma}

{\bf Proof}. We have to show that the following diagram is commutative:
\begin{equation}\label{7.7.04.13ww}
\begin{array}{ccccc}
{\cal A}_{\bf i}& \stackrel{\mu^{\sharp}_{k}}{\lra} &
{\cal A}_{{{\bf i}}}&\stackrel{\mu'_{k}}{\lra}&{\cal A}_{{{\bf i}'}}\\
p\downarrow && \downarrow p &&\downarrow p\\
{\cal X}_{\bf i} &\stackrel{\mu^{\sharp}_{k}}
{\lra}&{\cal X}_{{{\bf i}}}&
\stackrel{\mu'_{k}}{\lra} &
{\cal X}_{{{\bf i}'}}
\end{array}
\end{equation}
In the left square,
going up and to the left we get
$
X_i\lms \prod_{j\in I}A_j^{\varepsilon_{ij}} \lms 
\prod_{j\in I}A_j^{\varepsilon_{ij}} \cdot 
(A_k^{\sharp}/A_k)^{\varepsilon_{ik}}. 
$ 
 Going to the left and up we get the same:
$$
X_i \lms X_i(1+X_k)^{-\varepsilon_{ik}} \lms 
\prod_{j\in I}A_j^{\varepsilon_{ij}} 
(1+\frac{{\Bbb A}_k^+}{{\Bbb A}_k^-})^{-\varepsilon_{ik}}.
$$
The claim for the right square is equivalent to the 
existence of the map $p^*$ of lattices, see (\ref{CMA}).

\vskip 3mm

\subsection{The motivic dilogarithm structure on the symplectic double} 
Recall that for a field $F$ the abelian group 
$K_2(F)$ is the quotient of the abelian group $\Lambda^2F^*$, the wedge square of the multiplicative group $F^*$ 
of $F$,  by the subgroup generated by the so-called 
Steinberg relations $(1-x)\wedge x$, $x \in F^*-\{1\}$. Further, if $F$ is the field of rational 
functions $\Q(X)$ on a variety $X$, there is a canonical map 
$$
d\log: \Lambda^2\Q(X)^* \lra \Omega_{\rm log}^2(X), \quad f\wedge g \lms d\log (f) \wedge d\log (g).
$$
where $\Omega_{\rm log}^2(X)$ is the space of $2$-forms with logarithmic singularities on $X$. 
It apparently kills the Steinberg relations, 
inducing a map $K_2(\Q(X)) \lra \Omega_{\rm log}^2(X)$. 

Working modulo $2$-torsion, one can replace 
the Steinberg relation by the following one: 
\begin{equation}\label{11.5.06.8}
(1+x) \wedge x = (1-(-x)) \wedge x = (1-(-x)) \wedge (-x) +  
(1-(-x)) \wedge -1.
\end{equation}

The symplectic structure on ${\cal D}$ comes from a class in $K_2$ as follows.  
For each feed ${\bf i}$, the 
symplectic form  is the image under the $d\log$ map of  
an element $W_{\bf i} \in \Lambda^2F_{\bf i}^*$, where $F_{\bf i}$ is the function field of the feed torus. 
For a mutation ${\bf i} \to {\bf i}'$ we present the difference $W_{\bf i} - W_{{{\bf i}'}}$ 
as a difference of Steinberg relations (\ref{11.5.06.8}). 
So the 
class of the element $W_{\bf i}$ in $K_2$ is mutation invariant. 
The difference of Steinberg relations  gives rise to an element of the Bloch group, which is 
the motivic avatar of  generating function of the 
symplectic map ${\cal D}_{\bf i}\to {\cal D}_{\bf i'}$, a mutation  
of the feed tori. Let us implement this program. 

\vskip 3mm
We use $\cdot$ for the product of an element 
of the wedge product by an integer.   
The element $\omega_{\cal D}$, see (\ref{zx101a}),  provides an element 
\begin{equation}\label{ELL}
W_{\bf i}= -\frac{1}{2}\sum_{i,j\in I} \widetilde 
\varepsilon_{ij}\cdot  B_i \wedge B_j - \sum_{i \in I} d_i \cdot B_i \wedge X_i
  \in \Lambda^2\Q({\cal D}_{\bf i})^*, 
  \quad \widetilde 
\varepsilon_{ij}:= d_i  
\varepsilon_{ij}.
\end{equation}
Applying to $W_{\bf i}$ the map $d\log$ 
we get the $2$-form $\Omega_{\bf i}$ on the torus ${\cal D}_{\bf i}$, see 
(\ref{zx101}).

\begin{proposition}\label{zx2f}
Given a mutation ${\bf i} \to {\bf i'}$ in the direction $k$, one has 
$$
\mu_k^*W_{\bf i'}  - W_{\bf i} = d_k\cdot \left( 
(1+\widetilde X_k) \wedge \widetilde X_k - 
(1+X_k) \wedge X_k\right).
$$ 
\end{proposition}

{\bf Proof}. 
It follows from the two claims, describing behavior of the elements 
$W_{\bf i}$ under the automorphism 
$\mu^{\sharp}_k$ and the map $\mu_k'$. It is convenient to 
set $W^{\sharp}_{\bf i}:= (\mu^{\sharp}_k)^*W_{\bf i}$. 

\begin{lemma} \label{11.18.06.1sr} One has 
$
W^{\sharp}_{\bf i} - W_{\bf i} = d_k\cdot \Bigl(
(1+\widetilde X_k) \wedge \widetilde X_k -
(1+X_k) \wedge X_k \Bigr).
$
\end{lemma}

{\bf Proof}. We have 
\begin{equation} \label{5.11.03.3a}
-2 \cdot (W^{\sharp}_{{\bf i}} - W_{\bf i}) = 
\sum_{i \not = k, j \not = k}d_i\varepsilon_{ij}\cdot B^{\sharp}_i \wedge B^{\sharp}_j + 
d_k\sum_{j}\varepsilon_{kj}\cdot  B^{\sharp}_k \wedge B_j^{\sharp} + 
d_i \sum_{i}\varepsilon_{ik}\cdot  B^{\sharp}_i \wedge B_k^{\sharp}
\end{equation}
\begin{equation} \label{5.11.03.3ayu}
-\Bigl(\sum_{i \not = k, j \not = k}d_i\cdot \varepsilon_{ij}\cdot B_i \wedge B_j + 
d_k\sum_{j}\varepsilon_{kj}\cdot  B_k \wedge B_j + 
d_i\sum_{i}\varepsilon_{ik}\cdot  B_i \wedge B_k\Bigr)
\end{equation}
\begin{equation} \label{5.11.03.3adt}
+ 2 \sum_i d_i \cdot B_i^{\sharp}\wedge X_i^{\sharp} - 
2 \sum_i d_i \cdot B_i\wedge X_i. 
\end{equation}
Since $B_i^{\sharp}= B_i$ for $i \not = k$, and $X_k^{\sharp} = X_k$, we have 
\begin{equation} \label{5.11.03.3d}
(\ref{5.11.03.3a}) + (\ref{5.11.03.3ayu}) = 
d_k\cdot \sum_{j}\varepsilon_{kj}\cdot  \frac{B^{\sharp}_k}{B_k} \wedge B_j 
+d_i\cdot \sum_{i}\varepsilon_{ik}\cdot  B_i \wedge \frac{B^{\sharp}_k}{B_k} = 
2d_k\cdot \frac{B^{\sharp}_k}{B_k} \wedge 
\frac{{\Bbb B}^+_k}{{\Bbb B}^-_k}.
\end{equation}

Further, formula (\ref{5.11.03.3adt}) equals
\begin{equation} \label{5.11.03.3rt}
2 \sum_{i \not = k}d_i \cdot B_i\wedge (X_i^{\sharp}/X_i) 
+ 2 d_k \frac{B^{\sharp}_k}{B_k}\wedge X_k = 
-2\sum_{i \not = k}d_i\varepsilon_{ik} \cdot B_i\wedge (1+X_k)
+ 2d_k \cdot \frac{B^{\sharp}_k}{B_k}\wedge X_k  
\end{equation}
\begin{equation} \label{5.11.03.3ui}
=2d_k \cdot \Bigl(\frac{{\Bbb B}_k^+}{{\Bbb B}_k^-} \wedge (1+X_k) 
+ \frac{B^{\sharp}_k}{B_k}\wedge X_k\Bigr). 
\end{equation}
Here we used $-d_i\varepsilon_{ik} = d_k\varepsilon_{ki}$ to get the third
 equality. 
Therefore we get
$$
W^{\sharp}_{\bf i}  - W_{\bf i} = -(\ref{5.11.03.3d}) - (\ref{5.11.03.3ui})=
-d_k\cdot \Bigl(\frac{B^{\sharp}_k}{B_k}
 \wedge 
X_k\frac{{\Bbb B}^+_k}{{\Bbb B}^-_k} + 
\frac{{\Bbb B}_k^+}{{\Bbb B}_k^-} \wedge (1+X_k) \Bigr) = 
$$
$$
-d_k\cdot \Bigl((1+X_k)(1+\widetilde X_k)^{-1} \wedge 
X_k\frac{{\Bbb B}^+_k}{{\Bbb B}^-_k} + 
\frac{{\Bbb B}_k^+}{{\Bbb B}_k^-} \wedge (1+X_k) \Bigr)=
d_k\cdot \left((1+\widetilde X_k) \wedge \widetilde X_k
- (1 + X_k) \wedge 
X_k \right). 
$$
The lemma is proved. The next Lemma is an immediate corollary of Corollary \ref{Cor}.

\begin{lemma} \label{11.18.06.12a} One has 
$
(\mu'_k)^{*}W_{\bf i'} = W^{\sharp}_{\bf i}.
$
\end{lemma}

Lemmas \ref{11.18.06.1sr}
 and \ref{11.18.06.12a} imply the proposition.  

\vskip 3mm
The element $W_{\bf i}$ gives rise to an element of $K_2$ of the 
field $\Q({\cal D}_{\bf i})$, which evidently lies in 
$K_2({\cal D}_{\bf i})$. Proposition \ref{zx2f} implies that 
mutations preserve these elements. 
So they give rise to a $\Gamma$-invariant 
element ${\bf W} \in K_2({\cal D})$. 
Thus the double ${\cal D}$ 
has a $\Gamma$-invariant symplectic structure  
$\Omega = d\log({\bf W})$. 
\vskip 3mm

{\it The case of the ${\cal A}$-space}.    
The element $\omega^*_{\cal X}$, see (\ref{AK2}), 
provides an element 
\begin{equation}\label{zx100}
W_{{\cal A}, {\bf i}}:= \frac{1}{2}\cdot \sum_{i,j\in I} \widetilde \varepsilon_{ij}A_i 
\wedge A_j \in \Lambda^2\Q({\cal A}_{\bf i})^*. 
\end{equation}

\begin{lemma}\label{12.4.06.2}
Given a mutation ${\bf i} \to {\bf i'}$ in the direction $k$, one has 
$$
(\mu_k^{\sharp})^*W_{{\cal A}, {\bf i'}}  - W_{{\cal A}, {\bf i}} =  d_k\cdot 
(1+X_k) \wedge X_k, \qquad (\mu'_k)^*W_{{\cal A}, {\bf i}}  =  W_{{\cal A}, {\bf i}}.
$$ 
\end{lemma}

{\bf Proof}. The second follows from Lemma \ref{12.12.04.1}. For the first, one has  
$$
(\mu_k^{\sharp})^*W_{{\cal A}, {\bf i}}  - W_{{\cal A}, {\bf i}} = 
\frac{1}{2}\sum_{i,j\in I}\widetilde \varepsilon_{ij}\cdot 
\left( A^{\sharp}_i\wedge A^{\sharp}_j - A_i\wedge A_j \right) = 
\sum_{i\in I}\widetilde \varepsilon_{ik}
\left( A_i\wedge A^{\sharp}_k/A_k \right) = 
\left(1+ {\Bbb A}_k^+/{\Bbb A}_k^-\right)\wedge {\Bbb A}_k^+/{\Bbb A}_k^-.
$$

\vskip 3mm
\begin{corollary}
There is a $\Gamma$-invariant 
element ${\bf W}_{\cal A} \in K_2({\cal A})$, and thus 
a $\Gamma$-invariant presymplectic structure  
$\Omega = d\log({\bf W}_{\cal A})$ on the space ${\cal A}$. 
\end{corollary}

\vskip 3mm
\subsection{The unitary part of the symplectic double.} 
Since ${\cal D}$ is a positive space, the set of its 
positive real points is well defined. It is a real 
symplectic subspace of ${\cal D}(\C)$, with the symplectic form 
${\rm Re}\Omega$. 

It turns out that 
there is another real 
subspace ${\cal D}^U$ of ${\cal D}(\C)$,
 with a real symplectic structure given by 
${\rm Im}\Omega$, which we call the 
unitary part of the complex double. 
It is defined  
by setting the $A$-coordinates unitary, 
and the phases of the $X$-coordinates expressed via the $A$-coordinates: 

 \begin{definition}
The unitary part ${\cal D}^U$ of ${\cal D}(\C)$ is the real subspace 
obtained by gluing the real subvarieties of the 
complex cluster tori defined by the following equations: 
\begin{equation} \label{op}
|B_i| = 1, \quad 
\frac{\overline X_i }{X_i} = \prod_{j\in I}B_j^{\varepsilon_{ij}}, \quad i \in  I.
\end{equation} 
\end{definition}

Let $\Sigma_{\cal X}$ be the ``antiholomorphic diagonal'' in 
${\cal X}(\C) \times {\cal X}(\C)$, defined as the set of stable points of the 
composition of the antiholomorphic involution on 
${\cal X}(\C) \times {\cal X}(\C)$ 
with the one $(x_1, x_2) \lms (x_2, x_1)$:
$$
\Sigma_{\cal X} = \{(x_1, x_2) \in {\cal X}(\C) \times {\cal X}(\C) \quad | \quad (x_1, x_2) = (\overline x_2, \overline x_1)\}.
$$

\begin{proposition} a) There is a commutative square 
$$
\begin{array}{ccc}
{\cal D}^U& \stackrel{j^U}{\hra} &{\cal D}(\C)\\
\rho \downarrow &&\downarrow \pi\\
\Sigma_{\cal X}& \hra& {\cal X}(\C) \times {\cal X}(\C)
\end{array}
$$
The modular group acts on ${\cal D}^U$ by real 
birational transformations. 

b) The space ${\cal D}^U$ has a symplectic structure given by 
${\rm Im}\Omega$. 

c) The space ${\cal D}^U$ is a real Lagrangian subspace for 
the real symplectic structure ${\rm Re}\Omega$ on ${\cal D}(\C)$. 
\end{proposition}

{\bf Proof}. 
a) The commutativity just means that 
$\pi_-^*\overline X_i = \pi_+^*X_i$, that is 
$
\overline X_i = \widetilde X_i.
$ 
This is equivalent to the second equation in (\ref{op}). 
Let us check that the gluing respects the equations (\ref{op}):
$$
|B_kB_k'| = \Big |\frac{{\Bbb B}_k^- + X_k {\Bbb B}_k^+}{1+X_k}\Big | = 
\Big |\frac{1 + X_k {\Bbb B}_k^+/{\Bbb B}_k^- }{1+X_k}\Big | = 
\Big |\frac{1 + \overline X_k}{1+X_k}\Big | =1.
$$
Since the second condition in (\ref{op}) just means that the image of the map 
$\pi$ lies in the antiholomorphic diagonal, it is also invariant. 

b, c) Set $B_j = |B_j| e^{\sqrt{-1}b_j}$ and $X_j = r_je^{ \sqrt{-1} x_j}$. 
Then the restriction of the holomorphic $2$-form $\Omega$ to the unitary 
part ${\cal D}^U$ is given by 
\begin{equation} \label{popo} 
\Omega_{|{{\cal D}^U}}= -2\sum_{i,j\in I}d_i \varepsilon_{ij} d b_i \wedge db_j + 
\sum_{i\in I}\sum_{j\in I} d_i \varepsilon_{ij} db_i \wedge db_j + 
\sqrt{-1}\sum_{i \in I}da_i \wedge dr_i = \sqrt{-1}\sum_{i \in I}db_i \wedge dr_i.
\end{equation}  
The proposition is proved.

\vskip 3mm
If $\varepsilon_{ij}$ is non-degenerate, 
there exists a complex structure on ${\cal D}^U$ for which the maps
$\pi_-, \pi_+: {\cal D}^U \lra {\cal X}(\C)$ are, respectively, holomorphic and  antiholomorphic.

%\subsection{The canonical $\Gamma$-equivariant line bundle 
%with connection over ${\cal D}(\C)$}

%\input{gerbdouble.tex}

\subsection{Appendix}

\paragraph{Relating the quantum spaces ${\cal X}_q$ and  
${\cal A}_q$ for nondegenerate $\varepsilon_{ij}$} 
Recall the quantum cluster algebras 
of  Berenstein and Zelevinsky \cite{BZq}. 
They are defined using some additional data (the companion matrix). 
However if ${\rm det}(\varepsilon_{ij})\not = 0$,  
there is a canonical 
non-commutative $q$-deformation of the cluster algebra related to a
feed with the cluster function $\varepsilon_{ij}$. 
It is the one which we will use. 
It is easyly translated into the languige  of 
quantum spaces.  Given a feed ${\mathbf i}$, 
a quantum torus ${\cal A}_{q, {\bf i}}$ is defined by the generators 
$A_i$, $i \in I$, and relations 
$
q^{-\lambda_{kl}}A_kA_l = q^{-\lambda_{lk}}A_lA_k, 
$ 
where matrix $\lambda_{ij}$ is the inverse to the one $\varepsilon_{ij}$. 
(The cluster matrix $b_{ij}$ from \cite{BZq} is the transposed to our $\varepsilon_{ij}$). 
These quantum tori are glued into a quantum space 
${\cal A}_q$ via the rule given in formula 4.23
 in \cite{BZq}. 

\begin{lemma} \label{9.7.05.10}
Given a feed $\mathbf i$, there is an injective map of noncommutative algebras 
\begin{equation} \label{8.26.05.3}
p^*_{\bf i}: \Z[{\cal X}_{q, {\bf i}}] \lra \Z[{\cal A}_{q, {\bf i}}], 
\qquad  p^*_{\bf i}X_i = 
q^{-\sum_{k<l}\lambda_{kl}\varepsilon_{ik}\varepsilon_{il}}A^{\varepsilon_{ik}}_k
\end{equation}
commuting with mutations. 
\end{lemma}

{\bf Proof}. It is an easy calculation using the definitions of our 
Section 4 and Chapter 4 of \cite{BZq}, which we left to the reader. 

So we get a map of quantum spaces $p: {\cal A}_q \lra {\cal X}_q$. 
\vskip 3mm
\paragraph{Relations between different types of cluster transformations}
A feed cluster transformation ${\bf c}$ 
 gives rise to the corresponding cluster transformations 
 in several set-ups: the classical ${\cal X}$- and ${\cal A}$-cluster transformation 
${\bf c}^x$ and  ${\bf c}^a$, and 
the quantum ${\cal X}$- and ${\cal A}$-cluster transformation ${\bf c}_q^x$ and ${\bf c}_q^a$. 
Quantum mutations, 
and hence cluster transformations, admit the $q=1$ specialization. 
So   ${\bf c}_q^x={\rm Id}$ implies  ${\bf c}^x={\rm Id}$, and similalrly in the $a$-version.

\begin{lemma} \label{8.24.05.1} Assume that 
${\rm det}\varepsilon_{ij} \not = 0$. 
Let ${\bf c}$ be a feed 
cluster transformation. 
Then  
${\bf c}^a = {\rm Id}$ implies  ${\bf c}_q^x = {\rm Id}$.
\end{lemma}

{\bf Proof}. By \cite{BZq}, Theorem 6.1, ${\bf c}^a = {\rm Id}$ implies  ${\bf c}_q^a = {\rm Id}$. Since ${\rm det}\varepsilon_{ij} \not = 0$, 
 the algebra map  (\ref{8.26.05.3})
is injective. Since  ${\bf c}_q^a = {\rm Id}$, this implies the claim.

\section{The quantum double}

\subsection{Mutation maps for the quantum double}

Recall the {quantum torus algebra} (\ref{QTAL}) provided by 
a lattice $\Lambda$ with a bilinear skew-symmetric form $(\ast, \ast) \in \frac{1}{N}\Z$. 
Morphisms between the quantum torus algebras related to the 
lattices $\Lambda_1$ and $\Lambda_2$ with skew-symmetric bilinear forms are in bijection 
with 
homomorphisms of lattices  $\Lambda_2 \to \Lambda_1$ respecting the forms. 

%\vskip 3mm
%{\bf Remark}. Relations (\ref{QTAL}) in the quasiclassical limit give a  Poisson 
%bracket on ${\rm Hom}(\Lambda, {\Bbb G}_m)$ which is 
% twice the  one provided by  $(\ast, \ast)$. 

\vskip 3mm
Let ${\mathbf i}$ be a feed. Set $q_i:= q^{1/d_i}$. The lattice 
$\Lambda_{\cal D}$ with the form $(\ast, \ast)_{\cal D}$ 
related to the feed ${\mathbf i}$ gives rise to the quantum torus algebra 
${\bf D}_{\bf i}$. A cluster basis $(e_i, f_i)$, $i \in I$, of the latter 
provides the generators 
$B_i, X_i$, satisfying   
the relations 
\begin{equation} \label{4.28.03.11x}
q_i^{-1} X_iB_i  = 
q_i B_iX_i, \qquad  B_i X_j = 
X_jB_i \quad \mbox{if $i \not = j$}, \qquad  
q^{-\widehat \varepsilon_{ij}} X_i X_j = 
q^{-\widehat \varepsilon_{ji}} X_jX_i.
\end{equation}

\vskip 3mm
Denote by $ {\Bbb D}_{\mathbf i} = {\Bbb D}^q_{\mathbf i} $ the (non-commutative) fraction field of ${\bf D}_{\mathbf i}$. 
Let $\mu_k: {\mathbf i} \to {\mathbf i}'$ 
be a mutation. 
Our goal is to define 
a {\it quantum mutation} map,  understood as a homomorphism of 
non-commutative fields
$
\mu^q_{k}: 
{\Bbb D}_{{\mathbf i}'} \lra {\Bbb D}_{\mathbf i}. 
$ 
Recall the quantum dilogarithm 
$
{\bf \Psi}^q(x)= 
\prod_{k=1}^{\infty} (1+q^{2k-1}x)^{-1}. 
$ 
The generators $X_i$ and $\widetilde X_j$ 
commute. Thus
${\bf \Psi}^q(X_k)$ commutes with ${\bf \Psi}^q(\widetilde X_k)$. 
\begin{definition} \label{11.18.06.111} 
(i) The automorphism $\mu_k^{\sharp}: {\Bbb D}_{\bf i} \lra {\Bbb D}_{\bf i}$ 
is the conjugation by 
${\bf \Psi}^{q_k}(X_k)/{\bf \Psi}^{q_k}(\widetilde X_k)$: 
$$
\mu_k^{\sharp}:= 
{\rm Ad}_{{\bf \Psi}^{q_k}(X_k)/{\bf \Psi}^{q_k}(\widetilde X_k)}.
$$ 

(ii) The isomorphism 
$\mu'_k: {\bf D}_{\bf i'} \lra {\bf D}_{\bf i}$ 
is induced by the lattice map  $t'_k: {\Lambda}_{{\cal D}, \bf i'} \lra {\Lambda}_{{\cal D}, \bf i}$. 

(iii) The mutation map $\mu^q_k: {\Bbb D}_{\bf i'} \lra {\Bbb D}_{\bf i}$ 
is the composition 
$
\mu^q_k = \mu^{\sharp}_k\circ \mu'_k. 
$\end{definition} 

{\bf Remarks}. 1. The map $\mu_k^{\sharp}$ acts on the $X$-coordinates 
conjugating them by 
${\bf \Psi}^{q_k}(X_k)$. 

2.  The map 
$\mu'_k: {\bf D}_{\bf i'} \lra {\bf D}_{\bf i}$ 
is a homomorphism of algebras 
thanks to  Corollary \ref{Cor}.

3. The definition of the automorphism
$\mu_k^{\sharp}$ via adjoint action of
the ratio of quantum dilogarithms {\it a priori} produces
a formal power series in the generators. However Lemma \ref{11.18.06.3} below 
shows that the range of
the automorphism $\mu_k^{\sharp}$ is as already
anticipated in Definition \ref{11.18.06.111}.

\vskip 3mm
Here is an explicit computation of the automorphism $\mu_k^{\sharp}$. 
It is a bit involved. However what is really used is 
not this formula but transparent 
 Definition \ref{11.18.06.111}.  
%or its version (\ref{hghghgw1}) - (\ref{hghghgw}) below. 

\begin{lemma} \label{11.18.06.3} The automorphism $\mu_k^{\sharp}$ is 
given on the generators by the formulas
\begin{equation} \label{11.18.06.1}
B_i \lms B^{\sharp}_{i} := \left\{\begin{array}{lll} B_i& \mbox{ if } & i\not =k, \\
    B_k(1+q_kX_k)(1+q_k\widetilde X_k)^{-1} & \mbox{ if } &  i= k. \\
\end{array} \right.
\end{equation}
\begin{equation} \label{f3*}
X_i \lms X^{\sharp}_{i} := \left\{\begin{array}{lll} %X_k& \mbox{ if } & i=k, \\
X_i (1+q_kX_k)(1+q^3_kX_k)
  \ldots (1+q_k^{2|\varepsilon_{ik}|-1}X_k)& 
\mbox{ if } &  \varepsilon_{ik}\leq 0, 
\\
    X_i \left((1+q^{-1}_kX_k)(1+q_k^{-3}X_k)
  \ldots (1+q_k^{1-2|\varepsilon_{ik}|}X_k)\right)^{-1}& \mbox{ if } &  \varepsilon_{ik}\geq 0. \\
\end{array} \right.
\end{equation} 
\end{lemma}

{\bf Proof}. 
For any formal power series $\varphi(x)$ the relation 
$
q^{- \widehat \varepsilon_{ki}}X_kX_i = q^{- \widehat \varepsilon_{ik}}X_iX_k
$ 
implies 
\begin{equation} \label{hghghgw3} 
\varphi(X_k)X_i = 
X_i\varphi(q^{-2\widehat \varepsilon_{ik}}X_k).
\end{equation}

Recall the difference equation characterising  
${\bf \Psi}^q(x)$ 
up to a constant:
\begin{equation} \label{11.19.06.20}
{\bf \Psi}^q(q^2x) = (1+qx){\bf \Psi}^q(x), \quad 
\mbox{or, equivalently,}\quad  
{\bf \Psi}^q(q^{-2}x) = (1+q^{-1}x)^{-1}{\bf \Psi}^q(x).
\end{equation} 
It implies that  
formulas (\ref{11.18.06.1})  and  (\ref{f3*}) can be rewritten as
\begin{equation} \label{hghghgw1} 
B^\sharp_i =
B_k \cdot {\bf \Psi}^{q_k}(q_k^{2\delta_{ik}}X_k)
{\bf \Psi}^{q_k}(X_k)^{-1}\cdot
{\bf \Psi}^{q_k}(q_k^{-2\delta_{ik}}\widetilde X_k)^{-1}
{\bf \Psi}^{q_k}(\widetilde X_k),  
\end{equation}
\begin{equation} \label{hghghgw} 
X^\sharp_i =
X_i \cdot {\bf \Psi}^{q_k}(q_k^{-2\varepsilon_{ik}}X_k)
{\bf \Psi}^{q_k}(X_k)^{-1} .
\end{equation}
Using   (\ref{hghghgw3}) and 
$q_k^{-2\varepsilon_{ik}} =q^{-2\widehat \varepsilon_{ik}}$, 
we get 
$$
{\bf \Psi}^{q_k}(X_k)X_i{\bf \Psi}^{q_k}(X_k)^{-1} =
X_i {\bf \Psi}^{q_k}(q^{-2\widehat \varepsilon_{ik}}X_k) 
{\bf \Psi}^{q_k}(X_k)^{-1}  \stackrel{(\ref{hghghgw})}{=} X^\sharp_i.
$$ 
 Formula (\ref{hghghgw1}) is proved similarly. 
The lemma is proved.

\subsection{The quantum double and its properties}

Recall that quantum spaces are functors 
from the modular groupoid $\widehat {\cal G}^o$ 
(the opposite to $\widehat {\cal G}$) to the category 
${\rm QTor}^*$ (Section 2.1). Maps between quantum spaces are monomial 
morphisms of the functors. 

To relate with a geometric language, recall that 
the category of affine schemes is dual to 
the category of commutative algebras: 
A map ${\cal Y}_1 \to {\cal Y}_2$ 
of affine schemes is the same as  
a map of the corresponding commutative algebras ${Y}_2 \to {Y}_1$, where 
$Y_i$ is the algebra corresponding to the space  ${\cal Y}_i$. 

To define a quantum space ${\cal D}_q$ we use gluing 
isomorphisms $\mu^q_k: {\Bbb D}_{\bf i'} \to {\Bbb D}_{\bf i}$, 
and show that they satisfy relations (\ref{K10}). 
We can talk about them geometrically, 
saying that gluing isomorphisms correspond
to birational maps of non-commutative space 
$ {\cal D}_{{\bf i}, q} \to {\cal D}_{{\bf i'}, q}$, and that the 
quantum scheme ${\cal D}_q$ is obtained by gluing them via these maps. 
However the only meaning we put into this is that there is a  functor 
$\widehat {\cal G}^o \to {\rm QTor}^*$.  

\vskip 3mm

The main properties of the quantum double are summarized in the following theorem: 
\begin{theorem} \label{8.15.05.10gg} $a$) There is a 
$\widehat 
\Gamma$-equivariant quantum 
space ${\cal D}_q$. It is equipped with an involution 
$\ast_{\R}: {\cal D}_{q}\to {\cal D}_{q^{-1}}^{\rm op}$,  
given in any cluster coordinate system by 
$$
\ast_{\R}(q) = q^{-1}, \quad \ast_{\R}(X_i) = X_i, \quad \ast_{\R}(B_i) = B_i.
$$
\vskip 3mm
$a'$) There is a 
$\widehat 
\Gamma$-equivariant quantum 
space ${\cal X}_q$. It is equipped with an involution 
$\ast_{\R}: {\cal X}_{q}\to {\cal X}_{q^{-1}}^{\rm op}$,  
given in any cluster coordinate system by 
$$
\ast_{\R}(q) = q^{-1}, \quad \ast_{\R}(X_i) = X_i.
$$
\vskip 3mm
$b$) There is a map of quantum spaces $\pi:{\cal D}_q  \lra {\cal X}_q\times {\cal X}^{\rm op}_q $ 
given in a cluster coordinate system by 
\begin{equation} \label{GALG}
\pi^*(X_i \otimes 1) = X_i, \quad \pi^*(1 \otimes X_i) = 
\widetilde X_i:= X_i\prod_{j}B_j^{\varepsilon_{ij}}.
\end{equation}

\vskip 3mm
$c$) There is an involutive 
isomorphism of quantum spaces $i: {\cal D}_q\to {\cal D}_{q}^{\rm op}$ 
interchanging the two components of the projection $\pi$,   
given in any cluster 
coordinate system by 
$$
i^*B_i = B^{-1}_i, \quad i^*X_i = \widetilde X_i. 
$$

When $q=1$, the  unitary part ${\cal D}_U$ of the symplectic double is given by 
$
\overline B_j = i^*(B_j)$,
 $\overline X_j = i^* (X_j)$.

\vskip 3mm
$d$) There is a canonical map of 
quantum spaces $\theta_q: {\cal X}_q \to H_{\cal X}$. 
\end{theorem}
\vskip 3mm

{\bf Proof}. 
$a$), $a'$). Conjugation is an automorphism, so this and 
Remark 2 in Section 3.1 
imply that the mutation map $\mu^q_k: 
 {\Bbb D}_{\bf i'} \lra {\Bbb D}_{\bf i} 
$ is a $\ast_\R$-algebra homomorphism. Similarly 
the mutation map $\mu^q_k: 
 {\Bbb X}_{\bf i'} \lra {\Bbb X}_{\bf i}$ is a $\ast_\R$-algebra homomorphism. 
So to prove that 
we get a  $\Gamma$-equivariant quantum 
spaces ${\cal D}_q$ and ${\cal X}_q$ it remains to prove the following Lemma.

\begin{lemma}\label{11.6.06.1}
Cluster transformation  maps for the quantum double 
satisfy relations (\ref{K10}).
\end{lemma}

{\bf Proof}. 
%One needs only to check this for the $B$-coordinates, 
%since for the $X$-coordinates this was done in \cite{FG2}. 
This can be done by explicit calculations. 
Here is a trick which allows to skip involved computations. 
Let us do first the commutative case. 
One can embed a feed ${\bf i} = (I, \varepsilon_{ij}, d_i)$ to a 
bigger feed $\widetilde {\bf i} = 
(\widetilde I, \widetilde \varepsilon_{ij}, \widetilde d_i)$ so that 
$I \subset \widetilde I$, and ${\rm det}\widetilde \varepsilon_{ij}=1$. We put tilde's over the spaces related to the feed 
$\widetilde {\bf i}$. Then the map $\varphi: \widetilde {\cal A}
\times \widetilde {\cal A} \to \widetilde {\cal D}$ is an isomorphism. 
We claim that it is enough to check the relations (\ref{K10}) for the space 
$\widetilde {\cal A}
\times \widetilde {\cal A}$. Indeed, if we put then 
$B_s =1$ for every $s \in \widetilde I-I$ we recover the mutation for the feed 
${\bf i}$.  Since for the  ${\cal A}$-space  
these identities are known, and much easier to check, we are done. 

In the quantum case 
we use the quantum cluster algebras of \cite{BZq}. 
If ${\rm det}(\varepsilon_{ij})= 1$,  
there is a canonical 
non-commutative $q$-deformation of the cluster algebra. Namely, 
given a feed ${\mathbf i}$, 
a quantum torus ${\cal A}_{{\bf i}, q}$ 
corresponds to the lattice $\Lambda_{\cal X}^*$ with the bilinear form 
dual to the one on $\Lambda_{\cal X}$. 
These quantum tori are glued into a quantum space 
${\cal A}_q$ via the rule given in formula 4.23
 in {\it loc. cit.}. 
Since ${\rm det}(\varepsilon_{ij})= 1$, 
there is a canonical isomorphism of quantum spaces $p: 
{\cal A}_q \lra {\cal X}_q$ (\cite{King} or   
Appendix in Section 2, which for every feed $\mathbf i$ 
is induced by the lattice isomorphism  $p^*$ from Section 2.3. 
%It amounts to  an isomorphism of algebras 
%\begin{equation} \label{8.26.05.3}
%p^*_{\bf i}: \Z[{\cal X}_{{\bf i}, q}] \lra \Z[{\cal A}_{{\bf i}, q}], 
%\qquad  p^*_{\bf i}X_i = 
%q^{-\sum_{k<l}\lambda_{kl}\varepsilon_{ik}\varepsilon_{il}}A^{\varepsilon_{ik}}_k.
%\end{equation} 
We embed  a feed ${\bf i}$ to a feed $\widetilde {\bf i}$ as above. 
Then there is 
an isomorphism  $\widetilde {\cal A}_q \times \widetilde {\cal A}^{\rm op}_q
\stackrel{\sim}{\lra} \widetilde {\cal D}_q$, given 
by a composition 
$\widetilde {\cal A}_q \times \widetilde {\cal A}^{\rm op}_q
\stackrel{\sim}{\to} {\cal X}_q\times {\cal X}^{\rm op}_q 
\stackrel{\sim}{\to} {\cal D}_q$. 
Thus we can reduce 
checking the quantum relations (\ref{K10}) to the case of rank $2$ 
quantum cluster algebras, where they are known \cite{BZq}. 
The lemma is proved. 

\vskip 3mm
$b$). For a feed ${\bf i}$, thanks to  Lemma \ref{LR}, 
the lattice map $\pi^*$ in Section 2.3 (ii) provides a homomorphism of algebras 
$\pi^*_{\bf i}: {\bf X}_{\bf i} \otimes {\bf X}_{\bf i}^{\rm op} \lra {\bf D}_{\bf i}$.   
It is given by  formulas (\ref{GALG}).

Let us show that the quantum map $\pi$ commutes with each 
component of the mutation map. 

For $\mu_k^{\sharp}$ this is clear: 
Since $X_i$ and $\widetilde X_j$ commute, 
conjugation by ${\bf \Psi}^{q_k}(X_k)/{\bf \Psi}^{q_k}(\widetilde X_k)$
acts on $\pi^*(X_i\otimes 1)$ as  
${\rm Ad}_{{\bf \Psi}^{q_k}(X_k)}$, and on $\pi^*(1\otimes X_i)$ as 
${\rm Ad}^{-1}_{{\bf \Psi}^{q_k}(\widetilde X_k)}$. 

The maps $\mu_k'$ and $\pi$ are  monomial by definition, and  Poisson  
by Corollary \ref{Cor} and 
Lemma \ref{LR}. On the classical level 
the map $\pi$ intertwines the monomial parts 
of the mutation maps for the spaces ${\cal D}$ and ${\cal X}$ since the map of lattices 
$\pi^*$ in Section 2.3 (ii) is defined without 
using a basis. 
This implies the claim in the quantum case.

\vskip 3mm
$c$).  
Since $i^*$ is  a monomial map, and its classical version is Poisson by Section 2.3 (iv), 
the quantum map 
$i^*: {\bf D}^{\rm op}_{\bf i} \to {\bf D}_{\bf i}$ is an involutive automorphism. 
It commutes with the automorphism $\mu_k^{\sharp}$ since
$$
i^*\circ {\rm Ad}_{{\Psi}^{q_k}(X_k) /\Psi^{q_k}(\widetilde X_k) } \circ (i^*)^{-1} = 
{\rm Ad}^{-1}_{{\Psi}^{q_k}(\widetilde X_k)/  \Psi^{q_k}(X_k)}.  
$$
The claim that it commutes with the monomial part of the mutation 
follows from the fact that the map 
of lattices $i^*$  in Section 2.3 (iv) is defined without 
using a basis. 

$d$) For a given feed ${\bf i}$ the map of tori ${\cal X}_{{\bf i}, q} \to H_{\cal X}$ 
is dual to the canonical embedding ${\rm Ker} ~ p^* \hra \Lambda_{\cal X}$. 
Since  ${\rm Ker} ~ p^*$ is the kernel of the form $(\ast, \ast)_{\cal X}$, 
the automorphism part of the mutation map acts trivially on $H_{\cal X}$. The 
part d) is proved. The Theorem is proved.

\vskip 3mm
{\bf Remark}. The  obtained results include a construction 
of the quantum cluster ${\cal X}$-variety. It is  independent 
of the one given in \cite{FG2}. It is  
simpler and more transparent thanks to the decomposition of mutations and 
the fact that the automorphism part of mutation is given by conjugation by $\Psi^{q_k}(X_k)$. 

%\begin{lemma} \label{HTORUS}
%There is a canonical map of quantum spaces $\theta_q: {\cal X}_q \to H_{\cal X}$. 
%\end{lemma}

\vskip 3mm
{\it Connections between quantum ${\cal X}$-varieties.} 
 There are three ways to alter the space ${\cal X}_{q}$: 

(i) change $q$ to $q^{-1}$, 

(ii) change the quantum space ${\cal X}_{ q}$ to its chiral dual ${\cal X}^o_{ q}$, 

(iii) change the quantum space ${\cal X}_{ q}$ to the opposite quantum space 
${\cal X}^{\rm op}_{ q}$.

\noindent

\noindent
The resulting three quantum  spaces are canonically isomorphic (\cite{FG2II}, Lemma 2.1):
\begin{lemma} \label{6.9.03.11} 
There are canonical isomorphisms of quantum spaces
$$
\alpha_{\cal X}^q: {\cal X}_{q} \lra {\cal X}^{\rm op}_{ q^{-1}}, \qquad (\alpha_{\cal X}^q)^*
: X_i \lms X_i. 
$$
$$
i^q_{\cal X}: {\cal X}_{ q} \lra {\cal X}^o_{ {q^{-1}}}, \qquad (i^q_{\cal X})^*: X_i^o\lms X_i^{-1}.
$$
$$
\beta_{\cal X}^q:= \alpha_{\cal X}^q\circ i^q_{\cal X}: 
{\cal X}^o_{ q} \lra {\cal X}^{\rm op}_{ q}, \qquad X_i \lms {X_i^o}^{-1}. 
$$
\end{lemma}

 There are similar three ways to alter the space ${\cal D}_{q}$, and three similar 
isomorphisms, acting the same way on the $X$-coordinates, 
and identically on the $B$-coordinates. They are compatible with the projection 
$\pi$, while the involution $i$ from Theorem \ref{8.15.05.10} 
interchanges the two components of  $\pi$.

\subsection{Proof of Theorem  \ref{8.15.05.10}}

a), c), d), f). The claim in the part a) regarding the ${\cal A}$-space is known from the basic properties of the cluster algebras, see
\cite{FZI}. Its only non-trivial part is a proof of the 
$(h+2)$-gon relations in the  ${\cal A}$-case. The rest of a) plus 
c), d), f) follow from the similar quantum properties proved in  
Theorem \ref{8.15.05.10gg}.

b) 
For a feed ${\bf i}$, the projection $\varphi$ 
is given by a homomorphism $\varphi_{\bf i}: {\cal A}_{\bf i}\times 
{\cal A}_{\bf i} \lra {\cal D}_{\bf i}$, determined by the lattice map $\varphi^*$. 
We have to show that for a mutation $\mu_k: {\bf i} \to {\bf i'}$ the following diagram is commutative
$$
\begin{array}{ccc}
{\cal A}_{\bf i}\times {\cal A}_{\bf i}& \stackrel{\mu_k}{\lra} &
{\cal A}_{\bf i'}\times {\cal A}_{\bf {i'}}\\
%&&\\
\varphi_{\bf i}\downarrow &&\downarrow \varphi_{\bf i'}\\
%&&\\
{\cal D}_{\bf i}&\stackrel{\mu_k}{\lra} &{\cal D}_{\bf i'}
\end{array}
$$ 
We need to compare pull backs of the 
coordinates $B_i', X_i'$ obtained in two different ways. 
The claim for the coordinates $X_i'$ 
follows from the existence of the map $p: {\cal A} \to {\cal X}$. 
For the $B'$-coordinates 
it is non-obvious only for the coordinate $B_k'$, and  the map  $\mu_k^\sharp$. 
Going up and left 
in the diagram  we get
$$
B_k' \lra \frac{A^{'o}_k}{A^{'}_k}\lra \frac{A^o_k}{A_k}
\frac{(1+ {\Bbb A}_k^{+}/{\Bbb A}_k^{-})}{(1+ {\Bbb A}_k^{o+}/{\Bbb A}_k^{o-})}, \qquad 
{\Bbb A}_k^{o \pm} = \prod_{i}(A^o_i)^{[\pm \varepsilon_{ki}]_+}. 
$$
Going  left and up 
 in the diagram   we get the same:
$$
B_k' \lra B_k\frac{1+X_k}{1+\widetilde X_k} \lra 
\frac{A^o_k}{A_k}\frac{(1+ {\Bbb A}_k^{+}/{\Bbb A}_k^{-})}{(1+ {\Bbb A}_k^{o+}/{\Bbb A}_k^{o-})}.
$$
The claim that $\varphi^{*}\pi^{*} = p^* \times p^*$ 
and the claim about the $2$-forms follows from Proposition \ref{LR}(i).

e) Write $\pi = (\pi_-, \pi_+)$. Setting $B_i =1$ for all $i\in I$ we get subvariety of 
${\cal D}$ identified with ${\cal X}$ by $\pi_-$. The  restriction of the $2$-form $\Omega$ to 
this subvariety is clearly zero. 

g) It is straightforward linear algebra statement valid in every feed torus. 
We leave it to the reader. The Theorem is proved.

\vskip 3mm

\subsection{Duality conjectures for the double} 

The duality conjectures from \cite{FG2}, 
Section 4 should have the following version related to the double.

The semiring ${\Bbb L}_+({\cal D})$ of regular positive 
functions on the positive scheme ${\cal D}$ consists of all 
rational functions which are Laurent polynomials with 
positive integral coefficients in every cluster coordinate system on 
${\cal D}$. 
We conjecture that there exists a {\it canonical basis}
  in the space of regular 
positive functions on the Langlands dual space ${\cal D}^{\vee}$, 
parametrised by the set ${\cal D}(\Z^t)$ of the integral tropical 
points of ${\cal D}$. 
Moreover, it should admit  
a $q$-deformation to a basis in the space ${\Bbb L}_+({\cal D}_q)$ 
of regular positive 
functions on the space ${\cal D}_q$. The latter, 
 by definition, consists of all 
(non-commutative) rational functions which 
in every cluster coordinate system  are Laurent polynomials with 
positive integral coefficients in the cluster variables and $q$. 

The canonical basis should be compatible with  the canonical bases 
on the cluster ${\cal A}$- and ${\cal X}$-varieties 
from Section 4 of \cite{FG2}. 
Here is a precise conjecture. Recall that for a set $S$ we 
denote by $\Z_+\{S\}$ the semigroup of all $\Z_+$-valued 
function on the set $S$ with finite support.

\begin{conjecture} \label{11.7.06.1}
a) There exists a canonical isomorphism 
$ %\begin{equation} \label{11.7.06.2}
{\Bbb I}_{\cal D}: \Z_+\{{\cal D}(\Z^t)\} \stackrel{\sim}{\lra} 
{\Bbb L}_+({\cal D}^{\vee}),
$ % \end{equation}
which fits into the following commutative diagram: 
$$
\begin{array}{ccc}
\Z_+\{({\cal A}\times {\cal A})(\Z^t)\}&\lra & {\Bbb L}_+\Bigl({\cal X}^{\vee}\times 
{\cal X}^{\vee}\Bigl)\\
\downarrow &&\downarrow \\
\Z_+\{{\cal D}(\Z^t)\} & \stackrel{{\Bbb I}_{\cal D}}{\lra} 
&{\Bbb L}_+({\cal D}^{\vee})\\
\downarrow &&\downarrow \\
\Z_+\{({\cal X}\times {\cal X})(\Z^t)\} &\lra 
&{\Bbb L}_+\Bigl({\cal A}^{\vee}\times {{\cal A}}^{\vee}\Bigl)
\end{array}
$$
Here the vertical maps are induced by the maps in the diagram 
(\ref{11.5.06.1}). 
The top and the bottom horizontal maps are the isomorphisms predicted by the duality conjecture in {\it loc. cit}. 

b) The isomorphism ${\Bbb I}_{\cal D}$ admits  a $q$-deformation, given by an isomorphism
$% \begin{equation} \label{11.7.06.2}
{\Bbb I}^q_{\cal D}: {\cal D}(\Z^t) 
\stackrel{\sim}{\lra} {\Bbb L}_+({\cal D}_q^{\vee}).
%\end{equation}
$ \end{conjecture} 
The maps ${\Bbb I}_{\cal D}$ and ${\Bbb I}^q_{\cal D}$ should satisfy 
 additional properties just like the ones listed in {\it loc. cit}.

\section{The quantum dilogarithm and its properties}

%Recall the dilogarithm function
%$$
%{\rm Li}_2(x):= -\int_0^x\log(1-t)dt.
%$$
\subsection{The quantum logarithm function and its properties} 
We assume throughout this Section that $\hbar >0$. The quantum logarithm 
is the following function:
\begin{equation} \label{phi}
\phi^\hbar(z):= - 2\pi \hbar\int_{\Omega}\frac{e^{-ipz}}{(e^{\pi p} - 
e^{-\pi p})
(e^{\pi \hbar p}-e^{-\pi \hbar p}) }dp; 
\end{equation}
where the contour $\Omega$ goes along the real axes 
from $- \infty$ to $\infty$ bypassing the origin 
from above. 

\begin{proposition} \label{5.31.03.1} The function $\phi^{\hbar}(x)$ 
enjoys the following properties. 
$$ \lim_{\hbar
\rightarrow 0}\phi^\hbar(z) = \log(e^z + 1).\leqno{\bf (A1)} 
$$
$$
\phi^\hbar(z)-\phi^\hbar(-z)=z. \leqno{\bf (A2)}  
$$
$$
\overline{\phi^\hbar(z)} = \phi^\hbar(\overline{z}).\leqno{\bf (A3)}
$$
$$
\phi^\hbar(z)/\hbar = \phi^{1/\hbar}(z/\hbar).\leqno{\bf (A4)}
$$
$$
\phi^\hbar(z+i\pi \hbar)-\phi^\hbar(z-i\pi \hbar) = \frac{2\pi i
\hbar}{e^{-z}+1}, \qquad 
 \phi^\hbar(z+i\pi)-\phi^\hbar(z-i\pi)
= \frac{2\pi i}{e^{-z/\hbar}+1}.\leqno{\bf (A5)}
$$
$$ \phi^1(z) = \frac{z}{1-e^{-z}}.\leqno{\bf (A6)}
$$
{\bf (A7)} The form
$\phi^\hbar(z)dz$ is meromorphic with simple poles 
at the points of the upper half plane 
$$
\{\pi i \Bigl((2m-1)+
(2n-1)\hbar\Bigr)|m,n \in {\mathbb N}\} 
~~~\mbox{with residues $2 \pi i \hbar$},
$$
 and at the 
points of the lower half plane 
$$
\{-\pi i \Bigl((2m-1)+(2n-1)\hbar\Bigr)|m,n \in {\mathbb N}\} 
~~~\mbox{with residues 
$-2\pi i\hbar$}.
$$
$$
\phi^\hbar(z)= \frac{\hbar}{(\hbar +1)}\left(\phi^{\hbar+1}(z+i\pi) + 
\phi^\frac{\hbar+1}{\hbar}(z/\hbar-i\pi)\right) = \frac{\hbar}{(\hbar 
+1)}\left(\phi^{\hbar+1}(z-i\pi) + 
\phi^\frac{\hbar+1}{\hbar}(z/\hbar+i\pi)\right). \leqno{\bf (A8)}
$$
$$ \sum\limits_{l=\frac{1-r}{2}}^{\frac{r-1}{2}} 
\sum\limits_{m=\frac{1-s}{2}}^{\frac{s-1}{2}} \phi^\hbar(z+\frac{2\pi i}{r}l 
+ \frac{2\pi i\hbar }{s}m ) = 
s \phi^{\frac{r}{s}\hbar}(rz), \leqno{\bf (A9)}
$$ 
where the sum is taken over half-integers if the summation limits are 
half-integers.
\end{proposition}

{\bf Proof}. {\bf A1}. Since the contour 
$\Omega$ bypasses the origin from above, $\lim\limits_{z \rightarrow
-\infty}\phi^\hbar(z)=\lim\limits_{z \rightarrow
-\infty} \ln(e^z+1)=0$. Thus it is sufficient to prove that
$\lim\limits_{\hbar \rightarrow 0}\frac{\partial}{\partial z} \phi^\hbar(z) =
\frac{1}{e^{-z}+1}$. The l.h.s. of the latter is computed using
residues:
$$ \lim_{\hbar \rightarrow 0}\frac{\partial}{\partial z}
\phi^\hbar(z) = \lim_{\hbar \rightarrow 0}-\frac{\pi\hbar}{2}\int_\Omega
\frac{-ipe^{-ipz}}{{\rm sh}(\pi p){\rm sh}(\pi \hbar p)}dp =
$$
$$
=\frac{i}{2} \int_\Omega\frac{e^{-ipz}}{{\rm sh}(\pi p)}=
 \frac{i}{2} \frac{1}{1+e^{z}}(\int_{\Omega}
-\int_{\Omega+i})\frac{e^{-ipz}}{{\rm sh}(\pi p)}dp = \frac{-1}{1+e^{z}}{\mbox
Res}_{p=i } \frac{e^{-ipz}}{{\rm sh}(\pi p)} 
=\frac{1}{e^{-z}+1}.
$$

 {\bf A2}. It is verified by computing  using 
residues:
$$
\phi^\hbar(z)-\phi^\hbar(-z) 
=-\frac{\pi \hbar}{2}\int_\Omega\frac{e^{-ipz}-e^{ipz}}{{\rm sh}(\pi p){\rm sh}(\pi
\hbar p)}dp= -\frac{\pi
\hbar}{2}(\int_\Omega+\int_{-\Omega})\frac{e^{-ipz}}{{\rm sh}(\pi p){\rm sh}(\pi
\hbar p)}dp =
$$
$$
= \frac{\pi \hbar}{2} 2\pi i\, \mbox{Res}_{z=0}\frac{e^{-ipz}}{{\rm sh}(\pi
p){\rm sh}(\pi \hbar p)} = z.
$$

{\bf A3}. It is obtained by the change of the integration
variable $q=-\overline{p}$:

$$
\overline{\phi^\hbar(z)}=-\frac{\pi\hbar}{2}\int_{\overline{\Omega}}
\frac{e^{ip\overline{z}}}{{\rm sh}(\pi p){\rm sh}(\pi \hbar p)}dp =
$$
$$
= \frac{\pi\hbar}{2}\int_{-\Omega} \frac{e^{ip\overline{z}}}{{\rm sh}(\pi
p){\rm sh}(\pi \hbar p)}dp=-\frac{\pi\hbar}{2}\int_{\Omega}
\frac{e^{-ip\overline{z}}}{{\rm sh}(\pi p){\rm sh}(\pi \hbar p)}dp 
= \phi^\hbar(\overline{z}).
$$

{\bf A4}. It is  obtained by the change of the integration
variable $q = p/\hbar$.
%$$
%\phi^{1/\hbar}(z/\hbar) =
%-\frac{\pi}{2\hbar}\int_{\Omega}\frac{e^{-ipz/\hbar}}{{\rm sh}(\pi p){\rm sh}(\pi 
%p/\hbar)}dp 
%=-\frac{\pi}{2\hbar}\int_{\Omega}\frac{e^{-iqz}}{{\rm sh}(\pi \hbar q){\rm sh}(\pi 
%q)}d(q \hbar)=
%\phi^{\hbar}(z)/\hbar 
%$$

{\bf A5}  The proof is similar to the
proof of {\bf A1}. We give only the proof of the first identity:
$$
\phi^\hbar(z+i\pi\hbar)-\phi^\hbar(z-i\pi\hbar)
= -\frac{\pi\hbar}{2}\int_\Omega \frac{e^{-ipz}(e^{\pi\hbar p} - e^{-\pi\hbar
p})}{{\rm sh}(\pi p){\rm sh}(\pi \hbar p)}dp =
$$
$$
=-\pi\hbar\int_\Omega \frac{e^{-ipz}}{{\rm sh}(\pi p)}dp = \frac{-\pi 
\hbar}{e^z+1}
\mbox{~Res}_{p=i} \frac{e^{-ipz}}{{\rm sh}(\pi p)} = \frac{2 \pi i
\hbar}{e^{-z}+1}.
$$

{\bf A6}. It is done by explicit residue calculations:
$$ \phi^1(z) = -\frac{\pi}{2}\int_\Omega\frac{e^{-ipz}}{{\rm sh}^2(\pi p)}dp =
 \frac{\pi^2 i}{2}\frac{1}{1-e^{-z}}\mbox{Res}_{p=0} \frac{e^{-i p 
z}}{{\rm sh}^2(\pi p)} = \frac{z}{1-e^{-z}}. 
$$

{\bf A7}. The integral
(\ref{phi}) converges for $|\Im{z}| < \pi(1+\hbar)$. Using 
{\bf A5},  the function $\phi^h$ can be continued to the
whole complex plane. Then {\bf A7} is obvious. 
Property {\bf A.9} is an easy consequence of the identity
$$
\sum\limits_{l=\frac{1-p}{2}}^{\frac{p-1}{2}} e^{lx} = 
\frac{{\rm sh}((l+1/2)x)}{{\rm sh}((1/2)x)}.
$$

\vskip 3mm
\subsection{The quantum dilogarithm and its properties}
Recall the quantum dilogarithm function:
$$
\Phi^\hbar(z) := {\rm exp}\Bigl(-\frac{1}{4}\int_{\Omega}\frac{e^{-ipz}}{ {\rm sh} (\pi p)
{\rm sh} (\pi \hbar p) } \frac{dp}{p} \Bigr).
$$

\begin{proposition} \label{qdilog} The function $\Phi^{\hbar}(x)$
enjoys the following properties.
$$
2 \pi i \hbar\,d \log \Phi^\hbar(z) =  \phi^\hbar(z)dz. \leqno{\bf (B)}
$$
$$
\lim_{\Re z \rightarrow -\infty}\Phi^\hbar(z)=1.\leqno{\bf (B0)}
$$
Here the limit is taken along a line parallel to the real
axis. 
$$ \lim_{\hbar \rightarrow 0}\Phi^\hbar(z)/\exp \frac{L_2(e^{z})}{2\pi i
  \hbar} = 1. \leqno{\bf (B1)}
$$
$$
\Phi^\hbar(z)\Phi^\hbar(-z)=\exp\left(\frac{z^2}{4\pi i \hbar}\right)
e^{-\frac{\pi i}{12}(\hbar+\hbar^{-1})}. \leqno{\bf (B2)}
$$
$$
\overline{\Phi^\hbar(z)} = (\Phi^\hbar(\overline{z}))^{-1}.\mbox{ In
  particular } |\Phi^\hbar(z)|=1 \mbox{ for } z\in {\mathbb R}.\leqno{\bf (B3)}
$$
$$
\Phi^\hbar(z) = \Phi^{1/\hbar}(z/\hbar).\leqno{\bf (B4)}
$$
$$
\Phi^\hbar(z+ 2 \pi i \hbar) = \Phi^\hbar(z) (1+qe^z), \qquad
\Phi^\hbar(z+ 2 \pi i ) = \Phi^\hbar(z) (1+ q^{\vee}e^{z/\hbar}).  \leqno{\bf (B5)}
$$
$$ \Phi^1(z) = e^{(\pi^2/6-Li_2(1-e^z))/2\pi i}. \leqno{\bf (B6)}
$$
where the r.h.s. should be understood as the analytic continuation
from the origin.

\noindent({\bf B7}) The function
$\Phi^\hbar(z)dz$ is meromorphic. Its poles are simple poles 
located at the upper half plane at the points 
$$
\{-\pi i \Bigl((2m-1) - 
(2n-1)\hbar)\Bigr)~|~m,n \in {\mathbb N}\},
$$  
and its zeroes are located at the lower half plane at the
points $$\{\pi i \Bigl((2m-1)+(2n-1)\hbar\Bigr) ~| ~m,n \in {\mathbb N}\}.$$
$$
\Phi^\hbar(z)= \Phi^{\hbar+1}(z+i\pi)
\Phi^\frac{\hbar+1}{\hbar}(z/\hbar-i\pi)= \Phi^{\hbar+1}(z-i\pi)
\Phi^\frac{\hbar+1}{\hbar}(z/\hbar+i\pi). \leqno{\bf (B8)}
$$
$$ \prod\limits_{l=\frac{1-r}{2}}^{\frac{r-1}{2}}
\prod\limits_{m=\frac{1-s}{2}}^{\frac{s-1}{2}} \Phi^\hbar(z+\frac{2\pi i}{r}l
+ \frac{2\pi i\hbar }{s}m ) =
 \Phi^{\frac{r}{s}\hbar}(rz), \leqno{\bf (B9)}
$$
where the sum is taken over half-integers if the summation limits are
half-integers.
\end{proposition}

{\bf Proof.} Property {\bf B0} follows from the Cauchy
lemma for ${\rm Im} z \leq 1+\hbar$ and can be extended to any line parallel to
the real line using the relation {\bf B5}. Properties {\bf B3}-{\bf B5} 
and {\bf B7}-{\bf B9} follow immediately from the corresponding properties of 
the function $\phi^\hbar$ and property {\bf B0}. 
Indeed, since their logarithmic 
derivatives give exactly the relations for the functions $\phi$ one needs only 
to verify that they are valid in the limit $z \rightarrow -\infty$, what is 
obvious in these cases.

The proof of Property {\bf B1} is also analogous to the proof of
{\bf A1}. According to the latter we have:
$$
\int_\Omega\frac{e^{-ipz}}{{\rm sh}(\pi p)}=\frac{-2i}{e^{-z}+1}.
$$
By integrating both sides twice one has
$$
\int_\Omega\frac{e^{-ipz}}{-p^2{\rm sh}(\pi p)}=2iLi_2(e^z).
$$
It implies that
$$
\lim_{\hbar \rightarrow 0}-\frac{1}{4}\int_{\Omega}\frac{e^{-ipz}}{ {\rm sh} (\pi p)
{\rm s\hbar} (\pi h p) } \frac{dp}{p} - \frac{Li_2(e^z)}{2 \pi i h} =0
$$
since this difference is an odd function of $\hbar$ and the poles of the both terms cancel each other.

Property {\bf B2} can be proved by residue computations:

$$\Phi^\hbar(z)\Phi^\hbar(-z)=\exp\Bigl(-\frac{1}{4}\int_{\Omega}\Bigl(\frac{e^{-ipz}}{
  {\rm sh} (\pi p) {\rm sh} (\pi \hbar p) } - \frac{e^{ipz}}{ {\rm sh} (\pi p)
{\rm sh} (\pi \hbar p) }\Bigr)\frac{dp}{p}\Bigr) =
$$
$$
=\exp\Bigl(\frac{\pi i}{2}{\rm Res}_{p=0} \frac{e^{-ipz}}{ {\rm sh} (\pi p)
{\rm sh} (\pi \hbar p) }\frac{dp}{p}\Bigr)=
$$
$$
=\exp\bigl(-\frac{1}{2\pi i \hbar}{\rm Res}_{p=0} \frac{e^{-ipz}}{p^3(1+(\pi p)^2/6)(1+(\pi \hbar p)^2/6)}dp\bigr)=
$$
$$
=\exp\Bigl(-\frac{1}{2\pi i \hbar}(-z^2/2-\pi^2/6-\pi^2\hbar^2/6)\Bigr)=\exp\left(\frac{z^2}{4\pi i \hbar}\right)e^{-\frac{\pi i}{12}(\hbar+\hbar^{-1})}.
$$

Property {\bf B6} can be proven by a direct calculation:
$$
\Phi^1(z) = \exp(\frac{1}{2\pi i} \int_{-\infty}^{z} \phi^1(t)dt) =
\exp\left( \frac{1}{2\pi i}\int_{-\infty}^{z}\frac{z d(1-e^z)}{1-e^z}\right)=
\exp\left(\frac{Li_2(1)-Li_2(1-e^z)}{2\pi i}\right).
$$

%$$ \lim_{h \rightarrow 0}\Phi^h(x)/\exp \frac{Li_2(-e^{-x})}{2\pi i
%  h} = 1, \qquad x \in \R \leqno{\bf B6}
%$$
%$$
%{\Phi^1(z)} = \exp \frac{Li_2(-e^{-z})}{2\pi i} WRONG?
 %\leqno{\bf B7}
%$$
%{\rm Asymptotics}. Let $z =x+iy$. Then one  has   
%$$
%\lim_{x \to - \infty}\Phi^h(z) = 1, \qquad 
%\Phi^h(z) \sim e^{\frac{xy}{2\pi h}}, \quad x \to + \infty \leqno{\bf B5} 
%$$ 
%\end{proposition}

\section{Representations of 
quantum cluster varieties}

\subsection{The intertwiner}

\paragraph{The logarithmic coordinates} 
These are the coordinates  $a_i:= \log A_i$ on the 
 set ${\cal A}_{\bf i}(\R_{>0})$ of positive real points of the 
feed torus ${\cal A}_{\bf i}$. 
Here is a relevant 
algebraic notion:  
{\it The logarithmic feed 
${\cal A}$-space} ${\cal A}_{\log, {\bf i}}$ is the affine space 
with the coordinates $a_i$. So there is a canonical isomorphism
$$
{\cal A}_{\log, \bf i}(\R) = 
{\cal A}_{\bf i}(\R_{>0}).  
$$
Below we use the notation ${\cal A}_{\bf i}^+$ for ${\cal A}_{\log, \bf i}(\R)$,
 and $L^2({\cal A}_{\bf i}^+)$
for $L^2({\cal A}_{\log, \bf i}(\R))$.

It is easy to see that
the volume form $da_1 \ldots da_n$ on ${\cal A}^+$ 
changes the sign under the mutations. Thus mutations provide isomorphisms 
of the Hilbert spaces 
$L^2({\cal A}^+_{\bf i})$ for different feeds. 

\paragraph{The set-up} Let us set
$$
\hbar_k:= \hbar d_k^{-1}, \quad q_k:= e^{\pi i \hbar_k}, 
\qquad \hbar^{\vee}:= 1/\hbar, \quad \hbar_k^{\vee}:= 1/\hbar_k, \qquad q_k^{\vee}:= e^{\pi i /\hbar_k} = 
e^{\pi i \hbar^{\vee}_k}.
$$
It is handy to introduce the following notation:
\begin{equation} \label{handy}
\alpha_k^+:= \sum_{j\in I} 
[\varepsilon_{kj}]_+a_j, \qquad 
\alpha_k^-:= \sum_{j\in I} 
[-\varepsilon_{kj}]_+a_j; \qquad \mbox{so $\alpha_k^+ - \alpha_k^- = 
\sum_{j\in I} \varepsilon_{kj}a_j$.}
\end{equation}
Consider the following differential 
operators in ${\cal A}_{{\log}, \bf i}$: 
\begin{equation} \label{5.26.08.1}
\widehat x_p:= 2\pi i \hbar_p \frac{\partial}{\partial a_p}
-\alpha_p^+, \qquad \widehat b_p = a_p.
\end{equation}
They satisfy the commutation relations of the Heisenberg $\ast$-algebra 
${\cal H}^{\hbar}_{\bf i}$: 
$$
[\widehat x_p, \widehat x_q] = 2 \pi i \hbar \widehat \varepsilon_{pq}, 
\qquad [\widehat x_p, \widehat b_q] = 2 \pi i \hbar_p\delta_{pq}, \qquad
\ast \widehat x_p = \widehat x_p, \quad \ast \widehat b_p = \widehat b_p.
$$ 
Furthermore,  consider another collection of first order differential 
operators in ${\cal A}_{{\log}, \bf i}$: 
\begin{equation} \label{5.26.08.2}
\widehat {{\widetilde x}_p} = \widehat x_p 
+ \sum_{q}\varepsilon_{pq}  \widehat b_q = 
 2\pi i \hbar_p\frac{\partial}{\partial a_p}
-\alpha_p^-.
\end{equation}
They commute with the operators (\ref{5.26.08.1}), 
and, togerther with the operators $\widehat b_p$,  
satisfy the commutation relations 
of the opposite Heisenberg $\ast$-algebra 
${\cal H}^{\rm op}_{\bf i}$:
%We get a $\ast$-representation :
%provide a $\ast$-representation  by unbounded operators in the Hilbert
% space $L^2({\cal A}_{\bf i}^+)$:
$$
[\widehat {{\widetilde x}_p}, \widehat {{{\widetilde x}_q} }] 
= -2 \pi i \hbar \widehat \varepsilon_{pq}, \qquad 
[\widehat {{\widetilde x}_p}, \widehat b_q]  = 2 \pi i \hbar_p\delta_{pq}, 
\qquad [\widehat {{\widetilde x}_p}, \widehat x_p] = 0.\qquad
\ast \widehat {{\widetilde x}_p} = \widehat {{\widetilde x}_p}. 
$$ 
The exponentials $\widehat X_p:= {\rm exp}(\widehat x_p)$ and $\widehat B_p:= {\rm exp}(\widehat b_p)$ are 
 difference operators: 
\begin{equation} \label{intertwiner1}
\widehat X_p f (A_1, ..., A_n) = ({\Bbb A}^{+}_p)^{-1} \cdot f(A_1, ..., q^2_pA_p, 
\ldots A_n), \quad 
\widehat B_pf= A_pf.
\end{equation}
They satisfy the commutation relations  
of the quantum torus $\ast$-algebra ${\bf D}^q_{\bf i}$.
 
%by unbounded operators in the Hilbert
% space $L^2({\cal A}_{\bf i}^+)$. 

There is a canonical isomorphism ${\cal H}^{\rm op}_{\bf i} \to {\cal H}_{\bf i^o}$ 
acting as the identity on the generators.

\vskip 3mm
Finally, let us introduce the Langlands modular dual collection of operators
\begin{equation} \label{7.26.07.1}
\widehat {x}^{\vee}_p := \widehat {x}_p/\hbar_p, \quad 
\widehat {b}^{\vee}_p := \widehat {b}_p/\hbar_p, \qquad 
\widehat X^{\vee}_p:= {\rm exp}(\widehat x^{\vee}_p), \quad 
\widehat B^{\vee}_p:= {\rm exp}(\widehat b^{\vee}_p).
\end{equation} 
The operators $\{\widehat {x}^{\vee}_p, \widehat {b}^{\vee}_p\}$  
satisfy the commutation relations of the Heisenberg $\ast$-algebra 
${\cal H}^{\hbar^{\vee}}_{\bf i^{\vee}}$. 
The operators $\{\widehat {X}^{\vee}_p, \widehat {B}^{\vee}_p\}$  
satisfy the commutation relations of the $\ast$-algebra  
${\bf D}^{q^{\vee}}_{\bf i^{\vee}}$.

\vskip 3mm
{\bf Remark}. One can say that  operators (\ref{5.26.08.1}) 
(respectively (\ref{5.26.08.2}) and (\ref{intertwiner1})) 
provide a $\ast$-representation 
by unbounded operators of the Heisenberg $\ast$-algebra 
${\cal H}_{\bf i}$ (respectively 
${\cal H}^{\rm op}_{\bf i}$ and ${\bf D}^q_{\bf i}$)  
in the Hilbert
space $L^2({\cal A}_{\bf i}^+)$.   
\vskip 3mm

Since the operators $\widehat x_k$ and 
$\widehat {{\widetilde x}_k}$ are self-adjoint in $L^2({\cal A}_{\bf i}^+)$, 
one can apply to them any continuous function on the real line. 
If the function  is unitary, i.e. 
takes the values at the unit circle, we get a unitary operator. 
We employ below the quantum dilogarithm function 
$\Phi^{\hbar_k}(x)$, which is unitary by the property ${\bf B}3$. 

\begin{definition} \label{intertwiner}
Given a mutation $\mu_k: {\bf i} \to {\bf i}'$, 
the  intertwining operator 
\begin{equation} \label{12.01.06.1}
{\bf K}_{{\bf i}', {\bf i}}: L^2({\cal A}_{\bf i'}^+) \lra 
L^2({\cal A}_{\bf i}^+)
\end{equation}
is defined as a composition \footnote{We suppress the indices when 
this does not lead to a confusion.} 
$
{\bf K}_{{\bf i'}, {\bf i}}:= {\bf K}^{\sharp}\circ 
{\bf K}^{'} 
$, where:
\begin{itemize}
\item
{\it The operator ${\bf K}^{\sharp}$} is  the 
 ratio  of quantum dilogarithms of the operators $\widehat x_k$ and $\widehat {\widetilde x}_k$: 
$$
{\bf K}^{\sharp}:=  \Phi^{\hbar_k}(\widehat x_k)
\Phi^{\hbar_k}( \widehat {{\widetilde x}_k})^{-1}. 
$$

\item
{\it The operator ${\bf K}^{'}$} is induced by the logarithmic version 
$\mu_{\log, k}': {\cal A}^+_{\bf i} \lra {\cal A}^+_{\bf i'}$ of the 
cluster ${\cal A}$-transformation 
$\mu_k'$, see (\ref{11.18.06.10sdf}), which  is a linear map 
acting on the coordinates as follows:
\begin{equation} \label{11.18.06.10d}
{a}_i'\lms \left\{\begin{array}{lll} 
a_i& \mbox{ if } & i\not =k, \\
    \alpha_k^- - a_k
 & \mbox{ if } &  i= k. \\
\end{array} \right.
\end{equation}
\end{itemize}
\end{definition}

Let us write the operator ${\bf K}^{\sharp}$ 
by using the Fourier transform ${\cal F}_{a_k}$
along the $a_k$-coordinate: 
$$
{\cal F}_{a_k}(f)(c) = \widehat{f}(c)=\int e^{a_kc/2\pi i \hbar}f(a_k)da_k;~~
f(a_k)=\frac{1}{(2\pi i)^2 \hbar}\int e^{-a_kc_k/2\pi i \hbar}\widehat{f}(c)dc.
$$
Here $c$ is the variable dual to $a_k$. 
We omit the variables $a_1,\ldots, \widehat a_{k}, \ldots,a_n$
in both $f$ and $\widehat{f}$. 
Then 
$$
-2\pi i \hbar \widehat{\frac{\partial f}{\partial
    a_k}}=c\widehat{f},~~\widehat{a_kf}=2\pi i \hbar\frac{\partial
  \widehat{f}}{\partial c}.
$$
Equivalently:
\begin{equation} \label{1.05.06.14}
a_k = {\cal F}^{-1}_{a_k}\circ (2\pi i\hbar \frac{\partial}{\partial c})\circ 
{\cal F}_{a_k}, \qquad 
-2\pi i \hbar \frac{\partial}{\partial a_k}  = 
{\cal F}^{-1}_{a_k}\circ c \circ 
{\cal F}_{a_k}.
\end{equation}
Then 
\begin{equation} \label{1.05.06.14a}
{\bf K}^{\sharp} = {\cal F}_{a_k}^{-1} \circ \Phi^{\hbar_k}\Bigl(-d^{-1}_k c -\alpha_k^+\Bigr)\Phi^{\hbar_k}\Bigl( -d^{-1}_k c -\alpha_k^-\Bigr)^{-1} \circ {\cal F}_{a_k}. 
\end{equation}

The inverse ${\bf K}^{-1}_{{\bf i'}, {\bf i}}$ of the intertwiner  has a bit simpler  presentation as an integral operator:
\begin{equation}\label{convolution}
({\bf K}^{-1}_{{\bf i}', {\bf i}}f)(a_1,\ldots, a'_k, \ldots, a_n) := \int G(a_1,\ldots, 
a'_k+a_k, \ldots, a_n)f(a_1,\ldots, a_k,\ldots, a_n)da_k,
\end{equation}
where 
\begin{equation}\label{11.13.03.1}
G(a_1,\ldots,a_n):= %{\rm exp}\Bigl( \frac{(\sum_{j\in I}[\varepsilon_{kj}]_+a_j)^2 }{4\pi i h_k}\Bigr)\times
\frac{1}{(2\pi i)^2 \hbar}\int \Phi^{\hbar_k}(-{d}^{-1}_kc
-\alpha_k^+)^{-1}\Phi^{\hbar_k}(-{d}^{-1}_kc
-\alpha_k^-)^{-1}\exp\left(c
\frac{a_k-\alpha_k^-}{2\pi i \hbar}\right)dc.
\end{equation}

The kernel of the integral operator ${\bf K}^{\sharp}$ 
depends on a choice of a  representation of the Heisenberg algebra ${\cal H}_{\bf i}$ in 
$L^2({\cal A}_{\bf i}^+)$. 
Here is a different realization, employed in \cite{FG2II}: 
\begin{equation}\label{anotherR}
\widehat x^{\rm old}_p:= \pi i \hbar_p \frac{\partial}{\partial a_p} 
-\sum_{q}\varepsilon_{pq}  a_q, \qquad \widehat b^{\rm old}_p = 2a_p.
\end{equation}
Its disadvantage is that the operators $\widehat 
B_p = {\rm exp}(\widehat b_p)$ are given by multiplication by $A_p^2$. 
The inverse ${\bf K}^{-1}$ to the composition ${\bf K}= {\bf K}^{\sharp}{\bf K}^{'}$ 
is the intertwiner  given by formulas (14)-(15) in {\it loc. cit.}.

\subsection{Main results} 

{\it A  space $W_{\bf i}$}. Let $W_{\bf i} \subset L^2({\cal A}_{\bf i}^+)$ be the space of finite 
linear combinations of functions 
\begin{equation} \label{9.4.07.1}
P{\rm exp}{(-\alpha \sum_{i \in I} (a_i^2/2+ b_ia_i))}, \quad \mbox{
where $\alpha >0$, $b_i\in \C$, and $P$ is a polynomial in $a_i$, $i\in I$}. 
\end{equation}
It is invariant under the Fourier transform. 
The operators $\widehat B_i$, $\widehat B^{\vee}_i$, 
 $\widehat X_i$, $\widehat X^{\vee}_i$ are symmetric unbounded operators 
 defined on $W_{\bf i}$. They preserve $W_{\bf i}$, 
and satisfy, on $W_{\bf i}$, the standard commutation relations. 

\vskip 3mm
{\it The $\ast$-algebra ${\bf L}$ and its Schwartz space.} Denote by ${\Bbb L}_q$ the space of universally Laurent polynomials for the 
cluster double ${\cal D}_q$, and 
 by ${\bf L}$ the algebra ${\Bbb L}_q \otimes {\Bbb L}_{q^{\vee}}$. 
The algebra ${\bf L}$ is $\ast$-invariant. 

\begin{definition} \label{7.10.07.1} The Schwartz space ${\cal S}_{{\bf i}} = {\cal S}_{{\bf L}, {\bf i}}$ 
is a subspace of $L^2({\cal A}^+_{\bf i})$ consisting
of vectors $f$ such that the functional $w \to (f, \widehat Aw)$ on $W_{\bf i}$
is continuous for the $L^2$-norm, for all $A \in {\bf L}$.
\end{definition}
Denote by $(*,*)$ the scalar product in $L^2({\cal A}_{\bf i}^+)$.
The Schwartz space ${\cal S}_{\bf i}$
is
 the common domain of definition of
operators from ${\bf L}$ in $L^2({\cal A}_{\bf i}^+)$. Indeed,
since $W_{\bf i}$ is dense in  $L^2({\cal A}_{\bf i}^+)$,
the Riesz theorem implies that for any $f \in {\cal S}_{\bf i}$ there exists a
unique $g \in L^2({\cal A}_{\bf i}^+)$ 
such that $(g,w) = (f, \widehat Aw)$. We set $\widehat {A}^\ast f := g$.
Equivalently,
let $W_{\bf i}^*$ be the algebraic linear dual to $W_{\bf i}$.
So $L^2({\cal A}_{\bf i}^+) \subset W_{\bf i}^*$. Then
$$
{\cal S}_{\bf i}= \{v \in W_{\bf i}^*~ | ~ \widehat A^*v \in L^2({\cal A}_{\bf i}^+) 
\quad \mbox{for any $A \in {\bf L}$}\} \cap
L^2({\cal A}_{\bf i}^+). 
$$

The Schwartz space ${\cal S}_{{\bf i}}$ has a natural topology given by  seminorms
$$
\rho_B(f):= ||Bf||_{L^2}, \qquad \mbox{$B$ runs through a basis in ${\bf L}$}.
$$

\begin{definition} \label{7.10.07.1a} The distribution space ${\cal S}^*_{{\bf i}}$ 
is the topological dual of the space ${\cal S}_{{\bf i}}$.
\end{definition}

{\it Intertwiners for cluster transformations}. A 
feed cluster transformation  
${\bf c}: {\bf i} \to {\bf i'}$ provides a unitary operator 
$$
{\bf K}_{\bf c^o}: L^2({\cal A}_{{\bf i'}}^+) \lra L^2({\cal A}_{{\bf i}}^+).  
$$
Indeed,  we assigned to   a feed 
 mutation ${\bf i} \to {\bf i'}$ 
an intertwiner 
${\bf K}_{{\bf i'} \to {\bf i}}: L^2({\cal A}^+_{\bf i'}) \to 
L^2({\cal A}^+_{\bf i})$. Further, an automorphism $\sigma$ 
of a feed ${\bf i}$  gives rise 
to 
a unitary operator given by a permutation of  
coordinates in the space ${\cal A}_{\bf i}^+$. The  
feed cluster transformation  
 ${\bf c}$  is a composition of mutations and automorphisms. 
Taking the reverse 
composition of the corresponding intertwiners, we get the map 
${\bf K}_{\bf c^o}$.

A feed cluster transformation ${\bf c}$ gives rise to 
cluster transformations  ${\bf c}^x_q$ and 
${\bf c}^d_q$ of the quantum spaces ${\cal X}_q$  and ${\cal D}_q$. 
Denote by $\gamma_{\bf c^o}$ the $\ast$-algebra automorphism of 
${\bf L}$ corresponding to  ${\bf c}$. 

The following theorem is one of the main results of this paper. 

\begin{theorem} \label{K0}
(i) The operator ${\bf K}_{\bf c^o}$ provides  a map of 
 Schwartz spaces 
$$
{\bf K}_{\bf c^o}: {\cal S}_{\bf i'} \lra {\cal S}_{\bf i}.
$$ 
It
intertwines the automorphism $\gamma_{\bf c^o}$ of ${\bf L}$, i.e.
for any $A \in {\bf L}$ and $s \in {\cal S}_{\bf i}$ one has
\begin{equation} \label{tu31}
{\bf K}_{\bf c^o}\widehat A {\bf K}_{\bf c^o}^{-1} s = \widehat {\gamma_{\bf c^o}(A)} s.
\end{equation}

(ii) Suppose 
that the cluster transformations ${\bf c}^x_q$ and 
${\bf c}^d_q$ of the quantum spaces ${\cal X}_q$  and ${\cal D}_q$ 
are identity maps. Then the operator ${\bf K}_{\bf c^o}$ 
is proportional to the identity: ${\bf K}_{\bf c^o} = \lambda_{\bf c^o}{\rm Id}$, 
where $|\lambda_{\bf c^o}|=1$. 
\end{theorem} 

It is sufficient to prove the part (i) of the  theorem for mutations. 
The relation between the mutations and the intertwiners is 
summarized in the following diagram: 
$$
\begin{array}{ccccc}
{\bf L} & \stackrel{\mu_k'}{\lra}& {\bf L}
&\stackrel{\mu_k^{\sharp}}{\lra}&{\bf L} \\
&&&&\\
{\cal S}_{\bf i'} &\stackrel{{\bf K}'}{\lra} & {\cal S}_{\bf i}  &\stackrel{{\bf K}^{\sharp}}{\lra} 
& {\cal S}_{\bf i}  
\end{array}
$$ 
Here the algebra ${\bf L}$ 
acts on the Schwartz spaces ${\cal S}_{*}$ at the bottom. 
The operators ${\bf K}'$ and ${\bf K}^{\sharp}$ 
intertwine the isomorphisms at the top line. 

\vskip 3mm
Here are 
main issues and a strategy of the proof. 

The part (i) implies the part (ii) if we can show that the algebra ${\bf L}$ 
is big enough, e.g. of maximal functional dimension. 
So we start from the part (ii), and explain why it requires the part (i).

First, we do not know explicitly neither all 
trivial feed cluster transformations, nor feed cluster transformations 
satisfying the condition of the the part (ii) of the theorem.\footnote{We conjecture that these 
two types of cluster transformations coincide.} 
We suspect that in many 
cases 
relations (\ref{K10}) generate all of them,
 but this is not always the case, 
and is not known at the moment  even
for the cluster varieties ${\cal X}_{G,S}$ with $G = SL_m$, $m> 3$. 
This eliminates a hope for a proof  by an explicit calculation  
of the integrals providing the 
intertwiner. 

To prove the part (ii), 
we would like to show that the operator ${\bf K}_{\bf c^o}$ commutes with the operators from algebra 
${\bf D}^q_{\bf i}\otimes {\bf D}^{q^{\vee}}_{\bf i^{\vee}}$ of the $q$-difference 
operators acting on the space $W_{\bf i}$. 
This, for $\hbar\not \in \Q$, easily implies that ${\bf K}_{\bf c^o} = \lambda_{\bf c^o} {\rm Id}$. 
The algebra ${\bf D}^q_{\bf i}$ alone 
would not do the job, since there are many operators commuting with it, e.g. 
${\bf D}^{q^{\vee}}_{\bf i^{\vee}}$. Since ${\bf K}_{\bf c^o}$ 
depends continuously on $\hbar$, the claim follows then for all $\hbar$.

 Since the operator ${\bf K}_{\bf c^o}$ is a composition of 
  elementary intertwiners ${\bf K}$, we would like to  
have commutation relations with them. The problem is that 
we deal with unbounded operators,  
so $A{\bf K} = {\bf K}A$ means that 
$A{\bf K}f = {\bf K}Af$ for $f$'s from a certain dense subspace 
${\cal D}_A \subset L^2({\cal A}^+)$. However {\it a priori} 
it is not clear what the intertwiner  ${\bf K}$ 
does with the domain ${\cal D}_A$, 
so if both $A$ and $B$ commute with ${\bf K}$, {\it a priori} 
it is not clear why 
 $AB$ should commute with ${\bf K}$,  
or even whether $B{\bf K}Af$ is defined. To make sense out of this we restrict 
to the subalgebra ${\bf L}$ 
-- otherwise meaningless denominators would appear in the commutation formulas,  
and work with the corresponding Schwartz spaces -- 
the intertwiners do not respect the spaces $W_{\bf i}$. 
However then we have to show that the subalgebra ${\bf L}$ is big enough. 
This creates some issues of purely algebraic nature. 

\vskip 3mm
Here is  our scheme of the proof: 

1) First, we prove in Theorem \ref{12.13.06.10} 
that the operator ${\bf K}^\sharp$ intertwines the mutation 
automorphism 
of the 
algebra ${\bf L}$  on the subspace $W_{\bf i}$. Here we use 
the pair of difference relations for the function 
$\Phi^\hbar$, and in addition to this the following remarkable property  
of the function 
 $\Phi^\hbar(z)$: it is  analytic in the upper half plane, and 
its inverse is analytic in the lower half plane (!) -- see Property 
${\bf B}7$ in Section 4 -- plus at most exponential growth of $\Phi^\hbar(z)$
when $|{\rm Re } z|\to \infty$. 

We would like to stress that the pair of difference relations for the function 
$\Phi^\hbar$ is not sufficient for the proof. 
Indeed, it is clear from the proof,  that 
any pole of the function $\Phi^\hbar(z)$ in the upper 
half plane, or any zero in the lower half plane would be an obstraction 
to the intertwining property for general $\hbar$. 

2) Second, we prove that the space $W_{\bf i}$ is dense in the Freschet space 
$S_{\bf i}$. This eventually allows us to prove that 
the intertwiner respects the Schwartz spaces. 

3) Finally, we  show that the subalgebra ${\bf L}$ is big enough. 
This allows to use 2) to deduce (i) from (ii) in Theorem \ref{K0}. 

%Our proof follows \cite{Go2}, 
%where the simplest instance of this programm, not overshadowed by technical issues,  
%was implemented in detail.  It may be useful 
%to look at {\it loc. cit.} first. 

\vskip 3mm
Recall the dual %{ saturated cluster modular groupoid} 
groupoid $\widehat {\cal G}^o$ (Section 2.1).  
Its objects are feeds equivalent to a given feed ${\bf i}$. 
The set of morphisms ${\rm Hom}({\bf i}, {\bf i'})$ 
consists of feed cluster transformations ${\bf i'} \to {\bf i}$ 
modulo the ones (\ref{K10}). 
 The group $\widehat \Gamma$  acts 
by automorphisms of  the quantum double (Theorem \ref{8.15.05.10gg}a).

\vskip 3mm
{\bf The $(h+2)$-gon relations}. The 
automorphism  $\sigma_{ij}:I \to I$ permuting $i$ and $j$ 
induces an automorphism of $L^2({\cal A}_{\bf i}^+)$,
 also denoted by $\sigma_{ij}$. Denote by ${\bf K}_{\mu_i^o}$ 
the intertwiner for the mutation of a feed ${\bf i}$ at the direction $i$. Set 
$$
{\cal K}_{i,j}:= 
\sigma_{ij}\circ {\bf K}_{\mu_i^o}.
$$

\begin{theorem} \label{K1}
(i) Given a pair $\{i,j\}\subset I$ such that $\varepsilon_{ij}  = - p\varepsilon_{ji} = p 
\in \{0, 1, 2, 3\}$, one has 
\begin{equation} \label{K5}
{\cal K}_{i,j}^{h+2} = \lambda ~{\rm Id}, \qquad |\lambda|=1,
\end{equation}
 where $h=2,3,4,6$ for $p=0,1,2,3$ respectively. 

(ii) Assigning to a feed ${\bf i}$ the Hilbert space 
$L^2({\cal A}_{{\bf i}}^+)$, and to a feed cluster transformation ${\bf c}: 
{\bf i} \to {\bf i'}$ the operator ${\bf K}_{\bf c^o}$, 
we get a unitary projective representation 
of the saturated cluster modular groupoid $\widehat {\cal G}^o $, and hence 
 a unitary projective representation of the 
saturated cluster modular group $\widehat \Gamma_{\bf i}$ in the Hilbert space  
$L^2({\cal A}_{{\bf i}}^+)$.

(iii) If ${\rm det}\varepsilon_{ij} \not = 0$,  
we get 
a unitary projective representation 
of the cluster modular groupoid ${\cal G}^o $, and hence 
 a unitary projective representation of the 
cluster modular group $\Gamma_{\bf i}$ in  
$L^2({\cal A}_{{\bf i}}^+)$.
\end{theorem}

{\bf Proof}. The part 
(ii) follows immediately from (i). 
Since relations (\ref{K10}) are valid for the 
quantum cluster spaces ${\cal X}_q$  and ${\cal D}_q$ (Theorem \ref{8.15.05.10gg}a), 
the part (i) follows from Theorem \ref{K0}. 

If  ${\rm det}\varepsilon_{ij} \not = 0$, by 
Lemma \ref{8.24.05.1},  ${\bf c}^a = {\rm Id}$ implies 
${\bf c}^x_q = {\rm Id}$. Thanks to the quantum Laurent Phenomenon 
Theorem \cite{BZq}, 
 ${\bf c}^a = {\rm Id}$ implies ${\bf c}^a_q = {\rm Id}$. 
There is a 
canonical projection 
${\cal A}_q \times {\cal A}_q^o \lra {\cal D}_q$,  
defined for every feed by the canonical $q$-deformation of the 
monomial map $\varphi$. It is surjective thanks to ${\rm det}\varepsilon_{ij} \not = 0$. 
Using this, we deduce that 
  ${\bf c}^d_q = {\rm Id}$. So 
the part (iii) follows from Theorem \ref{K0}.
The theorem is proved.

%\vskip 3mm
%{\bf Quantum ${\cal D}$-varieties}.  
% Given a feed ${\bf i}$, we realize the algebra ${\bf L}$ as a 
% subalgebra ${\bf L}_{\bf i} \subset 
%{\bf D}^q_{\bf i}\otimes {\bf D}^{q^{\vee}}_{\bf i}$. 
%Thanks to Theorem \ref{8.15.05.10gg} there is a functor 
%$$
%{\cal L}: \mbox{the modular groupoid $\widehat 
%{\cal G}^o$} \lra \mbox{the category of non-commutative $\ast$-algebras}
%$$
%assigning to a feed ${\bf i}$ the $\ast$-algebra ${\bf L}_{\bf i}$, and to 
%a cluster transformation 
%${\bf c}: {\bf i} \to {\bf i'}$ an isomorphism 
%$\gamma_{\bf c^0}: {\bf L}_{\bf i'} \stackrel{\sim}{\to} 
%{\bf L}_{\bf i}$. 
%Thanks to Theorem \ref{K1} there is another functor
%$$
%{\cal Q}: \mbox{the modular groupoid $\widehat {\cal G}^o$} \lra 
%\mbox{the category of triples of topological vector spaces}
%$$
%assigning to a feed ${\bf i}$ a triple of topological spaces $S_{\bf i} \subset 
%L^2({\cal A}_{\bf i})\subset {\cal S}_{\bf i}^*$, and to a cluster transformation 
%${\bf c}: {\bf i} \to {\bf i'}$ an operator 
%${\bf K}_{\bf c^o}$ identifying the corresponding triples, and intertwining the action of the algebras 
%${\bf L}_{\bf i'}$ and ${\bf L}_{\bf i}$ on the Schwarz spaces ${\cal S}_{\bf i}$ and ${\cal S}_{\bf i'}$, 
%and hence on its topological duals ${\cal S}^*_{\bf i}$ and ${\cal S}^*_{\bf i'}$. 

%\begin{definition} \label{7.27.07.1}
%The quantum ${\cal D}$-variety is the pair of functors $({\cal L}, {\cal Q})$. 
%\end{definition}
%\vskip 3mm

\subsection{Commutation relations for the intertwiner}

Recall that a mutation ${\mu}_k: {\bf i} \to {\bf i}'$ provides  
automorphisms $\mu^\sharp_{k}$ and 
 $\mu'_{k}$ of ${\bf L}$. They are tensor products of the corresponding 
automorphisms of ${\Bbb L}_q$ and ${\Bbb L}_{q^{\vee}}$.

\begin{theorem} \label{12.13.06.10} 
For any $w \in W_{\bf i'}$,  $A \in {\bf L}$ one has
\begin{equation} \label{7.26.07.2}
{\bf K}^{\sharp}\widehat A w = \widehat {\mu_{k}^\sharp(A)}{\bf K}^{\sharp}w, \qquad 
{\bf K}'\widehat A w = \widehat {\mu_{k}'(A)}{\bf K}^{\sharp}w, 
\end{equation}
 Thus both $\widehat {\mu_{k}^\sharp(A)}{\bf K}^{\sharp}w$ and 
$\widehat {\mu_{k}'(A)}{\bf K}'w$ lie in $L^2({\cal A}_{\bf i}^+)$. 
\end{theorem}

{\bf Proof}. The proof of the right identity in (\ref{7.26.07.2}) is 
straightforward. 
Let us prove  the left identity in (\ref{7.26.07.2}). 
We do it for $A \in {\Bbb L}_q$. 
The proof for $A \in {\Bbb L}_{q^{\vee}}$ is completely similar, 
using the second 
difference relations for the function $\Phi^\hbar$. So from now on $A \in {\Bbb L}_q$.

Let ${\Bbb L}_q'$ be the space of Laurent $q$-polynomials $F$ in $B_i, X_i$ 
such that $\mu_{k}^\sharp(F)$ is again  a Laurent $q$-polynomial. 
It is easy to check that the following elements belong to ${\Bbb L}'_q$:
\begin{equation} \label{qp1001a}
[1] \quad B_i^{\pm 1}, i \not = k, \qquad [2] \quad B_k(1+q_k\widetilde X_k), \qquad 
[3] \quad 
(1+q_k X_k)B_k^{-1}, 
\end{equation}
\begin{equation} \label{qp1002}
[4] \quad X_k^{\pm 1}, ~~\qquad [5] \quad 
X_i\frac{\Psi^{q_k}(X_k)}{\Psi^{q_k}(q_k^{-2\varepsilon_{ik}}X_k)}
\mbox{~~if $\varepsilon_{ik} \geq 0$}, \qquad [6] \quad X_i^{-1} \mbox{~~if $\varepsilon_{ik} \geq 0$}, 
\end{equation}
\begin{equation} \label{qp1003}
\qquad [7] \quad
\frac{\Psi^{q_k}(q_k^{-2\varepsilon_{ik}}X_k)}
{\Psi^{q_k}(X_k)}X_i^{-1}\mbox{~~if $\varepsilon_{ik} \leq 0$}, \qquad 
[8] \quad X_i\mbox{~~if $\varepsilon_{ik} \leq 0$}. 
\end{equation}
We need the following simple algebraic lemma.

\begin{lemma} \label{12.13.06.1qw} Any element of ${\Bbb L}'_q$ is a sum of products of the 
expressions (\ref{qp1001a})-(\ref{qp1002}).
\end{lemma} 

{\it Strategy}. We prove the commutation relations from the 
Theorem for each of the expressions [1]-[8]. 
This easily implies the claim for their products and sums, i.e. for the algebra ${\Bbb L}'_q$. 
Thus we get Theorem  \ref{12.13.06.10}, since by 
the very definition ${\Bbb L}_q\subset {\Bbb L}'_q$.

Since $\widehat B_i$ and $\widehat x_k, \widehat {\widetilde x_k}$ commute, 
the commutation relations are obvious for [1]. 
 Since $\widehat x_k$ and 
$\widehat {\widetilde x_k}$ commute, 
they are obvious for  [4]. The commutation relations for [6] and [7] are deduced formally from the ones 
[5] and [8] using the canonical isomorphism of quantum cluster 
${\cal X}$-varieties related to the pairs $({\bf i}, q)$ and  $({\bf i}^o, q^{-1})$, where  
$\varepsilon^o_{ij} = - \varepsilon_{ij}$.   Finally, [3] is reduced to [2] via a similar trick 
using the involution $i: {\cal D}_q \to {\cal D}^0_q$ from Theorem \ref{8.15.05.10gg}b).  
So one needs to prove relations [2], [5],[8]. 

\vskip 3mm

[2]. Recall one of the two difference relations 
 for the function $\Phi^{\hbar_k}(z)$, see {\bf B5} in Section 4: 
\begin{equation} \label{12.13.06.1}
\Phi^{\hbar_k}(z+2\pi ih_k) = \Phi^{\hbar_k}(z)(1+q_ke^z) \quad <=> \quad 
\Phi^{\hbar_k}(z-2\pi ih_k) = \Phi^{\hbar_k}(z)(1+q_k^{-1}e^z)^{-1}.
\end{equation}
Recall  the notation $\alpha_k^\pm$, see (\ref{handy}), and 
$
B^{\sharp}_k = B_k(1+ q_kX_k)(1+ q_k\widetilde X_k)^{-1}$.
Set
$$
T_{2\pi i \hbar}\varphi(c)= T^{(c)}_{2\pi i \hbar}\varphi(c):= \varphi(c + 2\pi i \hbar). 
$$  

Here is an important point. 
For the shift in the imaginary direction, the usual identity
\begin{equation} \label{SUBTLE}
\widehat A_k {\cal F}^{-1}_{a_k}\varphi(c)
 = {\cal F}^{-1}_{a_k} T_{2\pi i\hbar } \varphi(c)
\end{equation}
requires an assumption that the function 
$\varphi(c)$ admits an analytic continuation 
to the strip $0 \leq {\rm Im}(c) \leq 2\pi \hbar$, and 
 decays sufficiently fast in this strip when $|{\rm Re}(c)|\to \infty$. 
Indeed, 
$$
{\cal F}^{-1}_{a_k} T_{2\pi i\hbar } \varphi(c) = 
\int_{C_0}e^{-a_kc/2\pi i \hbar}\varphi(c+2\pi \hbar)dc =
e^{a_k}\int_{C_0}e^{-a_k(c+2\pi i \hbar)/2\pi i \hbar}\varphi(c+2\pi \hbar)dc 
=
$$ 
$$
e^{a_k}\int_{C_{2\pi \hbar}}e^{-a_k(c)/2\pi i \hbar}\varphi(c)dc  \stackrel{?}{=} 
e^{a_k}\int_{C_{0}}e^{-a_k(c)/2\pi i \hbar}\varphi(c)dc = 
\widehat A_k {\cal F}^{-1}_{a_k}\varphi(c).
$$
To justify the marked by $?$ equality we should move the integration countour 
$C_{2\pi \hbar}:= \{c|{\rm Im}c=2\pi i \hbar\}$ to the one 
$C_{0}:=\{c|{\rm Im}c=0\}$. Below we use this in the situation when 
$\varphi(c)$ is analytic in the strip, decays as $e^{-c^2}$ at infinity, 
so one indeed can  move the countour.

Using (\ref{1.05.06.14})-(\ref{1.05.06.14a}), we have 
$$
{\bf K}^{\sharp} \left(\widehat A_k (1+q_k\widehat {\widetilde X_k}) w\right) = 
$$
\begin{equation} \label{7.24.07.1}
{\cal F}_{a_k}^{-1} \circ
\Phi^{\hbar_k}\Bigl(-d^{-1}_k c -  \alpha_k^+ \Bigr)
\Phi^{\hbar_k}\Bigl(-d^{-1}_k c - \alpha_k^-\Bigr)^{-1}\circ 
T_{2\pi i\hbar }\circ {\cal F}_{a_k}\circ (1+q_k
\widehat {\widetilde X_k})w = 
\end{equation}
$$%\begin{equation} \label{7.24.07.1q}
{\cal F}_{a_k}^{-1} \circ 
T_{2\pi i\hbar}\circ T_{-2\pi i\hbar}\circ
\Phi^{\hbar_k}\Bigl(-d^{-1}_k c -  \alpha_k^+ \Bigr)
\Phi^{\hbar_k}\Bigl(-d^{-1}_k c - \alpha_k^-\Bigr)^{-1}
 \left(1+q^{-1}_k {\rm e}^{-d^{-1}_k c - \alpha_k^-}\right)\circ 
T_{2\pi i\hbar }\circ {\cal F}_{a_k}w = 
$$%\end{equation}
\begin{equation} \label{7.24.07.1q}
\widehat A_k \circ {\cal F}_{a_k}^{-1} \circ 
T_{-2\pi i\hbar}\circ
\Phi^{\hbar_k}\Bigl(-d^{-1}_k c -  \alpha_k^+ \Bigr)
\Phi^{\hbar_k}\Bigl(-d^{-1}_k (c+2\pi i\hbar) - \alpha_k^-\Bigr)^{-1}\circ 
T_{2\pi i\hbar }\circ {\cal F}_{a_k}w.
\end{equation}
Here we used the following  facts, guaranteed by $\hbar >0$ 
and property ${\bf B}7$ from Section 4:

(i) The function $\Phi^{\hbar_k}(c)$ is analytic 
in the upper half plane.

(ii) The function $\Phi^{\hbar_k}(c -  2\pi i\hbar_k)^{-1}$ is analytic 
in the strip $0 \leq {\rm Im}(c) \leq 2\pi \hbar$ -- indeed, $\Phi^{\hbar_k}(c)^{-1}$ is analytic in the lower 
half plane. 

(iii) The functions $\Phi^{\hbar_k}(c)^{\pm 1}$  grow 
at most exponentially when 
$|{\rm Re}(c)| \to \infty$ , 
while the function $w(a)$ and its Fourier transform 
decay as $e^{-a^2}$. 

Thus we can use (\ref{SUBTLE}) to go from (\ref{7.24.07.1}) 
to (\ref{7.24.07.1q}). 
%The operator ${\rm exp}
%\Bigl(-2\pi i\hbar \frac{\partial}{\partial c}\Bigr)$ acts as a 
%shift $f(c) \lms f(c-2\pi i\hbar)$. Thus 
%$$
%{e}^{-2\pi i\hbar \partial/\partial c}F(c) {e}^{2\pi i\hbar \partial/\partial c}= F(c-2\pi i \hbar).
%$$ 
Identities (\ref{12.13.06.1}) imply 
$$
{\cal F}_{a_k}^{-1} \circ \Phi^{\hbar_k}\Bigl(-d^{-1}_k (c-2\pi i\hbar) -\alpha_k^+ \Bigr) \circ {\cal F}_{a_k}= (1+ q_k\widehat X_k)
\Phi^{\hbar_k}(\widehat x_k).
$$
%$$
%{\cal F}_{a_k}^{-1} \circ \Phi^{\hbar_k}\Bigl(-d^{-1}_k (c-2\pi i\hbar) -\alpha_k^-\Bigr)^{-1} \circ {\cal F}_{a_k} \circ (1+ q_k\widehat {\widetilde X_k}) = 
%\Phi^{\hbar_k}(\widehat {\widetilde x}_k)^{-1}.
%$$
Therefore 
\begin{equation} \label{7.24.07.2}
{\bf K}^{\sharp} \left(\widehat 
A_k (1+q_k\widehat {\widetilde X_k}) w\right) = e^{a_k} (1+ q_k{\widehat X_k})
\Phi^{\hbar_k}(\widehat x_k) 
\Phi^{\hbar_k}(\widehat {\widetilde x}_k)^{-1} w = 
\widehat A_k(1+ q_k\widehat X_k){\bf K}^\sharp w.
 \end{equation} 
 We proved the claim for the expression  [2].

\vskip 3mm
 [8]. Since 
$F(\widehat {{\widetilde x}_k})$ commutes with $X_j$ for any function $F$, 
the claim boils down to the following. 

Setting 
$\widetilde \omega:={\cal F}_{a_k}\omega$, and 
using difference equations (\ref{12.13.06.1}), 
we have to prove  an equality  
$$
{\cal F}_{a_k}^{-1}
\frac{\Phi^{\hbar_k}\Bigl(-d^{-1}_k c -\alpha_k^+ \Bigr)}
{\Phi^{\hbar_k}\Bigl(-d^{-1}_k c -\alpha_k^+ +\varepsilon_{ik}2\pi i\hbar_k\Bigr)}
{\cal F}_{a_k}
\widehat X_i
{\cal F}_{a_k}^{-1}
\Phi^{\hbar_k}\Bigl(-d^{-1}_k c -\alpha_k^+ \Bigr)
\widetilde \omega \stackrel{?}{=} 
$$
$$
{\cal F}_{a_k}^{-1}\Phi^{\hbar_k}\Bigl(-d^{-1}_k c -\alpha_k^+ \Bigr)
{\cal F}_{a_k}
\widehat X_i 
{\cal F}_{a_k}^{-1}
\widetilde \omega. 
$$
Equivalently, cancelling two factors on the left, one has to show that 
$$
\Phi^{\hbar_k}\Bigl(-d^{-1}_k c -\alpha_k^+ +\varepsilon_{ik}2\pi i\hbar_k\Bigr)^{-1}
{\cal F}_{a_k}
\widehat X_i
{\cal F}_{a_k}^{-1}
\Phi^{\hbar_k}\Bigl(-d^{-1}_k c -\alpha_k^+ \Bigr)
\widetilde \omega \stackrel{?}{=}
{\cal F}_{a_k}
\widehat X_i 
{\cal F}_{a_k}^{-1}
\widetilde \omega. 
$$
Denote by $(a;a_i, a_k)$ the arguments of $\omega$, and by $(a;a_i, c)$ the ones of 
$\widetilde \omega$. Here is the idea. The only variables 
participating non-trivially 
in the 
transformations are the ones $a_i$, and $a_k$ or $c$.  
So if we pay attention to them only, and do a formal comutation, 
the operator ${\cal F}_{a_k} \widehat X_i{\cal F}_{a_k}^{-1}$ acts on 
them as  shift 
$$
c \lms c-[\varepsilon_{ik}]_+2\pi i \hbar, \quad  a_i \lms a_i+ 2\pi i \hbar. 
$$
To justify this formal computation, we use (\ref{SUBTLE}) 
and  the fact that the function $\Phi^{\hbar_k}$ 
is analytic in the upper half plane, just like in the proof of [2]. 
Finally, 
since $\varepsilon_{ik}\hbar_k = - \varepsilon_{ki}\hbar_i$ and 
$[\widehat \varepsilon_{ik}]_+ - [-\widehat \varepsilon_{ik}]_+ 
= \widehat \varepsilon_{ik}$
$$
\Phi^{\hbar_k}\Bigl(-d^{-1}_k (c-[\varepsilon_{ik}]_+2\pi i\hbar) -\alpha_k^+  - [\varepsilon_{ki}]_+2\pi i\hbar_i \Bigr) = 
\Phi^{\hbar_k}\Bigl(-d^{-1}_k c -\alpha_k^+  +\varepsilon_{ik}2\pi i\hbar_k\Bigr).
$$
The statement is proved. 
\vskip 3mm

[5]. We have to compare two expressions: 
$$
{\cal F}_{a_k}^{-1}{\cal F}_{a_k}\widehat X_i
{\cal F}_{a_k}^{-1}\Phi^{\hbar_k}\widetilde \omega \quad \mbox{and} \quad 
{\cal F}_{a_k}^{-1}\Phi^{\hbar_k}{\cal F}_{a_k} \widehat X_i
\frac{\Psi^{q_k}(\widehat X_k)}{\Psi^{q_k}(q_k^{-2\varepsilon_{ik}}\widehat X_k)} {\cal F}_{a_k}^{-1}\widetilde \omega. 
$$
The computation is very similar to the one for [8], and thus omitted. 
The Proposition and hence Theorem \ref{12.13.06.10}  are proved.

\vskip 3mm

\subsection{The Schwartz spaces ${{\cal S}}_{\bf i}$ and the intertwiners}
The key property of the Schwartz space ${\cal S}_{{\bf i}}$ which we use below is the following.
\begin{proposition} \label{qp12}
The space $W_{\bf i}$ is dense in the Schwartz space ${\cal S}_{{\bf i}}$.
\end{proposition}

One interprets 
Proposition \ref{qp12} by saying that
the {\it $\ast$-algebra ${\bf L}$ is  essentially self-adjoint in $L^2({\cal A}_{\bf i}^+)$}.

\vskip 3mm

{\bf Proof}. Here is the scheme of the proof. 

\begin{itemize} 

\item \begin{lemma} \label{qp141}
For any $w \in W_{{\bf i}}, s \in {\cal S}_{{\bf i}}$, the convolution
$s \ast w $ lies in ${\cal S}_{{\bf i}}$.
\end{lemma}

\item 
Let $w_{\varepsilon}:= (2\pi)^{-\frac{n}{2}}\varepsilon^{-n} e^{-\frac{1}{2}|x/\varepsilon|^2}\in W_{\bf i}$ be a sequence
converging as $\varepsilon\to 0$ to  the $\delta$-function at $0$. 
By Lemma \ref{qp141},  $w_{\varepsilon} \ast s$ lies in ${\cal S}_{{\bf i}}$. 
Clearly one has in the topology of
${\cal S}_{{\bf i}}$
\begin{equation} \label{qp13}
\lim_{\varepsilon\to 0}w_{\varepsilon} \ast s  = s(x).
\end{equation}

\item \begin{lemma} \label{qp14}
For any $w \in W, s \in {\cal S}_{{\bf i}}$,
the Riemann sums for the integral 
\begin{equation} \label{qp131}
s \ast w (x) =
\int_{V}s(v) \omega(x-v)dv = \int_{V}s(v) T_vw(x)dv
\end{equation}
 converges in the topology of
${\cal S}_{{\bf i}}$ to the convolution $s \ast w$.
\end{lemma}

\item 
For every $\varepsilon >0$, 
we approximate the function $(s\ast w_{\varepsilon}) (x) \in {\cal S}_{{\bf i}}$ 
by the finite Riemann sums of the integral (\ref{qp131}). 
Each of these Riemann sums is an element of $W$. 
Then, letting $\varepsilon \to 0$ and using Lemma \ref{qp141} 
togerther with the triangle inequality,
 we get the Proposition.

\end{itemize}

\vskip 3mm
{\bf Proof of Lemma \ref{qp141}}. Let $V$ be a finite dimensional vector space with a norm $|v|$  
and the corresponding Lebesgue measure 
$dv$, and $v, x \in V$, set
 $T_v f(x) := f(x-v)$. Write
$$
s \ast w (x)= \int_{V}w(v) (T_vs)(x)dv.
$$
Below $V \stackrel{\sim}{=} \R^n \stackrel{\sim}{=} {\cal A}_{{\bf i}}^+.$ 
For any seminorm $\rho_B$ on ${\cal S}_{{\bf i}}$ the operator
$T_v: ({\cal S}_{{\bf i}}, \rho_B) \lra ({\cal S}_{{\bf i}}, \rho_B)$ is a bounded operator
with the norm bounded by $e^{c|v|}$ for some $c$ depending on $B$. Thus the operator
$\int_{V}w(v) T_vdv$ is a bounded operator on $({\cal S}_{{\bf i}}, \rho_B)$.
This implies the Lemma.

\vskip 3mm
{\bf Proof of Lemma \ref{qp14}}. 
Let us show first that (\ref{qp131})  is convergent in
$L^2({\cal A}_{\bf i}^+)$. The key fact is that a shift of
$w\in W_{\bf i}$  quickly becomes essentially orthogonal to $w$. More precisely,
in the important for us case when
 $w(x) = {\rm exp}(-\alpha(x\cdot x)^2/2+b\cdot x)$, $\alpha>0$,  
(this includes any $w \in W_0$)  we have
\begin{equation} \label{qp132}
( w(x), T_v w(x)) < C_we^{-\alpha(v \cdot v)^2/2+ (b-\overline b) \cdot v}.
\end{equation}
Therefore in this case
$$
\Bigl(\int_{V}s(v) T_vw(x)dv, \int_{V}s(v) T_vw(x)dv\Bigr)
$$
$$
\stackrel{(\ref{qp132})}{\leq} C_w\int_{V}\int_{V}
e^{-a(v_1-v_2)^2/2 + (b-\overline b)(v_1-v_2)} |s(v_1) s(v_2)| dv_1 dv_2 =
$$
$$
C_w\int_{V}e^{-\alpha (v \cdot v)^2/2+(b-\overline b)\cdot v}
\int_{V}|s(t) s(t+v)|  dtdv \leq
 C_w||s||_{L^2}^2\int_{-\infty}^{\infty}e^{-\alpha(v\cdot v)^2/2+(b-\overline b)\cdot v}  dv.
$$
We leave the case of an arbitrary $w$ to the reader: it is not used in the proof of the theorem. 
The convergence with respect to the seminorm $||Bf||$ is proved by 
the same argument. The Lemma is proved. The proof of the Proposition is fished.

\vskip 3mm
{\bf Remark}.  The same arguments show that the space $W_{\bf i}$ is dense
in the space ${\cal S}_{{\bf L'}, {\bf i}}$ defined for any subalgebra ${\bf L'}$ of ${\bf L}$.

\vskip 3mm
{\bf Proof  Theorem \ref{K0}(i)}. The case of ${\bf K}'$
is obvious. 
To prove it for  ${\bf K}^{\sharp}$
we use Theorem \ref{12.13.06.10}. 

We use shorthand ${\bf K}:= {\bf K}^{\sharp}$ and ${\gamma}:= {\gamma}^{\sharp}$. 
Let us show that ${\bf K}^{-1}s \in {\cal S}_{\bf i'}$ for an 
$s \in  {\cal S}_{\bf i}$. For this  we need to check that
for any $B \in {\bf L}$ the functional $w' \to ({\bf K}^{-1}s, \widehat B^*w')$ is 
continuous on $W_{\bf i'}$.

Since $W_{\bf i}$ is dense in ${\cal S}_{\bf i}$ by Proposition \ref{qp12}, 
there is a sequence $v_i \in W_{\bf i}$ converging to
$s$ in  ${\cal S}_{\bf i}$. This means that
\begin{equation} \label{tu32}
\lim_{i\to \infty}(\widehat Bv_i, w) 
= (\widehat Bs,w)\quad \mbox{for any $B \in {\bf L}$, $w \in W_{\bf i}$}.
\end{equation}
In particular, let us put for the record the fact that
\begin{equation} \label{tu32ss}
\mbox{The sequences $v_i$ and 
$\widehat {\gamma(B^*)}v_i$ converge in $L^2({\cal A}^+_{\bf i})$ to, respectively, 
 $s$ and $\widehat {\gamma(B^*)}s$}.
\end{equation}
One has 
$$
({\bf K}^{-1}s, \widehat B^*w') =
(s, {\bf K}\widehat B^*w')
\stackrel{\mbox{Th. \ref{12.13.06.10}}}{=}
(s, \widehat {\gamma(B^*)}{\bf K}w') \stackrel{(\ref{tu32ss})}{=}
\lim_{i\to \infty} (v_i, \widehat {\gamma(B^*)}{\bf K}w') 
$$
$$
\stackrel{\rm def}{=} 
\lim_{i\to \infty} (\widehat {\gamma(B)}v_i, {\bf K}w') \stackrel{(\ref{tu32ss})}{=}
(\widehat {\gamma(B)}s, {\bf K}w') =  ({\bf K}^{-1}\widehat {\gamma(B)}s, w').
$$
Since the functional on the right is continuous,  
${\bf K}^{-1}s \in {\cal S}_{\bf i'}$, 
and we have (\ref{tu31}). 
The theorem is proved.
\vskip 3mm

\begin{definition}
The  space of distributions  ${\cal S}^*_{\bf i}$ is the topological dual to the Schwartz space 
${\cal S}_{\bf i}$. 
\end{definition}

\begin{corollary}
The operator ${\bf K}_{\bf c^o}$ gives rise to an isomorphism  
of 
 topological spaces ${\bf K}^*_{\bf c}: {\cal S}^*_{\bf i} \lra {\cal S}^*_{\bf i'}$ 
intertwining the automorphism $\gamma_{\bf c^o}$ of ${\bf L}$. 
\end{corollary}

\vskip 3mm
\subsection{Relations for the intertwiners}

{\bf Proof of part (ii) of Theorem \ref{K0}}. 
{\it Step 1}. 

\begin{proposition} \label{P1}
Assume that ${\rm det}(\varepsilon_{ij})=1$. 
The intertwiner ${\bf K}_{\bf c^o}$ 
corresponding to a trivial cluster transformation 
${\bf c}$ is proportional to the identity. 
\end{proposition}

{\bf Proof}. 
Let us show first that the algebra ${\bf L}$ is big enough. 
Recall the quantum cluster algebras \cite{BZq}. 
Thanks to the quantum Laurent Phenomenon from {\it loc. cit.}, 
for every feed ${\bf i}$ each monomial  
$\prod_{i\in I}A_i^{n_i}$ where $n_i \geq 0$ is 
a  universally Laurent polynomial. 
Since ${\rm det}(\varepsilon_{ij})=1$, 
there are canonical isomorphisms of feed tori 
\begin{equation} \label{1.5.07.190}
{\cal A}_{\bf i} \times {\cal A}_{\bf i}^o \stackrel{\sim}{\lra} {\cal D}_{\bf i}
\stackrel{\sim}{\lra}  
{\cal X}_{\bf i} \times {\cal X}_{\bf i}^o 
\end{equation}
as well as the corresponding quantum feed tori algebras. 
Thanks to this  
there is a cone (i.e. a semigroup) 
$C_+$ of full rank
${\rm dim}{\cal D}$ in the group of characters $X^*({\cal D}_{{\bf i}})$
of the torus ${\cal D}_{\bf i}$  such that the quantum 
monomials  assigned to its vectors are belong to ${\bf L}$. 
This cone is 
the preimage of the square of the 
positive octant cone under the map 
$X^*({\cal D}_{\bf i}) \to X^*({\cal A}_{\bf i}) \times X^*({\cal A}^o_{\bf i})$ 
provided by  isomorphism (\ref{1.5.07.190}). 
This way we know that ${\bf L}$ is  big. 

\begin{lemma} \label{p1p} 
Let $A \in {\bf L}$ be a monomial. Then 
$\widehat A^{-1}{\bf K}_{\bf c^o}w = {\bf K}_{\bf c^o}\widehat A^{-1}w$ for any $w \in W_{\bf i}$. 
\end{lemma}

{\bf Proof}. By Theorem \ref{12.13.06.10} we have   
$\widehat A{\bf K}_{\bf c^o}w = {\bf K}_{\bf c^o}\widehat Aw$. Thus 
$$
{\bf K}_{\bf c^o}\widehat A^{-1}\widehat Aw = {\bf K}_{\bf c^o}w = \widehat A^{-1}\widehat A{\bf K}_{\bf c^o}w =\widehat A^{-1}{\bf K}_{\bf c^o}\widehat Aw. 
$$ 
Since 
for a monomial $A$ the operator $\widehat A$ is an 
automorphism of the space $W_{\bf i}$, the lemma follows. 
\vskip 3mm

\vskip 3mm
\begin{corollary} \label{1.03.07.3}
For any 
$A \in  {\bf D}_{\bf i}$ and $w \in W_{\bf i}$ 
one has  ${\bf K}_{\bf c^o}\widehat A w = \widehat A {\bf K}_{\bf c^o} w$. 
\end{corollary}

{\bf Proof}. Thanks to Theorem \ref{12.13.06.10}   and Lemma \ref{p1p} 
we know the claim for  monomials corresponding to 
the vectors of the cones $C_+$ and  $-C_+$. 
The semigroup generated by these cones is the group of characters of 
the torus ${\cal D}_{\bf i}$. 
Finally, if $A_1, A_2 \in {\bf D}_{\bf i}$ commute with ${\bf K}_{\bf c^o}$ on 
$W_{\bf i}$, then for any $w\in W_{\bf i}$ one has 
$
{\bf K}_{\bf c^o}\widehat A_1\widehat A_2 w = 
\widehat A_1\widehat A_2 {\bf K}_{\bf c^o}w. 
$ Indeed, the operators $\widehat A_i$ preserve $W_{\bf i}$. 
 The corollary is proved. 
\vskip 3mm

Now we can finish the proof of the theorem. 
Let 
$$
E = \{f\in L^2({\cal A}_{\bf i}^+)~| ~{\rm exp}(\sum_i n_ia_i)f 
\in L^2({\cal A}_{\bf i}^+) \mbox{ for any $n_i \geq 0$}\}.
$$ 
We define a topology in $E$ by using the seminorms related to the operators of multiplication by ${\rm exp}(\sum_i n_ia_i)$ just as for the Schwartz spaces. 
The same argument as in the proof of Proposition \ref{qp12} 
shows that that the space $W_{\bf i}$ is dense in $E$, 
see the Remark after the proof of 
the Proposition \ref{qp12}. 
 
\begin{lemma} \label{p1}
$
{\bf K}_{\bf c^o}(E) \subset E.  
$
\end{lemma}

{\bf Proof}. The elements $B_i \in {\bf D}_{\bf i}$ act as the 
operators of multiplication by ${\rm exp}(a_i)$. Thus by Corollary 
\ref{1.03.07.3} 
one has ${\bf K}_{\bf c^o}e^{na_i}w = 
e^{na_i}{\bf K}_{\bf c^o}w$ for any $n>0$, $w \in W_{\bf i}$. 
Since $W_{\bf i}$ is dense in $E$, we get the Lemma. 

\begin{lemma}\label{1.03.07.2}
${\bf K}_{\bf c^o}$ is the operator of multiplication by a function 
$F(a)$.  
\end{lemma}

{\bf Proof}. 
For a given point $a_0$,  
the value $({\bf K}_{\bf c^o}f)(a_0)$ depends only on the value 
$f(a_0)$. Indeed, for any 
$f_0 \subset E$ with $f_0(a_0) = f(a_0)$ we have 
 $f(a) = (e^a-e^{a_0})\phi(a) + f_0(a)$, where 
$
\phi= (f-f_0)/(e^a-e^{a_0})\in E. 
$ 
Thus by Corollary \ref{1.03.07.3}
$
{\bf K}_{\bf c^o}f = (e^a-e^{a_0}){\bf K}_{\bf c^o}\phi(a) + {\bf K}_{\bf c^o}f_0(a).
$ 
So ${\bf K}_{\bf c^o}f(a_0) = {\bf K}_{\bf c^o}f_0(a_0)$. Now define  $F(a_0)$ from 
${\bf K}_{\bf c^o}f_0(a_0) = F(a_0)f_0(a_0)$. The lemma is proved.

\begin{lemma} \label{const}
The function $F(a)$ is a constant.  
\end{lemma}

{\bf Proof}. 
Let us show that $F(a)$ extends to an analytic function in the product of the upper half planes  
${\rm Im}~a_i >0$. The operator $\widehat X_k$ acts as a shift by 
$2\pi i \hbar_k$ followed by multiplication by $p^*X_k$. The latter is  
an exponential in the logarithmic coordinates $a_k$. 
By Corollary \ref{1.03.07.3} it commutes 
with ${\bf K}_{\bf c^o}$. It follows that 
${\bf K}_{\bf c^o}$ commutes with the operator of shift by $2 \pi i \hbar_k$
along the variable $a_i$. 

\begin{lemma} \label{K3}
Suppose that positive powers of the shift 
operator $T_{2 \pi i \hbar}$, $\hbar>0$, 
 commute on the subspace $W\subset L^2(\R)$  
 with the operator of multiplication 
by a function $F(a)$. Then $F(a)$ is extended to 
an analytic function in the upper half space ${\rm Im}~a >0$, and is 
invariant 
by the shift by $2 \pi i \hbar$. 
\end{lemma}

{\bf Proof}. Since $T^n(Fw) \in L^2(\R)$, its Fourier transform 
${\cal F}(Fw)$ 
is in $L^2(\R)$. Thus $e^{nx}{\cal F}(Fw) \in L^2(\R)$.
Making the inverse Fourier transforms 
and using the Payley-Wiener argument we see that $Fw$ is 
analytic in the upper half plane for any $w \in W_{\bf i}$. 
The Lemma follows. 

\vskip 3mm

This implies that $F(a)$ 
is analytic in the product of the upper half planes  
${\rm Im}~a_k >0$, and is invariant by the shift by $2 \pi i \hbar_k$ 
along the variable $a_k$. 
Employing the Langlands dual 
family of the operators $\widehat X_k^{\vee}$ we conclude that $F(a)$ 
is invariant under the shift by $2\pi i$ along the variable $a_k$. 
So if $\hbar$ is irrational then 
$F$ is a constant. Since ${\bf K}$ evidently depends 
continuously on $\hbar$, $F$ is a constant for all $\hbar$. 
Proposition \ref{P1} is proved.

\vskip 3mm {\it Step 2}. 
\begin{proposition} \label{P2} One has 
\begin{equation} \label{6.14.07.1}
{\bf K}_{{\bf c^o}}\widehat X_iw = \widehat X_i{\bf K}_{{\bf c^o}}w, \qquad 
{\bf K}_{{\bf c^o}}\widehat B_iw = \widehat B_i{\bf K}_{{\bf c^o}}w, \qquad w \in W, i\in I.
\end{equation} 
\end{proposition} 

{\bf Proof}. Let us prove the first claim. 

{\it Reduction to the case when 
${\rm det}(\varepsilon_{ij})=1$}. 
Take a set $\overline I$ containing $I$
and a skewsymmetrizable function $\overline \varepsilon$ on 
$\overline I\times \overline I$ 
extending the function $\varepsilon$, 
with  ${\rm det}(\overline \varepsilon) = 1$. 
Denote by ${\bf \overline i}= 
(\overline I, \overline \varepsilon_{ij}, \overline d_i)$ the 
obtained feed.  There is a canonical projection: 
\begin{equation} \label{TCT}
\overline p: {\cal A}_{|\overline {\bf i}|}  \lra {\cal X}_{|{\bf i}|},
\quad X_i \lms \overline p^*X_i = \prod_{j\in I}A_j^{\varepsilon_{ij}}, \quad i\in I, 
\end{equation} 
as well as its $q$-deformed version $\overline p_q: 
{\cal A}_{q, |\overline {\bf i}|} \lra {\cal X}_{q, |{\bf i}|}$. 
Write $\overline p_q^*X_i = {\Bbb A}_i^+/{\Bbb A}_i^-$. 
Here ${\Bbb A}_i^\pm$ are monomials (with non-negative exponents), so 
thanks to the quantum Laurent Phenomenon theorem 
$L^\pm_i:= \overline {\bf c}^q ({\Bbb A}_i^\pm)$ are quantum Laurent polynomials. 
By the assumption of Theorem \ref{K0} we have ${\bf c}^q (X_i) = X_i$. 
Therefore $\overline {\bf c}_q (X_i) = X_i$. 
Thus 
\begin{equation} \label{TCTq}
L^+_i(L^-_i)^{-1} = \overline {\bf c}_q ({\Bbb A}_i^+({\Bbb A}_i^-)^{-1}) = 
\overline {\bf c}_q (\overline p_q^*X_i) =  
\overline p_q^*X_i = {\Bbb A}_i^+({\Bbb A}_i^-)^{-1} \quad => \quad
L^+_i = {\Bbb A}_i^+({\Bbb A}_i^-)^{-1}L^-_i.
\end{equation}
We interpret below ${\Bbb A}_i^\pm$ as monomials in the 
quantum feed torus algebra for  
${\cal D}_{\overline {\bf i}}$ via the quantum version of (\ref{1.5.07.190}). 
Thus they act on the space $W_{\overline {\bf i}}$ by the operators 
  $\widehat {\Bbb A}_i^\pm$. 
Denote by $\overline {\bf K}_{{\bf c^o} }$, or simply by $\overline {\bf K}$,  
 the intertwiner related to the 
cluster transformation ${\bf c}$ for the feed $\overline {\bf i}$. Since the operator $(\widehat {\Bbb A}_i^-)^{-1}$  is an isomorphism on 
$W_{\overline {\bf i}}$, 
we have 
\begin{equation} \label{1.4.06.1sd}
\overline  {\bf K}_{{\bf c^o} } \widehat {\Bbb A}_i^+(\widehat {\Bbb A}_i^-)^{-1}w \stackrel{\mbox{Prop}. \ref{P1}}{=} 
\widehat {L}_i^+\overline  {\bf K}_{{\bf c^o} } (\widehat {\Bbb A}_i^-)^{-1}w \stackrel{(\ref{TCTq})}{=} 
\widehat {\Bbb A}_i^+(\widehat {\Bbb A}_i^-)^{-1}\widehat L^-_i \overline  {\bf K}_{{\bf c^o} } (\widehat {\Bbb A}_i^-)^{-1}w \stackrel{\mbox{Prop}. \ref{P1}}{=} 
\end{equation}
\begin{equation} \label{1.4.06.1sdt}
\widehat {\Bbb A}_i^+(\widehat {\Bbb A}_i^-)^{-1} \overline  {\bf K}_{{\bf c^o} } \widehat {\Bbb A}_i^-(\widehat {\Bbb A}_i^-)^{-1}w =
\widehat {\Bbb A}_i^+(\widehat {\Bbb A}_i^-)^{-1} \overline  {\bf K}_{{\bf c^o} } w.
\end{equation}
In particular the right hand side is 
in $L^2({\cal A}_{\overline {\bf i}})$. 
We are going to deduce from this the left relation in (\ref{6.14.07.1}) 
by using the restriction  to the feed ${\bf i}$. 
Observe that setting $A_l =1$ for $l \in \overline I-I$, the space 
 $W_{\overline {\bf i}}$ restricts to $W_{{\bf i}}$. 
There are two issues:

\vskip 3mm
{\it Issue (i)}. Extending the feed ${\bf i}$, we change 
the action 
of a mutation $\mu_k$ on $B_j$-coordinates, $j \in I$, although we 
do not change its  action  
on the $X_j$-coordinates, $j \in I$.

{\it Issue (ii)}. Extending the feed ${\bf i}$, we change the intertwiner 
${\bf K}_{\bf i}$. 

\vskip 3mm
To handle them, we employ the following facts:
\vskip 3mm
{\it Fact (i)}. The ${\cal D}$-cluster transformation 
corresponding to the mutation $\overline \mu_k$ of the 
feed $\overline {\bf i}$ in the direction $k \in I$
preserves the  functions 
$B_l$, $l\in \overline I-I$. So we can restrict it  
to the subvariety of the double given by the equations 
$B_l=1$, $l\in \overline I-I$. Then restriction of the function 
$B_i$, $i \in I$,  to this subvariety mutates the same 
way as under  
the ${\cal D}$-cluster transformation corresponding to $\mu_k$.

{\it Fact (ii)}. The natural embedding 
${\cal A}_{\bf i}^+ \hra {\cal A}_{\overline {\bf i}}^+$ gives rise to 
an injective  map of linear spaces
\begin{equation} \label{1.4.06.1}
i_{{\bf i}\to \overline {\bf i}}: 
{\cal L}^2({\cal A}_{\bf i}^+) \lra {\rm Dist}
({\cal A}_{\overline {\bf i}}^+), \qquad 
f \lms f\prod_{l \in \overline I-I}\delta(a_l),
\end{equation}
where on the right stands  the dual 
 to the classical Schwarz space on ${\cal A}_{\overline {\bf i}}^+$. 
The intertwiner $\overline {\bf K}$ provides a map of 
the classical Schwarz spaces: indeed, it is a composition of Fourier transforms and 
operators of multiplication by a 
smooth function of absolute value $1$.  Therefore 
it gives rise to a map of the distribution spaces. 
We claim that it restricts to an operator between the images of map 
(\ref{1.4.06.1}), which coincides with the  
intertwiner ${\bf K}$. In other words, there is a commutative diagram
$$
\begin{array}{ccc}
{\cal L}^2({\cal A}_{\bf i'}^+) & \lra &{\rm Dist}
({\cal A}_{\overline {\bf i'}}^+)\\
%&&\\
{\bf K}\downarrow && \overline {\bf K}_{}\downarrow \\
%&&\\
{\cal L}^2({\cal A}_{\bf i}^+) & \lra &{\rm Dist}
({\cal A}_{\overline {\bf i}}^+)
\end{array}
$$

Fact (i) is obvious. Fact (ii) is clear from the definition of the intertwiner:
Restricting the differential operator $\widehat x_k$, $k \in I$, 
 acting on ${\cal A}^+_{\overline {\bf i}}$ to 
the distributions annihilated by operators of multiplication by $a_l$, $l \in \overline I - I$, we get  
the operators $\widehat x_k$ acting on ${\cal A}^+_{{\bf i}}$. Thus 
the operator $\overline {\bf K}^\sharp$ restricts to  ${\bf K}^\sharp$. 
The operator $\overline {\bf K}'$ obviously restricts to ${\bf K}'$. 
Thus (\ref{1.4.06.1sd})-(\ref{1.4.06.1sdt})
 imply the left identity in (\ref{6.14.07.1}).

The second claim in (\ref{6.14.07.1}) is proved similarly. 
Although $\overline {\bf c}_q^d(B_i)$ may differ from $B_i$, 
we run the same argument as in (\ref{1.4.06.1sd})-(\ref{1.4.06.1sdt}) with the following modification: 
we apply the operators to the distribution $i_{{\bf i} \to \overline {\bf i}}(w) \in 
{\rm Dist}
({\cal A}_{\overline {\bf i}}^+)$ corresponding to 
$w \in W_{\bf i}$, and use the fact that 
the restriction of $\overline {\bf c}_q^d(B_i)$ 
to the subspace given by the equations $B_l=1, l\in \overline I-I$ 
coincides with  $B_i$. 
Proposition \ref{P2} is proved. 
\vskip 3mm

Theorem \ref{K0}(ii) follows from 
Proposition \ref{P2}. Indeed, since ${\bf K}_{\bf c^o}$ commutes with the operators $\widehat B_i$, $i \in I$, 
we conclude, just like in the proof of Lemma \ref{1.03.07.2}, that ${\bf K}_{\bf c^o}$ 
is an operator of multiplication by a function. Then, since  ${\bf K}_{\bf c^o}$ 
commutes with the operators $\widehat X_i$ and $\widehat X^{\vee}_i$, arguing as  in 
the end of the proof of Proposition \ref{P1} we conclude that this function is a constant. 
Theorem \ref{K0} is proved.

\vskip 3mm
{\bf A result providing assumptions 
of Theorem \ref{K0}(ii)}. 
Let $\widetilde {\bf i}:= 
(I, I_0, \widetilde \varepsilon_{ij}, \widetilde d_i)$ be a feed. 
Restricting to $I-I_0$, we get a feed 
${\bf i}:= (I-I_0, \varepsilon_{ij}, d_i)$. There is a canonical projection, 
given as a composition
\begin{equation} \label{TCT}
{\cal A}_{|\widetilde {\bf i}|} 
\lra {\cal X}_{|\widetilde {\bf i}|} \lra {\cal X}_{|{\bf i}|},
\quad X_i \lms \prod_{j\in I}A_j^{\varepsilon_{ij}}, i\in I-I_0. 
\end{equation}
The precursor of the following result is Lemma \ref{8.24.05.1}. 

\begin{proposition} \label{6.11.07.1} Assume that a feed 
cluster transformation ${\bf c}: {\bf i} \to{\bf i}$, considered as a 
cluster  transformation 
$\widetilde {\bf c}: \widetilde{\bf i} \to \widetilde{\bf i}$, is  trivial.  
Assume that the composition (\ref{TCT}) is surjective. 
Then the induced by ${\bf c}$ cluster maps ${\bf c}^x_q$  and 
${\bf c}^d_q$ of the $q$-deformed  
spaces ${\cal X}_{q, |{\bf i}|}$ and ${\cal D}_{q, |{\bf i}|}$ 
 are trivial. 
\end{proposition}

{\bf Proof}. Since $\widetilde {\bf c}$ is trivial, by \cite{BZq}, the cluster transformation ${\bf c}^a_q$ 
of the $q$-deformed cluster algebra related to the feed 
$\widetilde{\bf i}$ is trivial.  
Recall the surjective monomial map 
$
\widetilde p^q: 
{\cal A}_{q, |\widetilde {\bf i}|}
\stackrel{}{\lra}  {\cal X}_{q, |{\bf i}|}. 
$ 
It commutes with mutations (Appendix to Section 2). 
This implies the proposition for the ${\cal X}$-space. 

{\it The ${\cal D}$-space}. A feed tori ${\cal D}_{\widetilde {\bf i}}$ 
has cluster coordinates 
$(B_i, X_i)$, $i \in I$. Consider 
a torus  ${\cal D}^\times_{\widetilde {\bf i}}$
with coordinates $(B_i, X_j)$, where $i \in I, j \in I-I_0$.
There is a canonical projection 
$\pi_{\widetilde {\bf i}}: {\cal D}_{\widetilde {\bf i}} \lra 
{\cal D}^\times_{\widetilde {\bf i}}$, given by dropping 
the ${\Bbb G}_m$-components of 
${\cal D}_{\widetilde {\bf i}}$ with the coordinates $X_j$, $j \in I-I_0$. 
Mutations in the directions of the set $I-I_0$ 
of the ${\cal X}$-coordinates of a 
${\cal D}$-space parametrised by $I-I_0$
do not involve the ${\cal X}$-coordinates parametrised by $I_0$. 
Thus, using the same formulas as in the definition of the 
${\cal D}$-space, we define a new positive space ${\cal D}^\times_{|\widetilde {\bf i}|}$. 
There is a canonical projection 
$
\pi: {\cal D}_{|\widetilde {\bf i}|} \lra 
{\cal D}^{\times}_{|\widetilde {\bf i}|}, 
$ 
and a unique Poisson structure on 
${\cal D}^{\times}_{|\widetilde {\bf i}|}$ for which the map $\pi$ 
is Poisson. It is 
given by the same formulas 
as for ${\cal D}_{|\widetilde {\bf i}|}$. 
Similarly 
there is a $q$-deformation ${\cal D}^{\times}_{q, |\widetilde {\bf i}|}$ 
of the space 
${\cal D}^\times_{|\widetilde {\bf i}|}$, as well as of the projection $\pi$. 

There is a projection  
\begin{equation} \label{6.13.07.3}
\widetilde \varphi: 
{\cal A}_{|\widetilde {\bf i}|} \times {\cal A}^o_{|\widetilde {\bf i}|} 
\lra {\cal D}^\times_{|\widetilde {\bf i}|}, \qquad \widetilde \varphi:= 
\pi\circ \varphi.
\end{equation}

\begin{lemma} \label{6.13.07.2} 
Assume that the rectangular matrix 
$\varepsilon_{ij}$, where $i \in I, j\in I-I_0$, is non-degenerate. 
For any $q$-deformation $ {\cal A}_{q, |\widetilde {\bf i}|}$ 
of the $ {\cal A}$-space as in \cite{BZq}, there is a map of quantum spaces
$$
\widetilde \varphi^q: 
{\cal A}_{q, |\widetilde {\bf i}|} \times {\cal A}^o_{q, |\widetilde {\bf i}|} 
\stackrel{}{\lra}   
{\cal D}^\times_{q, |\widetilde {\bf i}|}.
$$ 
\end{lemma} 

{\bf Proof}. The map 
$\widetilde\varphi^q$ is a monomial map, which 
for the feed tori is determined by its $q=1$ counterpart 
$\widetilde\varphi_{\bf i}^q$. 
One checks that 
it commutes with the mutations, see also \cite{King}. 

\vskip 3mm

\begin{lemma} \label{6.13.07.1} If the matrix 
$\varepsilon_{ij}$, $i \in I, j\in I-I_0$, is non-degenerate, then 
the map (\ref{6.13.07.3}) is surjective. 
\end{lemma}

{\bf Proof}. Assume that $\widetilde \varphi^* 
\prod_{i, j} X_i^{x_i}B_j^{b_i} = 1$. Then $b_j=0$ since 
$\widetilde \varphi^*B_j = A^o_j/A_j$, and $\widetilde \varphi^*X_j$ 
does not involve $A_j^o$'s. Thus the 
 assumption of the lemma implies $x_i=0$. 
The lemma is proved. 
\vskip 3mm

There is an embedding ${\cal D}_{q, |{\bf i}|} 
\hra {\cal D}^\times_{q, |\widetilde {\bf i}|}$, given in any cluster coordinate system by setting $B_i=1$ for $i \in I-I_0$. 
Now $\widetilde {\bf c}^a = {\rm Id}$ implies 
$\widetilde {\bf c}_q^a = {\rm Id}$ by \cite{BZq}, the latter implies 
by Lemma \ref{6.13.07.1} 
that the cluster transformation 
$\widetilde {\bf c}_q^d$ of the space $
{\cal D}^\times_{q, |\widetilde {\bf i}|}$ is the identity. 
Thus the restriction $\widetilde {\bf c}_q^d$ to ${\cal D}_{q, |{\bf i}|}$ 
is  the identity map.  The proposition is proved. 
\vskip 3mm

\subsection{A proof of Theorem \ref{MT}}
{\it Canonical representation}. Let $\Lambda$ be a lattice with a 
symplectic form $(\ast, \ast)$ with values 
in $\frac{1}{N}\Z$. Denote by ${\rm T}_{\Lambda}$ 
the corresponding quantum torus $\ast$-algebra 
with the parameter $q = e^{\pi i \hbar}$, 
and by ${\rm T}_{\Lambda}^{\vee}$ the algebra with the dual 
parameter $q^{\vee} = e^{\pi i /\hbar}$. So ${\rm T}_{\Lambda} \otimes 
{\rm T}_{\Lambda}^{\vee}$ is the modular double. 
Set $\Lambda_\Q:= \Lambda\otimes \Q$. Choose  a decomposition 
$\Lambda_{\Q} = P \oplus Q$ into a direct sum 
of two Lagrangian subspaces. Set $P_\R:= P\otimes \R$. 
One defines, similarly to (\ref{9.4.07.1}), a space 
$W(P_\R)$ with a canonical representation 
of the $\ast$-algebra ${\rm T}_{\Lambda}\otimes {\rm T}_{\Lambda}^{\vee}$: 
elements $p\in P_{\Q}$ and $q\in Q_{\Q}$
provide operators 
$$
\widehat p \lms T_{2\pi i \hbar p}, \quad \widehat q \lms e^{(q, \ast)}, \quad 
\widehat p^{\vee} \lms T_{2\pi i p}, \quad \widehat q^{\vee} \lms e^{(q/\hbar , \ast)}.
$$ 
Here $T_{a}f(p):= f(p-a)$, and  $e^{(q, \ast)}$ is 
the operator of multiplication by the function $e^{(q, p)}$. 
The Schwartz space ${\cal S}_{P,Q}$ 
is the maximal domain 
of the $\ast$-algebra ${\rm T}_{\Lambda}\otimes {\rm T}_{\Lambda}^{\vee}$ 
in $L^2(P_\R)$. 
It has a natural Freschet topology;  
$W(P_\R)$ is dense in this topology -- 
the proof is the same as of Proposition \ref{qp12}). 

It is easy to see that, given another 
decomposition $\Lambda_\Q:= P' \oplus Q'$, there is a 
canonical up to a constant unitary operator 
${\cal S}_{P,Q} \to {\cal S}_{P',Q'}$ intertwining 
the action of the  algebra ${\rm T}_{\Lambda}\otimes {\rm T}_{\Lambda}^{\vee}$. 
Its kernel 
is the exponential of an imaginary
 quadratic expression of the coordinates. 
If $\hbar \not \in \Q$, it is determined uniquely 
by the difference equations 
expressing the commutation relations.
Therefore the representation of the $\ast$-algebra 
${\rm T}_{\Lambda}\otimes {\rm T}_{\Lambda}^{\vee}$ is defined 
up to a canonical modulo constant unitary isomorphism. 
Having in mind this, 
we can say that there is 
a canonical representation $(V, {\cal S}_V)$ of 
${\rm T}_{\Lambda}\otimes {\rm T}_{\Lambda}^{\vee}$, where 
$V$ is a Hilbert space, and 
${\cal S}_V$ its subspace where the modular 
double ${\rm T}_{\Lambda}\otimes {\rm T}_{\Lambda}^{\vee}$ 
acts. 

If the form $(\ast, \ast)$ on $\Lambda$ is degenerate, denote by 
$Z \subset \Lambda$ its kernel. A model for the canonical 
representation is defined then by choosing a pair of subspaces 
$P, Q$ containing $Z$, which project to a pair of 
transversal Lagrangian subspaces in $\Lambda_\Q/Z_\Q$. 
One can decompose it further via the characters $\lambda$ of $Z_\R$ 
using the Fourier transform, getting a family of representations 
$V_\lambda$ and a canonical 
decomposition into an integral of Hilbert spaces $V = \int V_\lambda d\lambda$.  
 
\vskip 3mm
The representation of ${\rm T}_{{\Lambda}_{\cal D}}\otimes 
{\rm T}_{{\Lambda}_{\cal D}}^{\vee}$ in $L^2({\cal A}^+)$ 
from Section 5.1 is the canonical representation. 
%For technical reasons it is handy to pick 
%its realization given by the operators (\ref{anotherR}). 
Take the sublattice $\pi^*(\Lambda_{\cal X} \oplus 
\Lambda^{\rm opp}_{\cal X})$ of  $\Lambda_{\cal D}$. 
The following lemma is an easy linear algebra exercise. 

\begin{lemma} \label{PR} The canonical representation of ${\rm T}_{\pi^*(\Lambda_{\cal X} \oplus \Lambda^{\rm op}_{\cal X})}\otimes 
{\rm T}^{\vee}_{\pi^*(\Lambda_{\cal X} \oplus \Lambda^{\rm op}_{\cal X})}$ 
is identified with restriction of the canonical representation 
of ${\rm T}_{\Lambda_{\cal D}}\otimes 
{\rm T}^{\vee}_{\Lambda_{\cal D}}$. 
\end{lemma}

The group $H_{\cal A}(\R_{>0})$ acts on  ${\cal A}^+$ 
and hence on $L^2({\cal A}^+)$. On the other hand, let $Z$ be the kernel 
of the restriction of the symplectic form $(\ast, \ast)_{\cal D}$ to 
the lattice $\pi^*(\Lambda_{\cal X} \oplus \Lambda^{\rm opp}_{\cal X})\otimes \Q$. 
In realization (\ref{anotherR}) it gives rise to 
an action of the real torus with the 
group of cocharacters $Z \otimes \R$.  
Indeed, given a cluster basis $\{e_i\}$ of $\Lambda_{\cal X}$, 
 the sublattice  $Z$ consists of the vectors 
$(\sum c_ie_i, \sum c_ie_i)$ where 
$\sum c_i \varepsilon_{ij} =0$ for every $j \in I$. So in realization 
(\ref{anotherR}) the corresponding operator in the representation is just 
$\sum_{i} c_i \partial/\partial a_i$ --
 the linear term is gone thanks to the condition on $c_i$'s. 
These two real tori are canonically isomorphic, and the 
two actions coincide.  

Therefore the Fourier decomposition of $L^2({\cal A}^+)$ 
according to the characters of $H_{\cal A}(\R_{>0})$ is the same thing as 
the Fourier decomposition according to the center 
of ${\rm T}_{\pi^*(\Lambda_{\cal X} \oplus \Lambda^{\rm op}_{\cal X})}\otimes 
{\rm T}^{\vee}_{\pi^*(\Lambda_{\cal X} \oplus \Lambda^{\rm op}_{\cal X})}$. 
Thanks to theorem \ref{K0} 
 these two natural decompositions are respected by the intertwiners. 
Finally, 
a module over ${\rm T}_{\pi^*(\Lambda_{\cal X} \oplus \Lambda^{\rm op}_{\cal X})}\otimes 
{\rm T}^{\vee}_{\pi^*(\Lambda_{\cal X} \oplus \Lambda^{\rm op}_{\cal X})}$ with the central character 
$(\lambda, \lambda)$ is identified with the  representation ${\rm End}V_{\lambda}$. The theorem is proved.

\vskip 3mm
\subsection{The intertwiner and the quantum double}

Our goal is to show how the intertwiner  ${\bf K}_{{\bf i'} , {\bf i}}$ determines the 
coordinate transformations which we used to define the quantum double. 
The idea is similar to the ${\cal D}$-module approach to integral geometry 
from \cite{Go1}, used here for the difference equations.  
\vskip 3mm
Recall the (non-commutative) algebra ${\bf D}_{\bf i}$ of regular functions 
on the quantum double feed torus  
for the feed  ${\bf i}$. 
We realize it 
as the algebra of all $q$-difference operators on the 
torus ${\cal A}_{\bf i}$, see (\ref{intertwiner1}). Its 
fraction field ${\Bbb D}_{\bf i}$ is a non-commutative field. 

The ring 
$
{\bf D}_{\bf i} \otimes {\bf D}_{\bf i'}
$ 
as a bimodule:  a left ${\bf D}_{\bf i}$-module and a
right  ${\bf D}_{\bf i'}$-module. Denote by  
$
{\cal R}_{{\bf i}, {\bf i}'}\subset {\bf D}_{\bf i} \otimes {\bf D}^{\rm op}_{\bf i'}
$ 
be the ideal for this bimodule structure generated by all $q$-difference equations satisfied by the generalized function 
$K_{{\bf i'} , {\bf i}}(a,a')$ representing the kernel of the intertwiner 
${\bf K}_{{\bf i'} , {\bf i}}$. 
Consider the corresponding bimodule
$$
M_{{\bf i'} , {\bf i}}:= 
\frac{{\bf D}_{\bf i} \otimes {\bf D}^{\rm op}_{\bf i'}}{{\cal R}_{{\bf i'}, {\bf i}}}.
$$
Denote by ${\Bbb M}_{{\bf i} , {\bf i}'}$ its localization. 
It is a module over the tensor product 
${\Bbb D}_{\bf i} \otimes {\Bbb D}_{\bf i'}^{\rm op}$.

\begin{proposition} \label{11.14.06.1}
The restriction of the bimodule 
${\Bbb M}_{{\bf i'} , {\bf i}}$ to each of the factors in ${\Bbb D}_{\bf i} \otimes {\Bbb D}_{\bf i'}^{\rm op}$  
is a free rank one module 
with a canonical generator ${\Bbb K} = {\Bbb K}_{{\bf i'}, {\bf i}}$. 
\end{proposition}

{\bf Proof}. This is equivalent to the similar property for the kernel $K^{\sharp}$. 

\vskip 3mm
\begin{corollary} \label{11.14.06.3}
There is a homomorphism of fields 
$
{\kappa}: {\Bbb D}_{\bf i'} \lra {\Bbb D}_{\bf i}, \quad 
A{\Bbb K}= {\Bbb K}{\kappa}(A).
$ 
\end{corollary} 
This non-commutative 
field homomorphism is the coordinate transformation 
employed in the definition of the quantum double.

\section{Quantization of higher Teichm\"uller spaces}

\subsection{Principal series representations of 
quantum moduli spaces ${\cal X}_{G, S}$}

In this Section we denote by $S$ an oriented
 surface with holes and a finite collection of marked points on the boundary, considered modulo isotopy. We reserve the notation $\widetilde S$ to the 
surface with the marked points omitted. 
This differs from the notation we used in \cite{FG1} and the Introduction, 
where the surface with marked points was denoted by $\widehat S$.

\begin{figure}[ht]
\centerline{\epsfbox{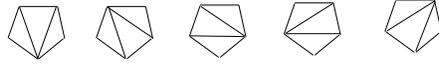}}
\caption{The composition of five flips is the identity.}
\label{double15}
\end{figure}

In Chapter 10 of \cite{FG1} we 
defined a cluster ensemble structure for the pair of moduli spaces 
$({\cal A}_{G, S}, {\cal X}_{G,S})$ where  $G=SL_m$. 
Namely,  let us shrink all holes without marked points at the boundary 
to puncutres. Let us pick a point, called a special point, in each 
 connected component of 
$$
\mbox{
$\{$the boundary of $\widetilde S$ (minus punctures) $\}$ -  $\{$the marked points$\}$).}
$$   
An ideal triangulation $T$ of 
$S$ ({\it loc. cit.}) is a triangulation of $\widetilde S$ 
with vertices at the special points and punctures. An ideal triangulation $T$ 
produces a feed 
${\bf i}_T$. To each edge $E$ of $T$ we assign a flip of the triangulation, 
producing a new ideal triangulation $T'$. The feed 
 ${\bf i}_{T'}$ is obtained from ${\bf i}_T$ by a sequence of 
$(m-1)^2$ mutations, described in {\it loc. cit.}. 
Their composition  is a feed cluster transformation $c_E: 
{\bf i}_T \to {\bf i}_{T'}$. The classical 
modular groupoid ${\Bbb G}_S$ of $S$ 
is a groupoid whose objects are ideal triangulations of $S$ 
up to an isomorphism. Its 
morphisms are generated by the flips. 
The relations are generated by ``rectangles'' and ``pentagons''. 
The former relation tells that flips at disjoint edges commute. 
The pentagon relation tells that the composition of five flips 
of the diagonal of an ideal pentagon, presented on Fig. \ref{double15}, 
is an automorphism. 

The cluster structure of the pair of moduli spaces 
$({\cal A}_{G,   S}, {\cal X}_{G,  S})$ provides 
a saturated cluster modular group $\widehat \Gamma_{G,  S}$, and 
a cluster modular group $\Gamma_{G,  S}$. 
Let $\widehat \Gamma_{G,  S} \to \Gamma_{G,  S}$ be the canonical 
epimorphism, 
We conjecture that it is an isomorphism.  Since 
a flip is decomposed into a composition of mutations, 
there is a canonical embedding 
$\Gamma_{S} \hra \Gamma_{G,  S}$ ({\cal loc. cit.}).  
The group $\widehat \Gamma_{S}$ is the preimage of $\Gamma_{S}$ in 
$\widehat \Gamma_{G,  S}$. There are similar maps of modular groupoids. 
So we are getting the diagrams
$$
\begin{array}{ccc}
\widehat \Gamma_{S}& \hra &\widehat \Gamma_{G,  S}\\
\downarrow &&\downarrow \\
\Gamma_{S}& \hra &\Gamma_{G,  S}
\end{array}; \qquad \qquad 
\begin{array}{ccc}
\widehat {\cal G}_{S}& \hra &\widehat {\cal G}_{G,  S}\\
\downarrow &&\downarrow \\
{\cal G}_{S}& \hra &{\cal G}_{G,  S}
\end{array}
$$

Applying the main construction of 
this paper to the cluster ensemble assigned to the pair 
$({\cal A}_{G,  S}, {\cal X}_{G,  S})$, we get a unitary projective 
representation of the modular 
groupoid $\widehat {\cal G}_{G,  S}$, and thus a representation  
of the groupoid $\widehat {\cal G}_{S}$. 
Our goal is to show that the latter descends to 
representations of the classical modular groupoid ${\cal G}_{S}$. 
This is equivalent to the part (iii) of following:

\begin{theorem} \label{MTH3} Let $S$ be an open surface, and $G=SL_m$. 
Then 

(i) The   
modular group $\Gamma_S$ acts by automorphisms 
of the 
$q$-deformed spaces  ${\cal X}^q_{G,  S}$ and ${\cal D}^q_{G,  S}$. 

(ii) The composition of the intertwiners 
corresponding to the pentagon relation in the classical 
modular groupoid is a multiple of the identity transformation. 

(iii) The group $\Gamma_S$ acts in the canonical 
unitary projective representation 
related to the pair of moduli spaces $({\cal A}_{G,  S}, {\cal X}_{G,  S})$. 
 \end{theorem}

{\bf Proof}.  Theorem \ref{K0} plus (i) implies (ii). 
Clearly (ii) implies (iii).

(i) Recall the moduli space ${\rm Conf}_m{\cal B}$ 
(respectively ${\rm Conf}_m{\cal A}$)
of configurations of $m$ flags (respectively affine flags) 
in $G$. There is a canonical 
surjective projection
\begin{equation} \label{PR}
p: {\rm Conf}_m{\cal A} \lra {\rm Conf}_m{\cal B}.
\end{equation}

We have to prove the quantum pentagon 
relation, asserting that 
the cluster 
transformation ${\bf c_P}$, corresponding to a pentagon $P$ of an ideal triangulation 
of $S$ is the identity map for the quantum spaces ${\cal X}^q_{G,S}$ and ${\cal D}^q_{G,S}$. The map ${\bf c_P}$ is the composition 
of quantum mutations corresponding to the sequence 
of $5$ flips of the  pentagon $P$. 
Since we work with a single ideal pentagon, 
it is enough to handle the case of the  
quantum configuration space ${\rm Conf}_5^q({\cal B})$ and 
the corresponding ${\cal D}$-space.

\begin{figure}[ht]
\centerline{\epsfbox{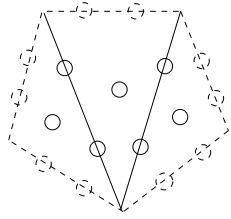}}
\caption{A feed for a configuration of $5$ affine flags in $SL_3$.}
\label{double10}
\end{figure}

Let $\widetilde {\bf i}:= (I, I_0, \varepsilon_{ij})$ be a feed 
corresponding to a triangulated pentagon which 
describes the cluster structure of the moduli space ${\rm Conf}_5({\cal A})$. 
The set $I$ (respectively the frozen set $I_0$) corresponds to the 
vertices of the $m$-triangulation assigned to 
the vertices (respectively external vertices) of the pentagon. The case $m=3$ 
is presented on Fig. \ref{double10}. 
Denote by ${\bf i}$ the restriction of the feed 
$\widetilde {\bf i}$ to $I-I_0$.
Then the configuration space ${\rm Conf}_5{\cal B}$ 
(respectively ${\rm Conf}_5{\cal A}$)
has a cluster ${\cal X}$- (respectively ${\cal A}$-) variety  
structure described by the feed  
${\bf i}$ (respectively $\widetilde {\bf i}$). The cluster structure 
of projection 
(\ref{PR}) is described as a composition (\ref{TCT}). 

The feed cluster transformation ${\bf c}_P: {\bf i} \to {\bf i}$ 
satisfies  conditions of Proposition
\ref{6.11.07.1}.  Indeed,  the corresponding 
cluster transformation $\widetilde {\bf c}_P^a$ of the cluster variety 
 ${\cal A}_{|\widetilde {\bf i}|}$ is trivial, 
since it is obviously so 
 for the corresponding moduli space ${\rm Conf}_5{\cal A}$, and 
projection (\ref{PR}) is surjective. 
The theorem is proved. 

%\vskip 3mm

\subsection{The modular functor conjecture} 
Let us cut a surface $S$ along a simple loop $\gamma$, getting a surface $S'$. 
It has two new boundary components $\gamma_{\pm}$. So a principal series reprsentation 
assigned to $S$ is parametrised by a central character $\lambda$, a one 
assigned to $S'$ is parametrised by a triple $(\lambda, \chi_-, \chi_+)$ where 
$\chi_\pm$ are character of the torus $H(\R_{>0})$. The Weyl group $W$ acts 
on the latter characters, and the dominant ones form a set of 
representatives.

\begin{conjecture} \label{MK} Given a surface $S$ with holes, a principal series 
representation ${\bf V}_{G,   S; \lambda}$ 
 with a central character $\lambda$ 
has a natural $\Gamma_{G,   S'}$-equivariant 
decomposition into an integral of Hilbert spaces
\begin{equation} \label{MK1}
{\bf V}_{G,   S; \lambda} \stackrel{\sim}{\lra} 
\int_{\chi}{\bf V}_{G,   S'; \lambda, \chi, \chi^{-1}} d\mu_\chi
\end{equation} 
where ${\bf V}_{G,   S'; \lambda, \chi, \chi^{-1}}$ is the principal series representation 
assigned to the surface $\widehat S'$ and the 
central character $(\lambda, \chi, \chi^{-1})$, where $\chi$ is a 
dominant character of the group $H(\R_{>0})$.
\end{conjecture}

%Any two (infinite-dimensional) 
%Hilbert spaces are isomorphic. So the conjecture is meaningless unless we 
Let us put an extra data on the  Hilbert spaces which rigidifies isomorphism (\ref{MK1}). 
Denote by ${\bf L}_S$ and ${\bf L}_{S'}$ the ${\bf L}$-algebras  
for the surfaces $S$ and $S'$. For every boundary component $\alpha$ of $S$, the algebra 
${\bf L}_S$ has a central subalgebra ${\bf L}_\alpha$ 
identified with the algebra of regular functions on the Cartan group $H$. 
Thus the algebra ${\bf L}_{S'}$ contains a subalgebra 
$\widetilde {\bf L}_{S'}$
given by the condition ``the monodromies around the boundary components $\gamma_-$ and 
$\gamma_+$ are opposite to each other'', i.e. by the (anti)diagonal subalgebra 
in ${\bf L}_{\gamma_-} \otimes {\bf L}_{\gamma_+}$. 
There should be (see the comments below) a $\Gamma_{S'}$-equivariant embedding 
$$
i_{S',S}: \widetilde {\bf L}_{S'} \hra {\bf L}_{S}. 
$$
The word ``natural'' in the formulation of Conjecture \ref{MK} 
means that

\begin{itemize}
%\item
%(i) The isomorphism (\ref{MK1}) is $\Gamma_{S'}$-equivariant, or even better,
% $\Gamma_{G, S'}$-equivariant. 
\item
{\it For (almost every) character $\chi$,  map (\ref{MK1}) 
leads to a map of the corresponding Schwartz spaces 
$$
\alpha_\chi: {\cal S}({\bf V}_S) \lra {\cal S}({\bf V}_{S', \chi}), 
$$
commuting with the action of the algebra  $\widetilde {\bf L}_{S'}$ -- 
the latter acts on the left via the map $i_{S',S}$}. 
\end{itemize}

{\it Comments}. 
The algebra generated by the monodromies $M_{\alpha}$ 
along  loops $\alpha$ on $S$ generate the algebra of regular functions on 
the moduli space ${\cal L}_{G,S}$ of $G$-local systems on $S$. 
The monodromy lies in $G/{\rm Ad}G$, so
any $W$-invariant function on the Cartan group gives 
rise to a regular function on ${\cal L}_{G,S}$. 

It was proved in \cite{FG1} that the monodromy $M_{\alpha}$ provides regular functions 
on the cluster ${\cal X}$-variety ${\cal X}_{G,S}$. 
Let us assume a similar claim for the $q$-deformed space. Then the algebra 
${\bf L}_{S}$ is of the ``right size'', and evidently containes 
$\widetilde {\bf L}_{S'}$ as a subalgebra: the latter is generated by the 
quantum monodromies  
along the loops  on $S$ which do not intersect the loop $\gamma$. 

If two loops on $S$ do not intersect, 
the corresponding monodromy functions on ${\cal L}_{G,S}$ commute for the 
Poisson bracket. Let us assume the corresponding quantum statement. 
 Then the subalgebra $\widetilde {\bf L}_{S'}$ centralizes 
the one ${\bf L}_{\gamma}$ generated by the monodromy along the loop $\gamma$. 
Thus the isomorphism (\ref{MK1}) is the spectral decomposition 
for the algebra ${\bf L}_{\gamma}$, understood as a $\Gamma_{S'}$-equivariant 
$\widetilde {\bf L}_{S'}$-module. 

\vskip 3mm
Finally, to justify the 
name ``modular functor'' one should add to Conjecture \ref{MK} 
relations between the spaces ${\bf V}_{G,S, \lambda}$ associated to 
surfaces $S$ with varying 
number of holes which correspond geometrically to closing a hole (i.e. filling the 
hole  by a disc by 
gluing a disc to the boundary component of $S$ corresponding to the hole), 
sometimes called `
`propagation of vacua'' in the context of 
the conformal field theory. 

\vskip 3mm
\paragraph{Representation of the modular groupoid as a combinatorial connection.} 
A rational conformal field theory provides  a bundle with 
a flat connections on the 
moduli space ${\cal M}_{g,n}$ -- 
the bundle of holomorphic conformal blocks. 

In our case we get a projective unitary representation of the cluster 
modular groupoid ${\Bbb G}_{G, S}$. 
Let us discuss first the case of $G = PGL_2$. 
Then the vertices of the modular groupoid ${\Bbb G}_{S}$ are 
ideal triangulations $T$ of $S$. Our construction assigns 
to $T$ a Hilbert 
vector space $V_T$.  The edges of ${\Bbb G}_{S}$ 
correspond to flips $T\to T'$. We assign to them 
unitary maps $V_T\to V_{T'}$, and  interpret this as a connection 
on the polyhedral complex ${\Bbb G}_{S}$. 
The pentagon equation just means that this connection 
is flat. 
In general we get a flat connection  on the polyhedral complex ${\Bbb G}_{G, S}$. 

The universal cover of the polyhedral complex ${\Bbb G}_{G, S}$ can be viewed 
as sitting at the projectivisation of the Thurston boundary 
of the higher Teichmuller space ${\cal X}_{G,S}(\R_+)$. 
Indeed, assume first that $G=PGL_2$. Then an ideal  triangulation $T$ of $S$ 
provides an integral ${\cal X}$-lamination 
in the terminology of \cite{FG1}, Section 12, by assigning the weight $1$ 
to every edge of $T$. Flips provide pathes connecting these boundary points, etc. 
In general a feed ${\bf i}$ 
provides a point ${\cal X}_{G,S}(\R^t_+)$ with all 
the tropical coordinates equal to $1$. 

It would be interesting to get  
this connection as a limit of a flat $\Gamma_{G,S}$-equivariant 
connection on 
the higher Teichmuller space ${\cal X}_{G,S}(\R_+)$.

\vskip 3mm

V.F.: ITEP, B. Cheremushkinskaya 25, 117259 Moscow, Russia.
fock@math.brown.edu

A.G.: Department of Math, Brown University, 
Providence RI 02906, USA; sasha@math.brown.edu

\end{document}